\documentclass[reqno]{amsart}
\usepackage{amsfonts}
\usepackage{amsmath}
\usepackage[english]{babel}
\usepackage{amsthm}
\usepackage{inputenc}
\usepackage{amssymb}
\usepackage{graphicx}
\usepackage{float}

\newtheorem{thm}{Theorem}[section]
\newtheorem{cor}[thm]{Corollary}
\newtheorem{lem}[thm]{Lemma}

\newtheorem{defn}[thm]{Definition}
\newtheorem{exam}[thm]{Example}
\newtheorem{rem}[thm]{Remark}
\numberwithin{equation}{section}

\input{tcilatex}

\begin{document}
\title[Nilpotent Jordan Algebras]{\textbf{Classification of five-dimensional
nilpotent Jordan algebras}}
\author{A. S. Hegazi and Hani Abdelwahab}
\address[A. S. Hegazi -- Hani Abdelwahab]{Department of Mathematics, Faculty
of science, Mansoura university. (Egypt)}
\email{hegazi@mans.edu.eg --- haniamar1985@gmail.com}

\begin{abstract}
The paper is devoted to classify nilpotent Jordan algebras of dimension up
to five over an algebraically closed field $\mathbb{F}$ of characteristic
not $2$. We obtained a list of $35$ isolated non-isomorphic $5$-dimensional
nilpotent non-associative Jordan algebras and $6$ families of non-isomorphic 
$5$-dimensional nilpotent non-associative Jordan algebras depending either
on one or two parameters over an algebraically closed field of
characteristic $\neq 2,3$. In addition to these algebras we obtained two
non-isomorphic $5$-dimensional nilpotent non-associative Jordan algebras
over an algebraically closed field of characteristic $3$.
\end{abstract}

\maketitle

\section{Introduction}

The classification of nilpotent Jordan algebras is one of the fundamental
problems. In this way nilpotent Jordan algebras of small dimension were
classified. Nilpotent associative Jordan algebras of dimension up to five
over algebraically closed fields of characteristic not $2,3$ were classified
in \cite{Mazzola}. Nilpotent associative Jordan algebras of dimension up to
five over algebraically closed fields were classified in \cite{Ponnen}.
Nilpotent Jordan algebras of dimension up to four over $%
\mathbb{C}
$ were classified in \cite{Ancochea}.

\bigskip

The aim of this paper is to complete the classification of $5$-dimensional
nilpotent Jordan algebras over an algebraically closed field of
characteristic not $2$. To this end, we classify $5$-dimensional nilpotent
non-associative Jordan algebras over an algebraically closed field of
characteristic not $2$.

\bigskip

In this paper we describe a method for classifying nilpotent Jordan
algebras. Subsequently we classify, up to isomorphism, all nilpotent Jordan
algebras of dimension up to four and all nilpotent non-associative Jordan
algebras of dimension $5$ over an algebraically closed field $\mathbb{F}$ of
characteristic not $2$. This method is analogous to the \emph{Skjelbred-Sund}
method for classifying nilpotent Lie algebras (see \cite{Graaf}, \cite{Gong}%
,\ \cite{Skjelbred}). We introduce a new invariant called the \emph{%
characteristic sequence of dimensions of radicals}. It offers us a very
effective way to distinguish two algebras.

\bigskip

The paper is organized as follows. In Section \ref{Annihilator Extension} we
introduce the notion of annihilator extension. In Sections \ref{clas.
method.}, \ref{construct} we describe a method for classification of
nilpotent Jordan algebras. In Section \ref{dim3} the classification of
nilpotent Jordan algebras up to dimension three is given. Section \ref{dim4}
contains a complete classification of $4$-dimensional nilpotent Jordan
algebras. Section \ref{dim5} is devoted to classify $5$-dimensional
nilpotent non-associative Jordan algebras.

\section{Annihilator extension of Jordan algebra}

\label{Annihilator Extension}

\begin{defn}
A \emph{Jordan algebra} $J$ is a commutative algebra over a field $\mathbb{F}
$ with a multiplication "$\circ $" satisfying the \emph{Jordan identity} 
\begin{equation}
x^{2}\circ (x\circ y)=x\circ (x^{2}\circ y)\mbox{ for all }x,y\in J.
\label{jor id1}
\end{equation}
\end{defn}

The linearization of the Jordan identity $\left( \ref{jor id1}\right) $ is%
\begin{equation}
\left( x,y,z\circ w\right) +\left( w,y,z\circ x\right) +\left( z,y,x\circ
w\right) =0  \label{lin.jor.id}
\end{equation}%
for all $x,y,z,w\in J$. Here $\left( x,y,z\right) =\left( x\circ y\right)
\circ z-x\circ \left( y\circ z\right) $ is the associator of $x,y,z.$

\begin{defn}
Let $S$ be a subset of a Jordan algebra $J$. The set $Ann\left( S\right)
=\left\{ x\in J:x\circ S=0\right\} $ is called the\emph{\ annihilator} of $S$%
.
\end{defn}

In a Jordan algebra $J$ we define inductively a series of subsets by setting 
$J^{\left\langle 1\right\rangle }=J$ and $J^{\left\langle n\right\rangle
}=J^{\left\langle n-1\right\rangle }\circ J$. The chain of subsets%
\begin{equation*}
J^{\left\langle 1\right\rangle }\supseteq J^{\left\langle 2\right\rangle
}\supseteq \cdots \supseteq J^{\left\langle n\right\rangle }\supseteq \cdots
\end{equation*}%
is a chain of ideals of the algebra $J$.

\begin{defn}
A Jordan algebra $J$ is said to be \emph{nilpotent} if there exists an
integer $n\in 
\mathbb{N}
$ such that $J^{\left\langle n\right\rangle }=0$, and the minimal such
number is called the \emph{nilindex} of $J$.
\end{defn}

For a given nilpotent Jordan algebra $J$ of nilindex $m$, we define the 
\emph{nilpotency type} of the algebra $J$ as the sequence $\left(
n_{1},n_{2},\ldots ,n_{m-1}\right) $, where $n_{i}=\dim \left(
J^{\left\langle i\right\rangle }/J^{\left\langle i+1\right\rangle }\right) $%
. The nilpotency type of $J$ is an invariant of the isomorphism class of $J$.

\begin{lem}
\label{same. Ann}Let $\phi :J_{1}\longrightarrow J_{2}$ be an isomorphism of
Jordan algebras. Then, for each $n\in 
\mathbb{N}
$, $\phi \left( Ann\left( J_{1}^{\left\langle n\right\rangle }\right)
\right) =Ann\left( J_{2}^{\left\langle n\right\rangle }\right) $.
\end{lem}

\begin{proof}
Since $\phi $ is an isomorphism, it follows that $\phi \left( x\right) \circ
\phi \left( J_{1}^{\left\langle n\right\rangle }\right) =0$ if and only if $%
x\circ J_{1}^{\left\langle n\right\rangle }=0$. Moreover, $\phi \left(
J_{1}^{\left\langle n\right\rangle }\right) =J_{2}^{\left\langle
n\right\rangle }$ for any $n\in 
\mathbb{N}
$. Then, for each $n\in 
\mathbb{N}
$, $\phi \left( Ann\left( J_{1}^{\left\langle n\right\rangle }\right)
\right) =Ann\left( J_{2}^{\left\langle n\right\rangle }\right) $.
\end{proof}

Therefore, for each $n\in 
\mathbb{N}
$, $\dim Ann\left( J^{\left\langle n\right\rangle }\right) $ is an invariant
of the isomorphism class of $J$.

\begin{exam}
\label{non-isom.}Suppose that $J_{5,10},J_{5,11},J_{5,13},J_{5,14}$ are $5$%
-dimensional nilpotent Jordan algebras defined as follows:%
\begin{table}[H] \centering%
\begin{tabular}[t]{|c|c|c|c|}
\hline
$J$ & Multiplication table (all other products are zero) & $J^{\left\langle
2\right\rangle }$ & $Ann\left( J^{\left\langle 2\right\rangle }\right) $ \\ 
\hline
$J_{5,10}$ & \multicolumn{1}{|l|}{$a\circ b=c,a\circ c=e,d^{2}=e,b\circ c=e.$%
} & $\left\langle c,e\right\rangle $ & $\left\langle c,d,e\right\rangle $ \\ 
\hline
$J_{5,11}$ & \multicolumn{1}{|l|}{$a\circ b=c,a\circ c=e,d^{2}=e.$} & $%
\left\langle c,e\right\rangle $ & $\left\langle b,c,d,e\right\rangle $ \\ 
\hline
$J_{5,13}$ & \multicolumn{1}{|l|}{$a\circ b=c,a\circ c=e,a\circ d=e,b\circ
c=e.$} & $\left\langle c,e\right\rangle $ & $\left\langle c,d,e\right\rangle 
$ \\ \hline
$J_{5,14}$ & \multicolumn{1}{|l|}{$a\circ b=c,a\circ c=e,b\circ d=e.$} & $%
\left\langle c,e\right\rangle $ & $\left\langle b,c,d,e\right\rangle $ \\ 
\hline
\end{tabular}%
\end{table}%
Then $J_{5,10}\ncong J_{5,11},J_{5,10}\ncong J_{5,14},J_{5,11}\ncong
J_{5,13} $ and $J_{5,13}\ncong J_{5,14}$.
\end{exam}

Let $J$ be a Jordan algebra, $V$ be a vector space over a field $\mathbb{F}$%
, and $\theta :J\times J\longrightarrow V$ be a bilinear map. Set $J_{\theta
}=J\oplus V$, and define a multiplication on $J_{\theta }$ by $\left(
x+v\right) \circ _{J_{\theta }}\left( y+w\right) =x\circ y+\theta \left(
x,y\right) $ for $x,y\in J$ and $v,w\in V$. Then $J_{\theta }$ is a
commutative algebra if and only if $\theta $ is a symmetric bilinear map.

\begin{defn}
\label{Coc}Let $J$ be a Jordan algebra and $V$ be a vector space over a
field $\mathbb{F}$. We define $Z^{2}\left( J,V\right) $ to be the set of all
symmetric bilinear maps $\theta :J\times J\rightarrow V$ such that%
\begin{equation}
\theta \left( x^{2},x\circ y\right) =\theta \left( x,x^{2}\circ y\right) %
\mbox{ for all }x,\text{$y\in $}J.  \label{linearize}
\end{equation}
\end{defn}

One can easily see that $Z^{2}\left( J,V\right) $ is a vector space if one
defines the vector space operations as follows. Let $\theta ,\vartheta \in
Z^{2}\left( J,V\right) $, then their linear combination $\alpha \theta
+\beta \vartheta $ ($\alpha ,\beta \in \mathbb{F}$) is defined by $\left(
\alpha \theta +\beta \vartheta \right) \left( x,y\right) =\alpha \theta
\left( x,y\right) +\beta \vartheta \left( x,y\right) $.

\begin{lem}
\label{Jor.Alg.}$J_{\theta }$ is a Jordan algebra if and only if $\theta \in 
$ $Z^{2}\left( J,V\right) $.
\end{lem}

\begin{proof}
Let $x,y\in J$ and $v,w\in V$. Then 
\begin{equation*}
\left( x+v\right) ^{2}\circ _{J_{\theta }}\left( \left( x+v\right) \circ
_{J_{\theta }}\left( y+w\right) \right) -\left( x+v\right) \circ _{J_{\theta
}}\left( \left( x+v\right) ^{2}\circ _{J_{\theta }}\left( y+w\right) \right)
=\theta \left( x^{2},x\circ y\right) -\theta \left( x,x^{2}\circ y\right) .
\end{equation*}%
Hence $J_{\theta }$ is a Jordan algebra if and only if $\theta $ is a
symmetric bilinear map with $\theta \left( x^{2},x\circ y\right) =\theta
\left( x,x^{2}\circ y\right) $ for all $x,y\in J$.
\end{proof}

We call $J_{\theta }$ an $s$\emph{-dimensional annihilator extension} of $J$
by $V$ if $\theta \in $ $Z^{2}\left( J,V\right) $ and $\dim V=s$.

\begin{lem}
\label{as sum}Let $e_{1},\ldots ,e_{s}$ be a basis of $V$, and $\theta \in
Z^{2}\left( J,V\right) $. Then $\theta $ can be uniquely written as $\theta
\left( x,y\right) =\underset{i=1}{\overset{s}{\sum }}\theta _{i}\left(
x,y\right) e_{i}$, where $\theta _{i}\in Z^{2}\left( J,\mathbb{F}\right) $.
\end{lem}

\begin{proof}
Consider any $x,y\in J$ then $\theta \left( x,y\right) $ can be uniquely
written as $\theta \left( x,y\right) =$ $\underset{i=1}{\overset{s}{\sum }}%
\alpha _{i}e_{i}$, where $\alpha _{i}\in \mathbb{F}$. For each $i=1,\ldots
,s $, define a bilinear form $\theta _{i}:J\times J\longrightarrow \mathbb{F}
$ by $\theta _{i}\left( x,y\right) =\alpha _{i}$. Then $\theta \left(
x,y\right) =\underset{i=1}{\overset{s}{\sum }}\theta _{i}\left( x,y\right)
e_{i}$. Moreover, $\underset{i=1}{\overset{s}{\sum }}\theta _{i}\left(
x,y\right) e_{i}=\underset{i=1}{\overset{s}{\sum }}\theta _{i}\left(
y,x\right) e_{i}$ and $\underset{i=1}{\overset{s}{\sum }}\theta _{i}\left(
x^{2},x\circ y\right) e_{i}=\underset{i=1}{\overset{s}{\sum }}\theta
_{i}\left( x,x^{2}\circ y\right) e_{i}$. Therefore, for every $i=1,\ldots ,s$%
, $\theta _{i}\in Z^{2}\left( J,\mathbb{F}\right) $. For uniqueness, let $%
\theta \left( x,y\right) =\underset{i=1}{\overset{s}{\sum }}\vartheta
_{i}\left( x,y\right) e_{i}$. Then $\underset{i=1}{\overset{s}{\sum }}\left(
\theta _{i}-\vartheta _{i}\right) \left( x,y\right) e_{i}=0$. Since $%
e_{1},\ldots ,e_{s}$ are linearly independent, it follows that $\left(
\theta _{i}-\vartheta _{i}\right) \left( x,y\right) =0$ for all $x,y\in J$, $%
i=1,\ldots ,s$. Hence $\theta _{i}-\vartheta _{i}=0$, $i=1,\ldots ,s$.
\end{proof}

Let $a_{1},a_{2},...,a_{n}$ be a basis of $J$. Then by $\delta _{a_{i},a_{j}}
$\ we denote the symmetric bilinear form $\delta _{a_{i},a_{j}}:J\times
J\longrightarrow \mathbb{F}$ with $\delta _{a_{i},a_{j}}\left(
a_{l},a_{m}\right) =1$ if $\left\{ i,j\right\} =\left\{ l,m\right\} $, and
takes the value $0$ otherwise. Denote the space of all symmetric bilinear
forms on $J$ by $Sym\left( J,\mathbb{F}\right) $. Then $Sym\left( J,\mathbb{F%
}\right) =\left\langle \delta _{a_{i},a_{j}}:1\leq i\leq j\leq
n\right\rangle $. Since $Z^{2}\left( J,\mathbb{F}\right) $ is a subspace of $%
Sym\left( J,\mathbb{F}\right) $, every $\theta \in Z^{2}\left( J,\mathbb{F}%
\right) $ can be uniquely written as $\theta =\underset{1\leq i\leq j\leq n}{%
\sum }c_{ij}\delta _{a_{i},a_{j}}$, where $c_{ij}\in \mathbb{F}$. Further,
let $\theta =\underset{1\leq i\leq j\leq n}{\sum }c_{ij}\delta
_{a_{i},a_{j}}\in Sym\left( J,\mathbb{F}\right) $. Then $\theta \in
Z^{2}\left( J,\mathbb{F}\right) $ if and only if the $c_{ij}$'s satisfy the
property $\left( \ref{linearize}\right) $. Note that this property\ is not
linear in $x$; it is better to linearize it. For that we have the following
lemma.

\begin{lem}
\label{lincocycle}Let $\theta \in $ $Z^{2}\left( J,V\right) $. Then 
\begin{equation}
T_{\theta }\left( x,y,z\circ w\right) +T_{\theta }\left( w,y,z\circ x\right)
+T_{\theta }\left( z,y,x\circ w\right) =0\mbox{ for all }x,y,z,w\in J
\label{lincoc}
\end{equation}%
where $T_{\theta }\left( x,y,z\right) =\theta \left( x\circ y,z\right)
-\theta \left( x,y\circ z\right) $.
\end{lem}

\begin{proof}
Let $\theta \in $ $Z^{2}\left( J,V\right) $. Then, by Lemma \ref{Jor.Alg.}, $%
J_{\theta }$ is a Jordan algebra. We denote\ the associator of $x,y,z$ in $%
J_{\theta }$ by $\left( x,y,z\right) _{J_{\theta }}$. Consider any $%
x,y,z,w\in J$ then%
\begin{eqnarray}
\left( x,y,z\circ _{J_{\theta }}w\right) _{J_{\theta }} &=&\left( x,y,z\circ
w\right) +T_{\theta }\left( x,y,z\circ w\right) ,  \label{Q1} \\
\left( w,y,z\circ _{J_{\theta }}x\right) _{J_{\theta }} &=&\left( w,y,z\circ
x\right) +T_{\theta }\left( w,y,z\circ x\right) ,  \label{Q2} \\
\left( z,y,x\circ _{J_{\theta }}w\right) _{J_{\theta }} &=&\left( z,y,x\circ
w\right) +T_{\theta }\left( z,y,x\circ w\right) .  \label{Q3}
\end{eqnarray}%
By the linearized Jordan identity (\ref{lin.jor.id}) we have%
\begin{eqnarray}
\left( x,y,z\circ w\right) +\left( w,y,z\circ x\right) +\left( z,y,x\circ
w\right) &=&0,  \label{Q4} \\
\left( x,y,z\circ _{J_{\theta }}w\right) _{J_{\theta }}+\left( w,y,z\circ
_{J_{\theta }}x\right) _{J_{\theta }}+\left( z,y,x\circ _{J_{\theta
}}w\right) _{J_{\theta }} &=&0.  \label{Q5}
\end{eqnarray}%
From Eqs. $\left( \ref{Q1}\right) $, $\left( \ref{Q2}\right) $, $\left( \ref%
{Q3}\right) $, $\left( \ref{Q4}\right) $ and $\left( \ref{Q5}\right) $, we
obtain%
\begin{equation*}
T_{\theta }\left( x,y,z\circ w\right) +T_{\theta }\left( w,y,z\circ x\right)
+T_{\theta }\left( z,y,x\circ w\right) =0.
\end{equation*}
\end{proof}

Note that $\left( \ref{linearize}\right) $ can be obtained from $\left( \ref%
{lincoc}\right) $ by taking $z=w=x$ in $\left( \ref{lincoc}\right) $
providing the characteristic of $\mathbb{F}$ is not three.

\begin{exam}
\label{Ex1}Let $J:a_{1}^{2}=a_{2}$ be a $2$-dimensional nilpotent Jordan
algebra. Here $0=J^{\left\langle 3\right\rangle }\subset J^{\left\langle
2\right\rangle }=\left\langle a_{2}\right\rangle \subset J=\left\langle
a_{1},a_{2}\right\rangle $. Let $\theta =c_{11}\delta
_{a_{1},a_{1}}+c_{12}\delta _{a_{1},a_{2}}+c_{22}\delta _{a_{2},a_{2}}\in
Z^{2}\left( J,\mathbb{F}\right) $. Then $c_{22}=\theta \left(
a_{2},a_{2}\right) =\theta \left( a_{1}^{2},a_{1}^{2}\right) =\theta \left(
a_{1},a_{1}^{2}\circ a_{1}\right) =0$. Hence $Z^{2}\left( J,\mathbb{F}%
\right) $ is spanned by $\delta _{a_{1},a_{1}},\delta _{a_{1},a_{2}}$.
\end{exam}

Let $\theta \in Z^{2}\left( J,V\right) $. The set $\theta ^{\perp }=\left\{
x\in J:\theta \left( x,y\right) =0\mbox{ for all }\text{$y\in $}J\right\} $
is called the \emph{radical} of $\theta $.

\begin{lem}
\label{Z(J+V)}Let $\theta \in Z^{2}\left( J,V\right) $. Then $Ann\left(
J_{\theta }\right) =\big(\theta ^{\bot }\cap Ann\left( J\right) \big)\oplus
V $. Furthermore, if $\theta \left( x,y\right) =\underset{i=1}{\overset{s}{%
\sum }}\theta _{i}\left( x,y\right) e_{i}$ with $\theta _{i}\in Z^{2}\left(
J,\mathbb{F}\right) $, then $Ann\left( J_{\theta }\right) =\big(\underset{i=1%
}{\overset{s}{\cap }}\theta _{i}^{\bot }\cap Ann\left( J\right) \big)\oplus
V $.
\end{lem}

\begin{proof}
Obviously, $V\subseteq Ann\left( J_{\theta }\right) $. So we can write $%
Ann\left( J_{\theta }\right) =W\oplus V$, where $W\subseteq J$. Hence $x\in W
$ if and only if $x\circ y=\theta \left( x,y\right) =0$ for all $y\in J$.
So, $W=\theta ^{\bot }\cap Ann\left( J\right) $. Further, let $\theta \left(
x,y\right) =\underset{i=1}{\overset{s}{\sum }}\theta _{i}\left( x,y\right)
e_{i}$. Then $\theta (x,y)=0$ if and only if $\theta _{1}\left( x,y\right)
=\cdots =\theta _{s}\left( x,y\right) =0$. So $\theta ^{\perp }=\theta
_{1}^{\bot }\cap \cdots \cap \theta _{s}^{\bot }$.
\end{proof}

\begin{cor}
\label{Z(J)=V}Let $\theta \in Z^{2}\left( J,V\right) $. Then $Ann\left(
J_{\theta }\right) =V$ if and only if $\theta ^{\bot }\cap Ann\left(
J\right) =0$.
\end{cor}

\begin{exam}
\label{rank2}Let $J$ be a Jordan algebra with a basis $a,b,c$ such that $%
a^{2}=b$ and all other products are zero. Then $Ann\left( J\right)
=\left\langle b,c\right\rangle $. Let $\theta =\alpha _{1}\delta
_{a,b}+\alpha _{2}\delta _{a,c}+\alpha _{3}\delta _{b,c}+\alpha _{4}\delta
_{c,c}\in Z^{2}\left( J,\mathbb{F}\right) $. If $x=\lambda _{1}b+\lambda
_{2}c\in \theta ^{\bot }$, then%
\begin{eqnarray*}
\alpha _{1}\lambda _{1}+\alpha _{2}\lambda _{2} &=&0, \\
\alpha _{3}\lambda _{2} &=&0, \\
\alpha _{3}\lambda _{1}+\alpha _{4}\lambda _{2} &=&0.
\end{eqnarray*}%
Therefore, $\theta ^{\bot }\cap Ann\left( J\right) =0$ if and only if the
matrix $%
\begin{bmatrix}
\alpha _{1} & \alpha _{2} \\ 
0 & \alpha _{3} \\ 
\alpha _{3} & \alpha _{4}%
\end{bmatrix}%
$ has rank $2$ or, equivalently, $\left( \alpha _{3},\alpha _{1}\alpha
_{4}\right) \neq \left( 0,0\right) $.
\end{exam}

Lemma \ref{Z(J+V)}\ shows that any annihilator extension of\ a lower
dimensional Jordan algebra has a non-trivial annihilator. Our next aim is to
prove the converse of this statement; that is to say, any Jordan algebra
with a non-trivial annihilator is an annihilator extension of a lower
dimensional Jordan algebra.

\begin{thm}
Let $J$ be a Jordan algebra with $Ann\left( J\right) \neq 0$. Then there
exist, up to isomorphism, a unique Jordan algebra $M$, and a $\theta \in
Z^{2}\left( M,Ann\left( J\right) \right) $ with $\theta ^{\bot }\cap
Ann\left( M\right) =0$ such that $J\cong M_{\theta }$ and $J/Ann\left(
J\right) \cong M$.
\end{thm}

\begin{proof}
Let $J$ be a Jordan algebra with $Ann\left( J\right) \neq 0$. Set $%
V=Ann\left( J\right) $. Let $M$ be a complement of $V$ in $J$; that is to
say, $J=M\oplus V$. Then every $z\in J$ can be uniquely written as $z=x+v$,
where $x\in M$ and $v\in V$. Let $P:J\longrightarrow M$ be the projection of 
$J$ onto $M$ along $V$ (i.e., $P\left( x+v\right) =x$ for $x\in M$ and $x\in
V$). Define a multiplication on $M$ by $x\circ _{M}y=P\left( x\circ y\right) 
$ for $x,y\in M$. Then%
\begin{eqnarray*}
P\left( x\circ y\right) &=&P\left( \left( P\left( x\right) +x-P\left(
x\right) \right) \circ \left( P\left( y\right) +y-P\left( y\right) \right)
\right) \\
&=&P\left( P\left( x\right) \circ P\left( y\right) \right) \\
&=&P\left( x\right) \circ _{M}P\left( y\right)
\end{eqnarray*}%
for all $x,y\in $$J$. Hence $P$ is a homomorphism of algebras. So $P\left(
J\right) =M$ is a Jordan algebra and $J/V\cong M$. Define a symmetric
bilinear map $\theta :M\times M\longrightarrow V$ by $\theta \left(
x,y\right) =$ $x\circ y-x\circ _{M}y$ for $x,y\in M$. Then $M_{\theta }$ is
a commutative algebra. Moreover, $\left( x+v\right) \circ _{M_{\theta
}}\left( y+w\right) =x\circ y$ for all $x,y\in M$ and $v,w\in V$. Hence $%
M_{\theta }$ is the same Jordan algebra as $J$. So by Lemma \ref{Jor.Alg.}, $%
\theta \in Z^{2}(M,V)$. Moreover, by Corollary \ref{Z(J)=V},$\ \theta ^{\bot
}\cap Ann\left( M\right) =0$ since $Ann\left( M_{\theta }\right) =Ann\left(
J\right) =V$.
\end{proof}

Let $J$ be a Jordan algebra and $V$ be a vector space. For a linear map $f$
from $J$ to $V$, if we define $\delta f\colon J\times J\rightarrow V$ by $%
\delta f\left( x,y\right) =f(x\circ y)$, then $\delta f\in Z^{2}\left(
J,V\right) $. Consequently, the map $\delta :f\longrightarrow \delta f\ $is
a map from $Hom\left( J,V\right) $ to $Z^{2}\left( J,V\right) $. Define $%
B^{2}\left( J,V\right) =\left\{ \theta \in Z^{2}\left( J,V\right) :\theta
=\delta f\ \mbox{ for some }f\in Hom\left( J,V\right) \right\} $. One can
easily check that $B^{2}(J,V)$ is a subspace of $Z^{2}(J,V)$.

\begin{lem}
\label{as sum of coboundaries}Let $\theta \in Z^{2}\left( J,V\right) $.
Suppose that $\theta \left( x,y\right) =\underset{i=1}{\overset{s}{\sum }}%
\theta _{i}\left( x,y\right) e_{i}$, where $\theta _{i}\in Z^{2}\left( J,%
\mathbb{F}\right) $. Then $\theta \in B^{2}\left( J,V\right) $\ if and only
if all $\theta _{i}\in B^{2}\left( J,\mathbb{F}\right) $.
\end{lem}

\begin{proof}
Let $\theta \in B^{2}\left( J,V\right) $. Then $\theta =\delta f$ for some $%
f\in Hom\left( J,V\right) $. Consider any $x\in J$\ then $f\left( x\right) $%
\ can be uniquely written as $f\left( x\right) =\underset{i=1}{\overset{s}{%
\sum }}\alpha _{i}e_{i}$, where $\alpha _{i}\in \mathbb{F}$. For each $%
i=1,\ldots ,s$, define a linear map $f_{i}:J\longrightarrow \mathbb{F}$ by $%
f_{i}\left( x\right) =\alpha _{i}$ for $x\in J$. Hence $f\left( x\right) =%
\underset{i=1}{\overset{s}{\sum }}f_{i}\left( x\right) e_{i}$. Moreover, we
have $\underset{i=1}{\overset{s}{\sum }}\theta _{i}\left( x,y\right) e_{i}=%
\underset{i=1}{\overset{s}{\sum }}\delta f_{i}\left( x,y\right) e_{i}$.
Therefore, for every $i=1,\ldots ,s$, $\theta _{i}=\delta f_{i}$.
Conversely, let $\theta _{1},\theta _{2},\ldots ,\theta _{s}\in B^{2}\left(
J,\mathbb{F}\right) $. Then there exist maps $f_{1},f_{2},\ldots ,f_{s}\in
Hom(J,\mathbb{F})$ such that $\theta _{i}=\delta f_{i},$ $i=1,\ldots ,s$.
Define a linear map $f:J\longrightarrow V$ by $f\left( x\right) =\underset{%
i=1}{\overset{s}{\sum }}f_{i}\left( x\right) e_{i}$ for $x\in J$. Then $%
\theta \left( x,y\right) =\delta f\left( x,y\right) $. Hence $\theta =\delta
f$.
\end{proof}

A further useful fact is the following.

\begin{lem}
\label{cobound}Let $J$ be an $n$-dimensional Jordan algebra, and let $%
a_{1},a_{2},\ldots ,a_{m}$ be a basis of $J^{\left\langle 2\right\rangle }$.
Then $B^{2}\left( J,\mathbb{F}\right) =\left\langle \delta a_{1}^{\ast
},\delta a_{2}^{\ast },\ldots ,\delta a_{m}^{\ast }\right\rangle $, where $%
a_{i}^{\ast }(a_{j})=\delta _{ij}$ and $\delta _{ij}$ is the \emph{Kronecker
delta}.
\end{lem}

\begin{proof}
Extend a basis of $J^{\left\langle 2\right\rangle }$ to a basis $%
a_{1},\ldots ,a_{m},a_{m+1},\ldots ,a_{n}$ of $J$. Then $a_{1}^{\ast
},\ldots ,a_{n}^{\ast }$ form a basis of $Hom\left( J,\mathbb{F}\right) $.
Consider any $\delta f$ $\in B^{2}\left( J,\mathbb{F}\right) $, and let $f=%
\underset{i=1}{\overset{n}{\sum }}\alpha _{i}a_{i}^{\ast }$ for some $\alpha
_{i}\in \mathbb{F}$. Then%
\begin{equation*}
\delta f\left( a_{j},a_{k}\right) =\underset{i=1}{\overset{n}{\sum }}\alpha
_{i}a_{i}^{\ast }\left( a_{j}\circ a_{k}\right) =\underset{i=1}{\overset{n}{%
\sum }}\alpha _{i}a_{i}^{\ast }\left( \underset{l=1}{\overset{m}{\sum }}%
\beta _{l}a_{l}\right) =\underset{i=1}{\overset{m}{\sum }}\alpha
_{i}a_{i}^{\ast }\left( a_{j}\circ a_{k}\right) =\underset{i=1}{\overset{m}{%
\sum }}\alpha _{i}\delta a_{i}^{\ast }\left( a_{j},a_{k}\right) .
\end{equation*}%
Hence $\delta f=\underset{i=1}{\overset{m}{\sum }}\alpha _{i}\delta
a_{i}^{\ast }$. Further, let $\alpha _{1},\cdots ,\alpha _{m}\in \mathbb{F}$
be such that $\underset{i=1}{\overset{m}{\sum }}\alpha _{i}\delta
a_{i}^{\ast }=0$. Then $\underset{i=1}{\overset{m}{\sum }}\alpha _{i}\delta
a_{i}^{\ast }\left( J,J\right) =\underset{i=1}{\overset{m}{\sum }}\alpha
_{i}\delta a_{i}^{\ast }\left( J^{\left\langle 2\right\rangle }\right) =0$.
This implies that $\alpha _{1}=\alpha _{2}=\cdots =\alpha _{m}=0$. So $%
B^{2}\left( J,\mathbb{F}\right) =\left\langle \delta a_{1}^{\ast },\delta
a_{2}^{\ast },\ldots ,\delta a_{m}^{\ast }\right\rangle $.
\end{proof}

\begin{exam}
Let $J$ be as in Example \ref{Ex1}. Then $Z^{2}\left( J,\mathbb{F}\right)
=\left\langle \delta _{a_{1},a_{1}},\delta _{a_{1},a_{2}}\right\rangle $. By
Lemma$\ $\ref{cobound}, $B^{2}\left( J,\mathbb{F}\right) =\left\langle
\delta a_{2}^{\ast }\right\rangle $. Since $B^{2}\left( J,\mathbb{F}\right) $
is a subspace of $Z^{2}\left( J,\mathbb{F}\right) $, $\delta a_{3}^{\ast
}\in Z^{2}\left( J,\mathbb{F}\right) $. Let $\delta a_{2}^{\ast }=\alpha
\delta _{a_{1},a_{1}}+\beta \delta _{a_{1},a_{2}}$. Then $\alpha =\delta
a_{2}^{\ast }\left( a_{1},a_{1}\right) =a_{2}^{\ast }\left( a_{1}\circ
a_{1}\right) =a_{2}^{\ast }\left( a_{2}\right) =1$ and $\beta =\delta
a_{2}^{\ast }\left( a_{1},a_{2}\right) =a_{2}^{\ast }\left( a_{1}\circ
a_{2}\right) =0$. Therefore, $B^{2}\left( J,\mathbb{F}\right) $ is spanned
by $\delta _{a_{1},a_{1}}$.
\end{exam}

\begin{defn}
\label{coho}Let $J$ be a Jordan algebra and $V$ be a vector space. We define 
$H^{2}\left( J,V\right) $ to be the quotient space $Z^{2}\left( J,V\right) %
\big/B^{2}\left( J,V\right) $. The equivalence class of $\theta \in
Z^{2}\left( J,V\right) $ is denoted $\left[ \theta \right] \in H^{2}\left(
J,V\right) $.
\end{defn}

Note that if we view $V$ as a trivial bimodule for $J$, then $H^{2}\left(
J,V\right) $ will be the \emph{second cohomology group} of $J$ with
coefficients in $V$ (see \cite{Jacobson68}).

\begin{exam}
Let $J:a_{1}^{2}=a_{3},a_{2}^{2}=a_{3}$ be a $3$-dimensional nilpotent
Jordan algebra. Let $\theta =c_{11}\delta _{a_{1},a_{1}}+c_{12}\delta
_{a_{1},a_{2}}+c_{22}\delta _{a_{2},a_{2}}+c_{13}\delta
_{a_{1},a_{3}}+c_{23}\delta _{a_{2},a_{3}}+c_{33}\delta _{a_{3},a_{3}}\in
Z^{2}\left( J,\mathbb{F}\right) $. Then $c_{33}=\theta \left(
a_{1}^{2},a_{1}^{2}\right) =\theta \left( a_{1},a_{1}^{2}\circ a_{1}\right)
=0$. Hence $Z^{2}\left( J,\mathbb{F}\right) $ is spanned by $\delta
_{a_{1},a_{1}},\delta _{a_{1},a_{2}},\delta _{a_{2},a_{2}},\delta
_{a_{1},a_{3}},\delta _{a_{2},a_{3}}$. By Lemma$\ $\ref{cobound}, $%
B^{2}\left( J,\mathbb{F}\right) =\left\langle \delta a_{3}^{\ast
}\right\rangle $. Let $\delta a_{3}^{\ast }=\lambda _{1}\delta
_{a_{1},a_{1}}+\lambda _{2}\delta _{a_{1},a_{2}}+\lambda _{3}\delta
_{a_{2},a_{2}}+\lambda _{4}\delta _{a_{1},a_{3}}+\lambda _{5}\delta
_{a_{2},a_{3}}$. Then $\lambda _{2}=\lambda _{4}=\lambda _{5}=0,\lambda
_{1}=\lambda _{3}=1$. Therefore, $B^{2}\left( J,\mathbb{F}\right) $ is
spanned by $\delta _{a_{1},a_{1}}+\delta _{a_{2},a_{2}}$. Hence $\left[
\delta _{a_{1},a_{1}}\right] ,\left[ \delta _{a_{1},a_{2}}\right] ,\left[
\delta _{a_{1},a_{3}}\right] ,\left[ \delta _{a_{2},a_{3}}\right] $ form a
basis of $H^{2}\left( J,\mathbb{F}\right) $.
\end{exam}

\begin{cor}
\label{ann. comp.}Let $\theta \left( x,y\right) =\underset{i=1}{\overset{s}{%
\sum }}\theta _{i}\left( x,y\right) e_{i}$ and $\vartheta \left( x,y\right) =%
\underset{i=1}{\overset{s}{\sum }}\vartheta _{i}\left( x,y\right) e_{i}$ be
two elements of $Z^{2}\left( J,V\right) $. Then $\left[ \theta \right] =%
\left[ \vartheta \right] $ if and only if $\left[ \theta _{1}\right] =\left[
\vartheta _{1}\right] ,\left[ \theta _{2}\right] =\left[ \vartheta _{2}%
\right] ,\ldots ,\left[ \theta _{s}\right] =\left[ \vartheta _{s}\right] $.
\end{cor}

\begin{proof}
It is an immediate consequence of Lemma \ref{as sum of coboundaries}.
\end{proof}

Next we show that the isomorphism type of $J_{\theta }$ only depends on the
element $\left[ \theta \right] $ of $H^{2}\left( J,V\right) $.

\begin{lem}
\label{eq. cocycles}Let $\theta ,\vartheta \in $ $Z^{2}\left( J,V\right) $
such that $\left[ \theta \right] =\left[ \vartheta \right] $. Then $\theta
^{\bot }\cap Ann\left( J\right) =\vartheta ^{\bot }\cap Ann\left( J\right) $
or, equivalently, $Ann\left( J_{\theta }\right) =Ann\left( J_{\vartheta
}\right) $. Furthermore, $J_{\theta }\cong J_{\vartheta }$.
\end{lem}

\begin{proof}
Since $\left[ \theta \right] =\left[ \vartheta \right] $, $\vartheta =\theta
+\delta f$ for some $f\in Hom\left( J,V\right) $. So $\vartheta \left(
x,y\right) =\theta \left( x,y\right) +f\left( x\circ y\right) $ for all $%
x,y\in J$. Hence $\theta \left( x,y\right) =x\circ y=0$ if and only if $%
\vartheta \left( x,y\right) =x\circ y=0$. Therefore $\theta ^{\bot }\cap
Ann\left( J\right) =\vartheta ^{\bot }\cap Ann\left( J\right) $. Then, by
Lemma \ref{Z(J+V)}, $Ann\left( J_{\theta }\right) =Ann\left( J_{\vartheta
}\right) $. Further, define a linear map $\sigma :J_{\theta }\longrightarrow
J_{\vartheta }$ by $\sigma \left( x+v\right) =x+f\left( x\right) +v$ for $%
x\in J$ and $v\in V$. Then $\sigma $ is an isomorphism.
\end{proof}

\section{The classification method}

\label{clas. method.}Let $Aut\left( J\right) $ be the automorphism group of
a Jordan algebra $J$. Let $\phi \in Aut\left( J\right) $. For $\theta \in
Z^{2}\left( J,V\right) $ define $\phi \theta \left( x,y\right) =\theta
\left( \phi \left( x\right) ,\phi \left( y\right) \right) $. Then $\phi
\theta \in Z^{2}\left( J,V\right) $. So, $Aut\left( J\right) $ acts on $%
Z^{2}\left( J,V\right) $.

\begin{lem}
Let $\phi \in Aut\left( J\right) $ and $\theta \in Z^{2}\left( J,V\right) $.
Then $\phi \theta \in B^{2}\left( J,V\right) $ if and only if $\theta \in
B^{2}\left( J,V\right) $.
\end{lem}

\begin{proof}
Let $\theta =\delta f$ for some $f\in Hom\left( J,V\right) $. Then $\phi
\theta \left( x,y\right) =\theta \left( \phi \left( x\right) ,\phi \left(
y\right) \right) =f\left( \phi \left( x\right) \circ \phi \left( y\right)
\right) =\delta \left( f\circ \phi \right) \left( x,y\right) $. Hence $\phi
\theta \in B^{2}\left( J,V\right) $. Conversely, let $\phi \theta =\delta f$
for some $f\in Hom\left( J,V\right) $. Then $\theta \left( x,y\right) =\phi
\theta \left( \phi ^{-1}\left( x\right) ,\phi ^{-1}\left( y\right) \right)
=f\left( \phi ^{-1}\left( x\right) \circ \phi ^{-1}\left( y\right) \right)
=\delta \left( f\circ \phi ^{-1}\right) \left( x,y\right) $. Hence $\theta
\in B^{2}\left( J,V\right) $.
\end{proof}

Consequently, the automorphism group $Aut\left( J\right) $ acts on $%
H^{2}\left( J,V\right) $.

\bigskip

Let $\phi =\big(a_{ij}\big)\in Aut\left( J\right) $ and $\theta \in
Z^{2}\left( J,\mathbb{F}\right) $. Let $C=\big(c_{ij}\big)$ be the matrix
representing $\theta $ and $C^{\prime }=\big(c_{ij}^{\prime }\big)$ be the
matrix representing $\phi \theta $. Then $C^{\prime }=\phi ^{t}C\phi $.

\begin{exam}
Let $J$ be as in Example \ref{Ex1}. Then $Z^{2}\left( J,\mathbb{F}\right)
=\left\langle \delta _{a_{1},a_{1}},\delta _{a_{1},a_{2}}\right\rangle $ and 
$H^{2}\left( J,\mathbb{F}\right) =\left\langle \left[ \delta _{a_{1},a_{2}}%
\right] \right\rangle $. The automorphism group $Aut\left( J\right) $\
consists of the invertible matrices of the form%
\begin{equation*}
\phi =%
\begin{bmatrix}
a_{11} & 0 \\ 
a_{21} & a_{11}^{2}%
\end{bmatrix}%
.
\end{equation*}%
Let $\theta =\alpha \delta _{a_{1},a_{1}}+\beta \delta _{a_{1},a_{2}}\in
Z^{2}\left( J,\mathbb{F}\right) $. Write $\phi \theta =\alpha ^{\prime
}\delta _{a_{1},a_{1}}+\beta ^{\prime }\delta _{a_{1},a_{2}}$. Then%
\begin{equation*}
\begin{bmatrix}
\alpha ^{\prime } & \beta ^{\prime } \\ 
\beta ^{\prime } & 0%
\end{bmatrix}%
=%
\begin{bmatrix}
a_{11} & a_{21} \\ 
0 & a_{11}^{2}%
\end{bmatrix}%
\begin{bmatrix}
\alpha & \beta \\ 
\beta & 0%
\end{bmatrix}%
\begin{bmatrix}
a_{11} & 0 \\ 
a_{21} & a_{11}^{2}%
\end{bmatrix}%
.
\end{equation*}%
Hence $\phi \theta =a_{11}\left( \alpha a_{11}+2\beta a_{21}\right) \delta
_{a_{1},a_{1}}+\beta a_{11}^{3}\delta _{a_{1},a_{2}}$, and therefore $\left[
\phi \theta \right] =\beta a_{11}^{3}\left[ \delta _{a_{1},a_{2}}\right] $.
\end{exam}

\begin{lem}
Let $\phi \in Aut\left( J\right) $ and $\theta \in Z^{2}\left( J,V\right) $.
Then $\dim \theta ^{\perp }=\dim \left( \phi \theta \right) ^{\perp }$.
\end{lem}

\begin{proof}
Let $x\in J$. Then $\phi \theta \left( x,J\right) =0$ if and only if $\theta
\left( \phi \left( x\right) ,J\right) =0$. So, $x\in \left( \phi \theta
\right) ^{\perp }$ if and only if $\phi \left( x\right) \in \theta ^{\perp }$%
. Define a linear map $\sigma :\left( \phi \theta \right) ^{\perp
}\longrightarrow \theta ^{\perp }$ by $\sigma \left( x\right) =\phi \left(
x\right) $ for $x\in \left( \phi \theta \right) ^{\perp }$. Then $\sigma $
is a bijective map. Therefore, $\dim \theta ^{\perp }=\dim \left( \phi
\theta \right) ^{\perp }$.
\end{proof}

For $\theta _{1},\theta _{2},\ldots ,\theta _{s}\in Z^{2}\left( J,\mathbb{F}%
\right) $, let $\Psi \left( \theta _{1},\theta _{2},\ldots ,\theta
_{s}\right) =\left( m_{1},m_{2},\ldots ,m_{s}\right) $ be the ordered
descending sequence of $\dim \theta _{1}^{\perp },\dim \theta _{2}^{\perp
},...,\dim \theta _{s}^{\perp }$.

\begin{cor}
\label{inv}Let $\phi \in Aut\left( J\right) $ and $\theta _{1},\ldots
,\theta _{s}\in Z^{2}\left( J,\mathbb{F}\right) $. Then $\Psi \left( \theta
_{1},\ldots ,\theta _{s}\right) =\Psi \left( \phi \theta _{1},\ldots ,\phi
\theta _{s}\right) $.
\end{cor}

Let $GL\left( V\right) $\ be the set of all bijective linear maps $%
V\longrightarrow V$. Let $\psi \in GL\left( V\right) $. For $\theta \in
Z^{2}\left( J,V\right) $ define $\psi \theta \left( x,y\right) =\psi \left(
\theta \left( x,y\right) \right) $. Then $\psi \theta \in Z^{2}\left(
J,V\right) $. So, $GL\left( V\right) $ acts on $Z^{2}\left( J,V\right) $.
Also, $\theta \in B^{2}\left( J,V\right) $ if and only if $\psi \theta \in
B^{2}\left( J,V\right) $. Hence $GL\left( V\right) $ acts on $H^{2}\left(
J,V\right) $.

\bigskip

Now let $\theta ,\vartheta \in Z^{2}\left( J,V\right) $ such that $Ann\left(
J_{\theta }\right) =Ann\left( J_{\vartheta }\right) =V$. Suppose that $%
J_{\theta }\cong J_{\vartheta }$. Then there exists an isomorphism $\phi :$ $%
J_{\theta }\longrightarrow J_{\vartheta }$. Since $\phi \left( V\right)
=\phi \left( Ann\left( J_{\theta }\right) \right) =Ann\left( J_{\vartheta
}\right) =V$, $\psi =\phi |_{V}\in GL\left( V\right) $. Let $%
x_{1},x_{2},...,x_{n}$ be a basis of $J$, and let $\phi \left( x_{i}\right)
=y_{i}+v_{i}$, where $y_{i}\in J$ and $v_{i}$ $\in V$. Then $\phi $ induces
an isomorphism $\phi _{0}:J\longrightarrow J$ defined by $\phi _{0}\left(
x_{i}\right) =y_{i}$,\ and a linear map $\varphi :J\longrightarrow V$
defined by $\varphi \left( x_{i}\right) =v_{i}$. So we can realize $\phi $
as a matrix of the form%
\begin{equation*}
\phi =%
\begin{bmatrix}
\phi _{0} & 0 \\ 
\varphi & \psi%
\end{bmatrix}%
:\phi _{0}\in Aut\left( J\right) ,\psi =\phi |_{V}\in GL\left( V\right) %
\mbox{ and }\varphi \in Hom\left( J,V\right) .
\end{equation*}%
Furthermore, for any $x,y\in $$J$ we have%
\begin{eqnarray*}
\phi \left( x\circ _{J_{\theta }}y\right) &=&\phi \left( x\circ y+\theta
\left( x,y\right) \right) =\phi _{0}(x\circ y)+\varphi (x\circ y)+\psi
(\theta \left( x,y\right) ), \\
\phi \left( x\right) \circ _{J_{\vartheta }}\phi \left( y\right) &=&\left(
\phi _{0}\left( x\right) +\varphi \left( x\right) \right) \circ
_{J_{\vartheta }}\left( \phi _{0}\left( y\right) +\varphi \left( y\right)
\right) =\phi _{0}\left( x\right) \circ \phi _{0}\left( y\right) +\vartheta
\left( \phi _{0}\left( x\right) ,\phi _{0}\left( y\right) \right) .
\end{eqnarray*}%
Since $\phi $ is an isomorphism, it follows that%
\begin{equation}
\vartheta \left( \phi _{0}\left( x\right) ,\phi _{0}\left( y\right) \right)
=\varphi \left( x\circ y\right) +\psi \left( \theta \left( x,y\right)
\right) \mbox{ for all }\text{$x,y\in $}J.  \label{auto}
\end{equation}%
This is equivalent to $\phi _{0}\vartheta =\delta \varphi +\psi \theta $.
Hence $\left[ \phi _{0}\vartheta \right] =\left[ \psi \theta \right] $. Then
we have the following lemma.

\begin{lem}
\label{iso4}Let $\theta ,\vartheta \in Z^{2}\left( J,V\right) $ such that $%
Ann\left( J_{\theta }\right) =Ann\left( J_{\vartheta }\right) =V$. Then $%
J_{\theta }\cong J_{\vartheta }$ if and only if there exist a map $\phi \in
Aut\left( J\right) $ and a map $\psi \in GL(V)$ such that $\left[ \phi
\theta \right] =\left[ \psi \vartheta \right] $.
\end{lem}

In case of $\theta =\vartheta $, we obtain from $\left( \ref{auto}\right) $
the following description of $Aut(J_{\theta })$.

\begin{lem}
Let $J$ be a Jordan algebra and $\theta \in Z^{2}\left( J,V\right) $ such
that $\theta ^{\bot }\cap Ann\left( J\right) =0$. Then the automorphism
group $Aut(J_{\theta })$ of the Jordan algebra $J_{\theta }$ consists of all
linear maps of the matrix form 
\begin{equation*}
\phi =%
\begin{bmatrix}
\phi _{0} & 0 \\ 
\varphi & \psi%
\end{bmatrix}%
:\phi _{0}\in Aut\left( J\right) ,\psi \in GL\left( V\right) \mbox{ and }%
\varphi \in Hom\left( J,V\right) \text{ }
\end{equation*}%
such that $\theta (\phi _{0}\left( x\right) ,\phi _{0}\left( y\right)
)=\varphi (x\circ y)+\psi \left( \theta \left( x,y\right) \right) \ $for all 
$x,y\in J$.
\end{lem}

\begin{defn}
\label{Annihilator component}Let $J=I\oplus \mathbb{F}x$ be the direct sum
of two ideals. Then $\mathbb{F}x$ is called an \emph{annihilator component}
of $J$ if $x\in Ann\left( J\right) $.
\end{defn}

Let $\phi :J_{1}\longrightarrow J_{2}$ be an isomorphism of Jordan algebras.
Let $x\in J_{1}$ and $y\in J_{2}$ such that $\mathbb{\phi }\left( x\right)
=y $. Then $\mathbb{F}x$ is an annihilator component of $J_{1}$ if and only
if $\mathbb{F}y$ is an annihilator component of $J_{2}$.

\begin{lem}
\label{cent.comp}Let $\theta \left( x,y\right) =\underset{i=1}{\overset{s}{%
\sum }}$ $\theta _{i}\left( x,y\right) e_{i}\in Z^{2}\left( J,V\right) $ and 
$\theta ^{\bot }\cap Ann\left( J\right) =0$ . Then $J_{\theta }$ has an
annihilator component if and only if $\left[ \theta _{1}\right] ,\left[
\theta _{2}\right] ,\ldots ,\left[ \theta _{s}\right] $ are linearly
dependent in $H^{2}\left( J,\mathbb{F}\right) $.
\end{lem}

\begin{proof}
Suppose that $J_{\theta }$ has an annihilator component. Then there exists
an element $v_{1}\in V$ such that $J_{\theta }=I\oplus \mathbb{F}v_{1}$.
Enlarge the set $\left\{ v_{1}\right\} $ to this set $\left\{
v_{1},v_{2},...,v_{s}\right\} $ to form a basis of $V$. Then there exists an
invertible matrix $\big(a_{ij}\big)$ such that for any $i=1,\ldots ,s$, $%
e_{i}=\underset{j=1}{\overset{s}{\sum }}a_{ij}v_{j}$. So we can write $%
\theta \left( x,y\right) =\underset{j=1}{\overset{s}{\sum }}\left( \overset{s%
}{\underset{i=1}{\sum }}a_{ij}\theta _{i}\left( x,y\right) \right) v_{j}$.
Since $J_{\theta }^{\left\langle 2\right\rangle }\subset I$, it follows that 
$\overset{s}{\underset{i=1}{\sum }}a_{i1}\theta _{i}\left( x,y\right) =0$
for all $x,y\in J$. Then $\overset{s}{\underset{i=1}{\sum }}a_{i1}\theta
_{i}=0$, and therefore $\overset{s}{\underset{i=1}{\sum }}a_{i1}\left[
\theta _{i}\right] =0$. Now if $a_{11}=a_{21}=\cdots =a_{s1}=0$, then $\det %
\big(a_{ij}\big)=0$ which is a contradiction. So $\left[ \theta _{1}\right] ,%
\left[ \theta _{2}\right] ,\ldots ,\left[ \theta _{s}\right] $ are linearly
dependent in $H^{2}\left( J,\mathbb{F}\right) $. On the other hand, suppose
that $\left[ \theta _{1}\right] ,\left[ \theta _{2}\right] ,\ldots ,\left[
\theta _{s}\right] $ are linearly dependent. Without loss of generality we
may assume that $\left[ \theta _{s}\right] =\underset{i=1}{\overset{s-1}{%
\sum }}\alpha _{i}\left[ \theta _{i}\right] $ for some $\alpha _{i}\in 
\mathbb{F}$. Define $\vartheta \left( x,y\right) =\underset{i=1}{\overset{s}{%
\sum }}\vartheta _{i}\left( x,y\right) $ by setting $\vartheta _{i}=\theta
_{i}$ for $i=1,\ldots ,s-1$ and $\vartheta _{s}=\underset{i=1}{\overset{s-1}{%
\sum }}\alpha _{i}\theta _{i}$. Then $\vartheta \in Z^{2}\left( J,V\right) $%
. Moreover, by Corollary \ref{ann. comp.}, $\left[ \theta \right] =\left[
\vartheta \right] $ and hence $J_{\theta }\cong J_{\vartheta }$. Easy
computation shows that $\vartheta \left( x,y\right) =\underset{i=1}{\overset{%
s-1}{\sum }}\theta _{i}\left( x,y\right) \left( e_{i}+\alpha
_{i}e_{s}\right) $. For $i=1,\ldots ,s-1$, set $w_{i}=e_{i}+\alpha _{i}e_{s}$%
. Then $\vartheta \left( x,y\right) =\overset{s-1}{\underset{i=1}{\sum }}%
\theta _{i}\left( x,y\right) w_{i}$. Hence $J_{\vartheta }^{\left\langle
2\right\rangle }\subset J\oplus \left\langle w_{1},w_{2},\ldots
,w_{s-1}\right\rangle $, so that $J_{\vartheta }$, and therefore also $%
J_{\theta }$, has an annihilator component.
\end{proof}

The statement in Lemma \ref{iso4} can be rephrased as follows.

\begin{lem}
\label{excl}Let $\theta \left( x,y\right) =\underset{i=1}{\overset{s}{\sum }}%
\theta _{i}\left( x,y\right) e_{i}$ and $\vartheta \left( x,y\right) =%
\underset{i=1}{\overset{s}{\sum }}\vartheta _{i}\left( x,y\right) e_{i}$ be
two elements of $Z^{2}\left( J,V\right) $. Suppose that $J_{\theta }$ has no
annihilator components and $\theta ^{\perp }\cap Ann\left( J\right)
=\vartheta ^{\perp }\cap Ann\left( J\right) =0$. Then $J_{\theta }\cong
J_{\vartheta }$ if and only if there exists a map $\phi \in Aut\left(
J\right) $ such that the $\left[ \phi \vartheta _{i}\right] $ span the same
subspace of $H^{2}\left( J,\mathbb{F}\right) $ as the $\left[ \theta _{i}%
\right] $.
\end{lem}

\begin{proof}
Suppose first that $J_{\theta }\cong J_{\vartheta }$. Then, by Lemma \ref%
{iso4}, there exist a map $\phi \in Aut\left( J\right) $ and a map $\psi \in
GL(V)$ such that $\left[ \phi \vartheta \right] =\left[ \psi \theta \right] $%
. Let $\psi \left( e_{i}\right) =\underset{j=1}{\overset{s}{\sum }}%
a_{ij}e_{j}$. Then $\left( \phi \vartheta -\psi \theta \right) \left(
x,y\right) =\underset{j=1}{\overset{s}{\sum }}\left( \phi \vartheta _{j}-%
\underset{i=1}{\overset{s}{\sum }}a_{ij}\theta _{i}\right) \left( x,y\right)
e_{j}$. Since $\left[ \phi \vartheta \right] =\left[ \psi \theta \right] $,
and according to Lemma \ref{as sum of coboundaries}, $\phi \vartheta _{j}-%
\underset{i=1}{\overset{s}{\sum }}a_{ij}\theta _{i}\in B^{2}\left( J,\mathbb{%
F}\right) $ for $j=1,\ldots ,s$. Therefore, for every $j=1,\ldots ,s$, $%
\left[ \phi \vartheta _{j}\right] =\underset{i=1}{\overset{s}{\sum }}a_{ij}%
\left[ \theta _{i}\right] $. Hence the $\left[ \phi \vartheta _{i}\right] $
span the same subspace of $H^{2}\left( J,\mathbb{F}\right) $ as the $\left[
\theta _{i}\right] $. On the other hand, let the $\left[ \phi \vartheta _{i}%
\right] $ span the same subspace of $H^{2}\left( J,\mathbb{F}\right) $ as
the $\left[ \theta _{i}\right] $. Then there exists an invertible matrix $%
\big(a_{ij}\big)$ such that for any $j=1,\ldots ,s$, $\left[ \phi \vartheta
_{j}\right] =\underset{i=1}{\overset{s}{\sum }}a_{ij}\left[ \theta _{i}%
\right] $. Define a linear map $\psi :V\longrightarrow V$ by $\psi \left(
e_{i}\right) =\underset{j=1}{\overset{s}{\sum }}a_{ij}e_{j}$. Then $\psi
\theta \left( x,y\right) =\underset{i=1}{\overset{s}{\sum }}\underset{j=1}{%
\overset{s}{\sum }}a_{ij}\theta _{i}\left( x,y\right) e_{j}$. Moreover, $%
\left( \phi \vartheta -\psi \theta \right) \left( x,y\right) =\underset{j=1}{%
\overset{s}{\sum }}\left( \phi \vartheta _{j}-\underset{i=1}{\overset{s}{%
\sum }}a_{ij}\theta _{i}\right) \left( x,y\right) e_{j}$. Then, by Lemma \ref%
{as sum of coboundaries}, $\left[ \phi \vartheta \right] =\left[ \psi \theta %
\right] $. Hence, by Lemma \ref{iso4}, $J_{\theta }\cong J_{\vartheta }$.
\end{proof}

\section{Analogue of the Skjelbred-Sund method}

\label{construct}Let $V$ be a finite-dimensional vector space over a $%
\mathbb{F}$. The \emph{Grassmannian} $G_{k}\left( V\right) $ is the set of
all $k$-dimensional linear subspaces of $V$. Let $G_{s}\left( H^{2}\left( J,%
\mathbb{F}\right) \right) $ be the Grassmannian of subspaces of dimension $s$
in $H^{2}\left( J,\mathbb{F}\right) $. There is a natural action of $%
Aut\left( J\right) $ on $G_{s}\left( H^{2}\left( J,\mathbb{F}\right) \right) 
$. Let $\phi \in Aut\left( J\right) $. For $W=\left\langle \left[ \theta _{1}%
\right] ,\left[ \theta _{2}\right] ,...,\left[ \theta _{s}\right]
\right\rangle \in G_{s}\left( H^{2}\left( J,\mathbb{F}\right) \right) $
define $\phi W=\left\langle \left[ \phi \theta _{1}\right] ,\left[ \phi
\theta _{2}\right] ,...,\left[ \phi \theta _{s}\right] \right\rangle $. Then 
$\phi W\in G_{s}\left( H^{2}\left( J,\mathbb{F}\right) \right) $. We denote
the orbit of $W\in G_{s}\left( H^{2}\left( J,\mathbb{F}\right) \right) $
under the action of $Aut\left( J\right) $ by $\mbox{Orb}\left( W\right) $.

\begin{defn}
Let $W\in G_{s}\left( H^{2}\left( J,\mathbb{F}\right) \right) $. We define $%
\Psi \left( W\right) $ to be the least upper bound for the lexicographic
order of the set%
\begin{equation*}
\left\{ \Psi \left( \theta _{1},\theta _{2},\ldots ,\theta _{s}\right)
:\theta _{1},\theta _{2},\ldots ,\theta _{s}\in Z^{2}\left( J,\mathbb{F}%
\right) \mbox{ and }W=\left\langle \left[ \theta _{1}\right] ,\left[ \theta
_{2}\right] ,...,\left[ \theta _{s}\right] \right\rangle \right\} .
\end{equation*}%
We call $\Psi \left( W\right) $ the \emph{characteristic sequence of
dimensions of radicals}.
\end{defn}

\begin{lem}
Let $\phi \in Aut\left( J\right) $ and $W\in G_{s}\left( H^{2}\left( J,%
\mathbb{F}\right) \right) $. Then $\Psi \left( W\right) =\Psi \left( \phi
W\right) $.
\end{lem}

\begin{proof}
Let $\Psi \left( W\right) =\Psi \left( \theta _{1},\ldots ,\theta
_{s}\right) $ for some $\theta _{1},\ldots ,\theta _{s}\in Z^{2}\left( J,%
\mathbb{F}\right) $. Then, by corollary \ref{inv}, $\Psi \left( \theta
_{1},\ldots ,\theta _{s}\right) =\Psi \left( \phi \theta _{1},\ldots ,\phi
\theta _{s}\right) $. Therefore, $\Psi \left( W\right) \leq \Psi \left( \phi
W\right) $. On the other hand, let $\Psi \left( \phi W\right) =\Psi \left(
\phi \vartheta _{1},\ldots ,\phi \vartheta _{s}\right) $ for some $\vartheta
_{1},\ldots ,\vartheta _{s}\in Z^{2}\left( J,\mathbb{F}\right) $. Then, by
corollary \ref{inv}, $\Psi \left( \vartheta _{1},\ldots ,\vartheta
_{s}\right) =\Psi \left( \phi \vartheta _{1},\ldots ,\phi \vartheta
_{s}\right) $. Therefore, $\Psi \left( \phi W\right) \leq \Psi \left(
W\right) $. So, $\Psi \left( W\right) =\Psi \left( \phi W\right) $.
\end{proof}

So, if $W_{1},W_{2}\in G_{s}\left( H^{2}\left( J,\mathbb{F}\right) \right) $
such that $\Psi \left( W_{1}\right) \neq \Psi \left( W_{2}\right) $ then $%
\mbox{Orb}\left( W_{1}\right) \cap \mbox{Orb}\left( W_{2}\right) =\emptyset $%
. Therefore, the characteristic sequence of dimensions of radicals $\Psi
\left( W\right) $ is an invariant of $\mbox{Orb}\left( W\right) $.

\begin{exam}
Let $J$ be a Jordan algebra with a basis $a,b,c,d$ such that $a^{2}=b$ and
all other products are zero. Then $B^{2}\left( J,\mathbb{F}\right)
=\left\langle \delta _{a,a}\right\rangle $. Let $W_{1},W_{2}\in G_{s}\left(
H^{2}\left( J,\mathbb{F}\right) \right) $ such that $W_{1}=\left\langle %
\left[ \delta _{d,d}\right] +\left[ \delta _{b,c}\right] \right\rangle $ and 
$W_{2}=\left\langle \left[ \delta _{a,d}\right] +\left[ \delta _{b,c}\right]
\right\rangle $. Then $\Psi \left( W_{1}\right) $ is the least upper bound
of the set $\left\{ \Psi \left( \delta _{d,d}+\delta _{b,c}+\alpha \delta
_{a,a}\right) :\alpha \in \mathbb{F}\right\} $. Therefore, $\Psi \left(
W_{1}\right) =\left( 1\right) $. Also, $\Psi \left( W_{2}\right) $ is the
least upper bound of the set $\left\{ \Psi \left( \delta _{a,d}+\delta
_{b,c}+\alpha \delta _{a,a}\right) :\alpha \in \mathbb{F}\right\} $.
Therefore, $\Psi \left( W_{2}\right) =\left( 0\right) $. Hence $\emph{%
\mbox{Orb}}\left( W_{1}\right) \cap \emph{\mbox{Orb}}\left( W_{2}\right)
=\emptyset $.
\end{exam}

\begin{exam}
Let $J$ be a Jordan algebra with a basis $a,b,c$ such that $a\circ b=c$ and
all other products are zero. Then $B^{2}\left( J,\mathbb{F}\right)
=\left\langle \delta _{a,b}\right\rangle $. Let $W_{1},W_{2}\in G_{s}\left(
H^{2}\left( J,\mathbb{F}\right) \right) $ such that $W_{1}=\left\langle %
\left[ \delta _{b,b}\right] +\left[ \delta _{a,c}\right] ,\left[ \delta
_{a,a}\right] +\left[ \delta _{b,c}\right] \right\rangle $ and $%
W_{2}=\left\langle \left[ \delta _{b,b}\right] +\left[ \delta _{a,c}\right] ,%
\left[ \delta _{a,a}\right] \right\rangle $. Then $\Psi \left( W_{1}\right) $
is the least upper bound of the set%
\begin{equation*}
\left\{ \Psi \left( \delta _{b,b}+\delta _{a,c}+\alpha \left( \delta
_{a,a}+\delta _{b,c}\right) +\beta \delta _{a,b},\alpha ^{\prime }\left(
\delta _{b,b}+\delta _{a,c}\right) +\delta _{a,a}+\delta _{b,c}+\beta
^{\prime }\delta _{a,b}\right) :\alpha ,\alpha ^{\prime },\beta ,\beta
^{\prime }\in \mathbb{F}\right\} .
\end{equation*}%
Therefore, $\Psi \left( W_{1}\right) =\left( 1,1\right) $. Also, $\Psi
\left( W_{2}\right) $ is the least upper bound of the set%
\begin{equation*}
\left\{ \Psi \left( \delta _{b,b}+\delta _{a,c}+\alpha \delta _{a,a}+\beta
\delta _{a,b},\alpha ^{\prime }\left( \delta _{b,b}+\delta _{a,c}\right)
+\delta _{a,a}+\beta ^{\prime }\delta _{a,b}\right) :\alpha ,\alpha ^{\prime
},\beta ,\beta ^{\prime }\in \mathbb{F}\right\} .
\end{equation*}%
Therefore, $\Psi \left( W_{2}\right) =\left( 2,0\right) $. Hence $\emph{%
\mbox{Orb}}\left( W_{1}\right) \cap \emph{\mbox{Orb}}\left( W_{2}\right)
=\emptyset $.
\end{exam}

\begin{lem}
Let $W_{1}=\left\langle \left[ \theta _{1}\right] ,\left[ \theta _{2}\right]
,...,\left[ \theta _{s}\right] \right\rangle ,W_{2}=\left\langle \left[
\vartheta _{1}\right] ,\left[ \vartheta _{2}\right] ,...,\left[ \vartheta
_{s}\right] \right\rangle \in G_{s}\left( H^{2}\left( J,\mathbb{F}\right)
\right) $. If $W_{1}=W_{2}$, then $\underset{i=1}{\overset{s}{\cap }}\theta
_{i}^{\bot }\cap Ann\left( J\right) =\underset{i=1}{\overset{s}{\cap }}%
\vartheta _{i}^{\bot }\cap Ann\left( J\right) $.
\end{lem}

\begin{proof}
Let $W_{1}=W_{2}$. Then there exists an invertible matrix $\big(a_{ij}\big)$
such that for any $i=1,\ldots ,s$, $\left[ \theta _{i}\right] =\underset{j=1}%
{\overset{s}{\sum }}a_{ij}\left[ \vartheta _{j}\right] $. Therefore, $\theta
_{i}=\underset{j=1}{\overset{s}{\sum }}a_{ij}\vartheta _{j}+\delta f_{i}$
for some $f_{i}\in Hom\left( J,\mathbb{F}\right) $. Then, for each $%
i=1,\ldots ,s$, $\theta _{i}\left( x,y\right) =\underset{j=1}{\overset{s}{%
\sum }}a_{ij}\vartheta _{j}\left( x,y\right) +f_{i}\left( x\circ y\right) $.
Hence $\theta _{1}\left( x,y\right) =\cdots =\theta _{s}\left( x,y\right)
=x\circ y=0$ if and only if $\vartheta _{1}\left( x,y\right) =\cdots
=\vartheta _{s}\left( x,y\right) =x\circ y=0$. Therefore $\underset{i=1}{%
\overset{s}{\cap }}\theta _{i}^{\bot }\cap Ann\left( J\right) =\underset{i=1}%
{\overset{s}{\cap }}\vartheta _{i}^{\bot }\cap Ann\left( J\right) $.
\end{proof}

This result allows us to define%
\begin{equation*}
T_{s}\left( J\right) =\left\{ W=\left\langle \left[ \theta _{1}\right] ,%
\left[ \theta _{2}\right] ,...,\left[ \theta _{s}\right] \right\rangle \in
G_{s}\left( H^{2}\left( J,\mathbb{F}\right) \right) :\underset{i=1}{\overset{%
s}{\cap }}\theta _{i}^{\bot }\cap Ann\left( J\right) =0\right\} .
\end{equation*}

\begin{lem}
\label{invariantT}The set $T_{s}\left( J\right) $ is stable under the action
of $Aut\left( J\right) $.
\end{lem}

\begin{proof}
Let $\phi \in Aut\left( J\right) $ and $W=\left\langle \left[ \theta _{1}%
\right] ,\left[ \theta _{2}\right] ,...,\left[ \theta _{s}\right]
\right\rangle \in G_{s}\left( H^{2}\left( J,\mathbb{F}\right) \right) $.
Then, for each $i=1,\ldots ,s$, $x\in (\phi \theta _{i})^{\perp }\cap
Ann\left( J\right) $ if and only if $\phi \left( x\right) \in \theta
_{i}^{\perp }\cap Ann\left( J\right) $. So $\underset{i=1}{\overset{s}{\cap }%
}(\phi \theta _{i})^{\perp }\cap Ann\left( J\right) =0$ if and only if $%
\underset{i=1}{\overset{s}{\cap }}\theta _{i}^{\bot }\cap Ann\left( J\right)
=0$. Hence $\phi W\in T_{s}\left( J\right) $ if and only if $W\in
T_{s}\left( J\right) $.
\end{proof}

Let $V$ be an $s$-dimensional vector space spanned by $e_{1},e_{2},\ldots
,e_{s}$. Given a Jordan algebra $J$, let $E\left( J,V\right) $ denote the
set of all Jordan algebras without annihilator components which are $s$\emph{%
-}dimensional annihilator extensions of $J$ by $V$ and have $s$-dimensional
annihilator. Then, by Lemma \ref{cent.comp} and Corollary \ref{Z(J)=V}, $%
E\left( J,V\right) =\left\{ J_{\theta }:\theta \left( x,y\right) =\underset{%
i=1}{\overset{s}{\sum }}\theta _{i}\left( x,y\right) e_{i}\mbox{
and }\left\langle \left[ \theta _{1}\right] ,\left[ \theta _{2}\right] ,...,%
\left[ \theta _{s}\right] \right\rangle \in T_{s}\left( J\right) \right\} $.
Given $J_{\theta }\in E\left( J,V\right) $, let $\left[ J_{\theta }\right] $
denote the isomorphism class of $J_{\theta }$.

\begin{lem}
Let $J_{\theta },J_{\vartheta }\in E\left( J,V\right) $. Suppose that $%
\theta \left( x,y\right) =\underset{i=1}{\overset{s}{\sum }}\theta
_{i}\left( x,y\right) e_{i}$ and $\vartheta \left( x,y\right) =\underset{i=1}%
{\overset{s}{\sum }}\vartheta _{i}\left( x,y\right) e_{i}$. Then $\left[
J_{\theta }\right] =\left[ J_{\vartheta }\right] $ if and only if $\emph{%
\mbox{Orb}}\left\langle \left[ \theta _{1}\right] ,\left[ \theta _{2}\right]
,...,\left[ \theta _{s}\right] \right\rangle =\emph{\mbox{Orb}}\left\langle %
\left[ \vartheta _{1}\right] ,\left[ \vartheta _{2}\right] ,...,\left[
\vartheta _{s}\right] \right\rangle $.
\end{lem}

\begin{proof}
Let $J_{\theta },J_{\vartheta }\in E\left( J,V\right) $. Then, by Lemma \ref%
{excl}, $J_{\theta }\cong J_{\vartheta }$ if and only if there exists a map $%
\phi \in Aut\left( J\right) $ such that $\left\langle \left[ \phi \vartheta
_{1}\right] ,\left[ \phi \vartheta _{2}\right] ,...,\left[ \phi \vartheta
_{s}\right] \right\rangle =\left\langle \left[ \theta _{1}\right] ,\left[
\theta _{2}\right] ,...,\left[ \theta _{s}\right] \right\rangle $. Hence $%
\left[ J_{\theta }\right] =\left[ J_{\vartheta }\right] $ if and only if $%
\mbox{Orb}\left\langle \left[ \theta _{1}\right] ,\left[ \theta _{2}\right]
,...,\left[ \theta _{s}\right] \right\rangle =\mbox{Orb}\left\langle \left[
\vartheta _{1}\right] ,\left[ \vartheta _{2}\right] ,...,\left[ \vartheta
_{s}\right] \right\rangle $
\end{proof}

Thus, each orbit of $Aut\left( J\right) $ on $T_{s}\left( J\right) $\
corresponds uniquely to an isomorphism class of $E\left( J,V\right) $, and
vice versa. This correspondence is defined by%
\begin{equation*}
\mbox{Orb}\left\langle \left[ \theta _{1}\right] ,\left[ \theta _{2}\right]
,...,\left[ \theta _{s}\right] \right\rangle \in \left\{ \mbox{Orb}\left(
W\right) :W\in T_{s}\left( J\right) \right\} \leftrightarrow \left[
J_{\theta }\right] \in \left\{ \left[ J_{\vartheta }\right] :J_{\vartheta
}\in E\left( J,V\right) \right\} \text{,}
\end{equation*}%
where $\theta \left( x,y\right) =\underset{i=1}{\overset{s}{\sum }}\theta
_{i}\left( x,y\right) e_{i}$. We call $J_{\theta }$ the Jordan algebra
corresponding to the representative $\left\langle \left[ \theta _{1}\right] ,%
\left[ \theta _{2}\right] ,...,\left[ \theta _{s}\right] \right\rangle $.

\begin{thm}
There exists a one-to-one correspondence between the set of $Aut\left(
J\right) $-orbits on $T_{s}\left( J\right) $ and the set of isomorphism
classes of $E\left( J,V\right) $.
\end{thm}

By this theorem, we may construct all the nilpotent Jordan algebras of
dimension $n$, given those algebras of dimension less than $n$, in the
following way:

\begin{enumerate}
\item For a given nilpotent Jordan algebra $J$ of dimension $n-s$, determine 
$H^{2}\left( J,\mathbb{F}\right) $, $Ann\left( J\right) $ and $Aut\left(
J\right) $.

\item Determine the set of $Aut\left( J\right) $-orbits on $T_{s}\left(
J\right) $.

\item For each orbit, construct the Jordan algebra corresponding to a
representative of it.
\end{enumerate}

This method gives all nilpotent (associative and non-associative) Jordan
algebras. We want to rewrite this method such that it only gives
non-associative algebras. Clearly, any annihilator extension of
non-associative algebra is non-associative. So, we only have to examine the
annihilator extensions of associative Jordan algebras. Let $J$ be an
associative Jordan algebra and $\theta \in Z^{2}\left( J,V\right) $. Then $%
J_{\theta }$ is an associative Jordan algebra if and only if $\theta \left(
x,y\circ z\right) =\theta \left( x\circ y,z\right) $ for all $x,y,z\in $$J$.
Define a subspace $\mathcal{Z}^{2}\left( J,V\right) $ of $Z^{2}\left(
J,V\right) $ by%
\begin{equation*}
\mathcal{Z}^{2}\left( J,V\right) =\left\{ \theta \in Z^{2}\left( J,V\right)
:\theta \left( x,y\circ z\right) =\theta \left( x\circ y,z\right) 
\mbox{ for
all }x,\text{$y,z\in $}J\right\} .
\end{equation*}%
Let $e_{1},\ldots ,e_{s}$ be a basis of $V$. Suppose that $\theta \left(
x,y\right) =\underset{i=1}{\overset{s}{\sum }}\theta _{i}\left( x,y\right)
e_{i}$. Then $\theta \in \mathcal{Z}^{2}\left( J,V\right) $ if and only if
all $\theta _{i}\in \mathcal{Z}^{2}\left( J,\mathbb{F}\right) $. Define $%
\mathcal{H}^{2}\left( J,V\right) =\mathcal{Z}^{2}\left( J,V\right) \big/%
B^{2}\left( J,V\right) $. Therefore, $\mathcal{H}^{2}\left( J,V\right) $ is
a subspace of $H^{2}\left( J,V\right) $.

\begin{lem}
\label{invariantR}Let $\phi \in Aut\left( J\right) $ and $\theta \in
Z^{2}\left( J,V\right) $. Then $\phi \theta \in \mathcal{Z}^{2}\left(
J,V\right) $ if and only if $\theta \in \mathcal{Z}^{2}\left( J,V\right) $.
\end{lem}

Define $R_{s}\left( J\right) =\left\{ W\in T_{s}\left( J\right) :W\in
G_{s}\left( \mathcal{H}^{2}\left( J,\mathbb{F}\right) \right) \right\} $.
Then $T_{s}\left( J\right) =R_{s}\left( J\right) \cup U_{s}\left( J\right) $
where $U_{s}\left( J\right) =T_{s}\left( J\right) -R_{s}\left( J\right) $.
By Lemma \ref{invariantT} and Lemma \ref{invariantR}, the sets $R_{s}\left(
J\right) $ and $U_{s}\left( J\right) $ are stable under the action of $%
Aut\left( J\right) $. Let $E_{R}\left( J,V\right) =\left\{ J_{\theta }\in
E\left( J,V\right) :J_{\theta }\mbox{ is associative }\right\} $. Then $%
E\left( J,V\right) =E_{R}\left( J,V\right) \cup E_{U}\left( J,V\right) $
where $E_{U}\left( J,V\right) =E\left( J,V\right) -E_{R}\left( J,V\right) $.

\begin{thm}
Let $J$ be an associative Jordan algebra.

\begin{enumerate}
\item There exists a one-to-one correspondence between the set of $Aut\left(
J\right) $-orbits on $R_{s}\left( J\right) $ and the set of isomorphism
classes of $E_{R}\left( J,V\right) $.

\item There exists a one-to-one correspondence between the set of $Aut\left(
J\right) $-orbits on $U_{s}\left( J\right) $ and the set of isomorphism
classes of $E_{U}\left( J,V\right) $.
\end{enumerate}
\end{thm}

By this theorem, we may construct all nilpotent non-associative Jordan
algebras of dimension $n$, given those algebras of dimension less than $n$,
in the following way:

\begin{enumerate}
\item For a given nilpotent Jordan algebra $J$ of dimension $n-s$, if $J$ is
non-associative then do the following:

\begin{enumerate}
\item Determine $H^{2}\left( J,\mathbb{F}\right) $, $Ann\left( J\right) $
and $Aut\left( J\right) $.

\item Determine the set of $Aut\left( J\right) $-orbits on $T_{s}\left(
J\right) $.

\item For each orbit, construct the Jordan algebra corresponding to a
representative of it.
\end{enumerate}

\item Otherwise, do the following:

\begin{enumerate}
\item Determine $\mathcal{H}^{2}\left( J,\mathbb{F}\right) ,H^{2}\left( J,%
\mathbb{F}\right) $, $Ann\left( J\right) $ and $Aut\left( J\right) $.

\item Determine the set of $Aut\left( J\right) $-orbits on $U_{s}\left(
J\right) $.

\item For each orbit, construct the Jordan algebra corresponding to a
representative of it.
\end{enumerate}
\end{enumerate}

\section{Nilpotent Jordan algebras of dimension$\boldsymbol{\leq 3}$}

\label{dim3}In this section the classification of nilpotent Jordan algebras
of dimension up to three is given. Throughout the paper we use some
notational conventions. We denote the $j$-th algebra of dimension $i$ by $%
J_{i,j}$, and by $M_{i,k}$ we denoted the $k$-th algebra of dimension $i$
which is not isomorphic to any of $J_{i,1},\ldots ,J_{i,j}$ if
characteristic $\mathbb{F}=3$. The basis elements of algebras will be
denoted by the letters $a,b,c,\ldots $. The multiplication of an algebra is
specified by giving only the non-zero products. We describe the automorphism
group of a Jordan algebra by giving the matrix of a general element. For
this we use the column convention: the action of a $\phi \in Aut\left(
J\right) $ on the $i$-th basis element of $J$ is given by the $i$-th column
of the matrix of $\phi $.

\subsection{Nilpotent Jordan algebras of dimension\textbf{\ }$\boldsymbol{1}$%
}

\label{dim1}Let $J$ be a Jordan algebra with a basis $a$. If $J$ is
nilpotent, then $Ann\left( J\right) =\left\langle a\right\rangle $.
Therefore there is only one nilpotent Jordan algebra of dimension $1$; that
is, $J_{1,1}:a^{2}=0$.

\subsection{Nilpotent Jordan algebras of dimension $\boldsymbol{2}$}

\label{dim2}First we construct nilpotent Jordan algebras that have
annihilator components. By Definition \ref{Annihilator component}, we get
the algebra $J_{2,1}=J_{1,1}\oplus J_{1,1}$. To complete the classification,
we need to classify nilpotent Jordan algebras without annihilator
components. For this, we consider $1$-dimensional annihilator extensions of $%
J_{1,1}$. We have $H^{2}\left( J_{1,1},\mathbb{F}\right) =\left\langle \left[
\delta _{a,a}\right] \right\rangle $ and $Ann\left( J_{1,1}\right) =J_{1,1}$%
. Therefore, $T_{1}\left( J_{1,1}\right) =\left\{ \left\langle \left[ \delta
_{a,a}\right] \right\rangle \right\} $. So we get only one algebra, namely $%
J_{2,2}:a^{2}=b$.

\subsection{Nilpotent Jordan algebras of dimension $\boldsymbol{3}$}

First we classify nilpotent Jordan algebras that have annihilator
components. By Definition \ref{Annihilator component}, we get the algebras $%
J_{3,1}=J_{2,1}\oplus J_{1,1},J_{3,2}=J_{2,2}\oplus J_{1,1}$. Left for us to
classify nilpotent Jordan algebras that have no annihilator components. This
will be discussed in detail in the following subsubsections.

\subsubsection{1-dimensional annihilator extensions of $J_{2,1}$}

Here $H^{2}(J_{2,1},\mathbb{F})=Sym\left( J_{2,1},\mathbb{F}\right) ,$ $%
Ann\left( J_{2,1}\right) =J_{2,1}$ and $Aut\left( J_{2,1}\right) =GL\left(
J_{2,1}\right) $. Therefore, $T_{1}\left( J_{2,1}\right) =\left\{
\left\langle \theta \right\rangle :\theta \in Sym\left( J_{2,1},\mathbb{F}%
\right) \mbox{ and }\theta ^{\bot }\cap J_{2,1}=0\right\} $. Hence $\theta $
is a nondegenerate symmetric bilinear form on $J_{2,1}$. Up to equivalence
there is only one nondegenerate symmetric bilinear form on $J_{2,1}$ (see 
\cite[Chapter V, Theorem 4]{Jacobson53}), so that $T_{1}\left(
J_{2,1}\right) =\mbox{Orb}\left( \left\langle \delta _{a,b}\right\rangle
\right) $. Therefore we get only one algebra, namely $J_{3,3}:a\circ b=c$.

\subsubsection{1-dimensional annihilator extensions of $J_{2,2}$}

Here $H^{2}(J_{2,2},\mathbb{F})=$ $\left\langle \left[ \delta _{a,b}\right]
\right\rangle $ and $Ann\left( J_{2,2}\right) =\left\langle a\right\rangle $%
. Therefore, $T_{2}\left( J_{1,1}\right) =\left\{ \left\langle \left[ \delta
_{a,b}\right] \right\rangle \right\} $. So we get only one algebra, namely $%
J_{3,4}:a^{2}=b,a\circ b=c$.

\subsubsection{2-dimensional annihilator extensions of $J_{1,1}$}

Here $H^{2}(J_{1,1},\mathbb{F})=\left\langle \left[ \delta _{a,a}\right]
\right\rangle $ and so $T_{2}\left( J_{1,1}\right) =\emptyset $.

\begin{thm}
\label{compare}Any nilpotent Jordan algebra of dimension less than four is
isomorphic to one of the following\ pairwise non-isomorphic algebras:%
\begin{table}[H] \centering%
\begin{tabular}[t]{|c|c|c|}
\hline
Algebra & Multiplication table & Comments \\ \hline
$J_{1,1}$ & \multicolumn{1}{|l|}{$\cdots $} & associative \\ \hline
$J_{2,1}$ & \multicolumn{1}{|l|}{$\cdots $} & associative \\ \hline
$J_{2,2}$ & \multicolumn{1}{|l|}{$a^{2}=b$} & associative \\ \hline
$J_{3,1}$ & \multicolumn{1}{|l|}{$\cdots $} & associative \\ \hline
$J_{3,2}$ & \multicolumn{1}{|l|}{$a^{2}=b$} & associative \\ \hline
$J_{3,3}$ & \multicolumn{1}{|l|}{$a\circ b=c$} & associative \\ \hline
$J_{3,4}$ & \multicolumn{1}{|l|}{$a^{2}=b,a\circ b=c$} & associative \\ 
\hline
\end{tabular}%
\caption{Nilpotent Jordan algebras of dimension less than four.}\label{tab
dim3}%
\end{table}%
\end{thm}

\section{Nilpotent Jordan algebras of dimension $\boldsymbol{4}$}

\label{dim4}In this section we give the list of all isomorphic types of
four-dimensional nilpotent Jordan algebras.

\subsection{Nilpotent Jordan algebras with annihilator components}

By Definition \ref{Annihilator component}, we get the algebras $%
J_{4,1}=J_{3,1}\oplus J_{1,1},J_{4,2}=J_{3,2}\oplus
J_{1,1},J_{4,3}=J_{3,3}\oplus J_{1,1},J_{4,4}=J_{3,4}\oplus J_{1,1}$.

\subsection{1-dimensional annihilator extensions of $J_{3,1}$}

Here $H^{2}\left( J_{3,1},\mathbb{F}\right) =Sym\left( J_{3,1},\mathbb{F}%
\right) ,Ann\left( J_{3,1}\right) =J_{3,1}$ and $Aut\left( J_{3,1}\right)
=GL\left( J_{3,1}\right) $. Therefore, $T_{1}\left( J_{3,1}\right) =\left\{
\left\langle \theta \right\rangle :\theta \in Sym\left( J_{3,1},\mathbb{F}%
\right) \mbox{ and }\theta ^{\bot }\cap J_{3,1}=0\right\} $. Hence $\theta $
is a nondegenerate symmetric bilinear form on $J_{3,1}$. Up to equivalence
there is only one nondegenerate symmetric bilinear form on $J_{3,1}$, so
that $T_{1}\left( J_{3,1}\right) =\mbox{Orb}\left( \left\langle \delta
_{a,b}+\delta _{c,c}\right\rangle \right) $. Therefore we get only one
algebra, namely $J_{4,5}:a\circ b=d,c^{2}=d$.

\subsection{1-dimensional annihilator extensions of $J_{3,2}$}

Here $H^{2}(J_{3,2},\mathbb{F})$ consists of $\left[ \theta \right] =\alpha
_{1}\left[ \delta _{a,b}\right] +\alpha _{2}\left[ \delta _{a,c}\right]
+\alpha _{3}\left[ \delta _{b,c}\right] +\alpha _{4}\left[ \delta _{c,c}%
\right] $ and $Ann\left( J_{3,2}\right) =\left\langle b,c\right\rangle $.
Then $\theta ^{\perp }\cap Ann\left( J_{3,2}\right) =0$ if and only if $%
\left( \alpha _{3},\alpha _{1}\alpha _{4}\right) \neq \left( 0,0\right) $
(see Example \ref{rank2}). Furthermore the automorphism group $Aut\left(
J_{3,2}\right) $\ consists of%
\begin{equation*}
\phi =\left[ 
\begin{array}{ccc}
a_{11} & 0 & 0 \\ 
a_{21} & a_{11}^{2} & a_{23} \\ 
a_{31} & 0 & a_{33}%
\end{array}%
\right] \text{ }:a_{11}a_{33}\neq 0.
\end{equation*}%
Choose an arbitrary subspace $W\in T_{1}\left( J_{3,2}\right) $. Then $%
W=\left\langle \left[ \theta \right] \mid \left( \alpha _{3},\alpha
_{1}\alpha _{4}\right) \neq \left( 0,0\right) \right\rangle $. Let $\phi =%
\big(a_{ij}\big)\in $ $Aut\left( J_{3,2}\right) $, and write $\theta =\alpha
_{1}^{\prime }\left[ \delta _{a,b}\right] +\alpha _{2}^{\prime }\left[
\delta _{a,c}\right] +\alpha _{3}^{\prime }\left[ \delta _{b,c}\right]
+\alpha _{4}^{\prime }\left[ \delta _{c,c}\right] $. Then 
\begin{eqnarray*}
\alpha _{1}^{\prime } &=&a_{11}^{3}\alpha _{1}+a_{31}a_{11}^{2}\alpha _{3},
\\
\alpha _{2}^{\prime } &=&a_{11}a_{23}\alpha _{1}+a_{11}a_{33}\alpha
_{2}+(a_{21}a_{33}+a_{31}a_{23})\alpha _{3}+a_{31}a_{33}\alpha _{4}, \\
\alpha _{3}^{\prime } &=&a_{11}^{2}a_{33}\alpha _{3}, \\
\alpha _{4}^{\prime } &=&2a_{23}a_{33}\alpha _{3}+a_{33}^{2}\alpha _{4}.
\end{eqnarray*}%
Since $a_{11}^{2}a_{33}\neq 0$, the equation for $\alpha _{3}^{\prime }$
implies that $\alpha _{3}^{\prime }\neq 0$ if and only if $\alpha _{3}\neq 0$%
. Thus $\mbox{Orb}\left( \left\langle \left[ \theta \right] :\alpha _{3}\neq
0\right\rangle \right) \cap \mbox{Orb}\left( \left\langle \left[ \theta %
\right] :\alpha _{3}=0\right\rangle \right) =\emptyset $ and hence $%
Aut(J_{3,2})$ has at least two orbits on $T_{1}\left( J_{3,2}\right) $. So
we distinguish two cases.

\begin{enumerate}
\item Suppose\ first that $\alpha _{3}\neq 0$. Let $\phi $\ be as follows:%
\begin{equation*}
\phi =%
\begin{bmatrix}
1 & 0 & 0 \\ 
\alpha _{1}\alpha _{3}^{-2}\alpha _{4}-\alpha _{2}\alpha _{3}^{-1} & 1 & -%
\frac{1}{2}\alpha _{3}^{-2}\alpha _{4} \\ 
-\alpha _{1}\alpha _{3}^{-1} & 0 & \alpha _{3}^{-1}%
\end{bmatrix}%
.
\end{equation*}%
Then $\phi W=\left\langle \left[ \delta _{b,c}\right] \right\rangle $.
Consequently, we get the algebra $J_{4,6}:a^{2}=b,b\circ c=d.$

\item Assume now that $\alpha _{3}=0$. This then implies that $\alpha
_{1}\alpha _{4}\neq 0$. Let $\phi $\ be as follows:%
\begin{equation*}
\phi =\allowbreak \allowbreak 
\begin{bmatrix}
\alpha _{1}^{-\frac{1}{3}} & 0 & 0 \\ 
0 & \alpha _{1}^{-\frac{2}{3}} & -\alpha _{1}^{-1}\alpha _{2}\alpha _{4}^{-%
\frac{1}{2}} \\ 
0 & 0 & \alpha _{4}^{-\frac{1}{2}}%
\end{bmatrix}%
.
\end{equation*}%
Then $\phi W=\left\langle \left[ \delta _{a,b}\right] +\left[ \delta _{c,c}%
\right] \right\rangle $. So we get the algebra $J_{4,7}:a^{2}=b,a\circ
b=d,c^{2}=d$.
\end{enumerate}

\subsection{1-dimensional annihilator extensions of $J_{3,3}$}

Here $H^{2}\left( J_{3,3},\mathbb{F}\right) $ consists of $\theta =\alpha
_{1}\left[ \delta _{a,a}\right] +\alpha _{2}\left[ \delta _{b,b}\right]
+\alpha _{3}\left[ \delta _{a,c}\right] +\alpha _{4}\left[ \delta _{b,c}%
\right] $ and $Ann\left( J_{3,3}\right) =\left\langle c\right\rangle $.
Furthermore the automorphism group $Aut\left( J_{3,3}\right) $\ consists of%
\begin{equation}
\phi =\left[ 
\begin{array}{ccc}
a_{11} & a_{12} & 0 \\ 
a_{21} & a_{22} & 0 \\ 
a_{31} & a_{32} & a_{11}a_{22}+a_{21}a_{12}%
\end{array}%
\right] :a_{11}^{2}a_{22}^{2}-a_{12}^{2}a_{21}^{2}\neq 0\mbox{ and }%
a_{11}a_{21}=a_{12}a_{22}=0.  \label{AutJ3,3}
\end{equation}%
Choose an arbitrary subspace $W\in T_{1}\left( J_{3,3}\right) $. Then$\
W=\left\langle \left[ \theta \right] \mid \left( \alpha _{3},\alpha
_{4}\right) \neq \left( 0,0\right) \right\rangle $. Let $\phi =\big(a_{ij}%
\big)\in $ $Aut\left( J_{3,3}\right) $. Write $\phi \theta =\alpha
_{1}^{\prime }\delta _{a,a}+\alpha _{2}^{\prime }\delta _{b,b}+\alpha
_{3}^{\prime }\delta _{a,c}+\alpha _{4}^{\prime }\delta _{b,c}$. Then 
\begin{eqnarray*}
\alpha _{1}^{^{\prime }} &=&a_{11}^{2}\alpha _{1}+a_{21}^{2}\alpha
_{2}+2a_{11}a_{31}\alpha _{3}+2a_{21}a_{31}\alpha _{4}, \\
\alpha _{2}^{^{\prime }} &=&a_{12}^{2}\alpha _{1}+a_{22}^{2}\alpha
_{2}+2a_{12}a_{32\alpha _{3}}+2a_{22}a_{32}\alpha _{4}, \\
\alpha _{3}^{^{\prime }} &=&a_{11}^{2}a_{22}\alpha
_{3}+a_{21}^{2}a_{12}\alpha _{4}, \\
\alpha _{4}^{^{\prime }} &=&a_{21}a_{12}^{2}\alpha
_{3}+a_{11}a_{22}^{2}\alpha _{4}.
\end{eqnarray*}%
From $\left( \ref{AutJ3,3}\right) $ we either have $a_{11}a_{22}\neq
0,a_{21}=a_{12}=0$ or $a_{11}=a_{22}=0,a_{21}a_{12}\neq 0$. Consequently, $%
\alpha _{3}^{\prime }\alpha _{4}^{^{\prime }}\neq 0$ if and only if $\alpha
_{3}\alpha _{4}\neq 0$. Thus $\mbox{Orb}\left( \left\langle \left[ \theta %
\right] :\alpha _{3}\alpha _{4}\neq 0\right\rangle \right) \cap \mbox{Orb}%
\left( \left\langle \left[ \theta \right] :\alpha _{3}\alpha
_{4}=0\right\rangle \right) =\emptyset $ and hence $Aut\left( J_{3,3}\right) 
$ has at least two orbits on $T_{1}\left( J_{3,3}\right) $. So we
distinguish two cases.

\begin{enumerate}
\item Suppose first that $\alpha _{3}\alpha _{4}=0$. Then $\alpha _{3}\neq 0$
and $\alpha _{4}=0$ or $\alpha _{3}=0$ and $\alpha _{4}\neq 0$. Let us
assume first that $\alpha _{3}\neq 0$ and $\alpha _{4}=0$. Choose $\phi $ as
follows:%
\begin{equation*}
\phi =%
\begin{bmatrix}
\epsilon ^{\frac{1}{4}}\alpha _{3}^{-\frac{1}{2}} & 0 & 0 \\ 
0 & \epsilon ^{-\frac{1}{2}} & 0 \\ 
-\frac{1}{2}\epsilon ^{\frac{1}{4}}\alpha _{1}\alpha _{3}^{-\frac{3}{2}} & 0
& \epsilon ^{-\frac{1}{4}}\alpha _{3}^{-\frac{1}{2}}%
\end{bmatrix}%
\end{equation*}%
where $\epsilon =1$ if $\alpha _{2}=0$, and otherwise $\epsilon =\alpha _{2}$%
. Then $\phi W=\left\langle \left[ \delta _{a,c}\right] \right\rangle $ if $%
\alpha _{2}=0$, while $\phi W=\left\langle \left[ \delta _{a,c}\right] +%
\left[ \delta _{b,b}\right] \right\rangle $ otherwise. So we have two
representatives, namely $W_{1}=\left\langle \delta _{a,c}\right\rangle $ and 
$W_{2}=\left\langle \left[ \delta _{a,c}\right] +\left[ \delta _{b,b}\right]
\right\rangle $. Assume now that $\alpha _{3}=0$ and $\alpha _{4}\neq 0$.
Choose $\phi $ as follows:%
\begin{equation*}
\phi =%
\begin{bmatrix}
0 & \varepsilon ^{-\frac{1}{2}} & 0 \\ 
\varepsilon ^{\frac{1}{4}}\alpha _{4}^{-\frac{1}{2}} & 0 & 0 \\ 
-\frac{1}{2}\varepsilon ^{\frac{1}{4}}\alpha _{2}\alpha _{4}^{-\frac{3}{2}}
& 0 & \varepsilon ^{-\frac{1}{4}}\alpha _{4}^{-\frac{1}{2}}%
\end{bmatrix}%
\end{equation*}%
where $\varepsilon =1$ if $\alpha _{1}=0$, and otherwise $\varepsilon
=\alpha _{1}$. Then $\phi W=W_{1}$ if $\alpha _{1}=0$, while $\phi W=W_{2}$
otherwise. Therefore, we have the same representatives. Furthermore the
subspaces $W_{1},W_{2}$ are not in the same orbit since $\Psi \left(
W_{1}\right) =\left( 1\right) $ and $\Psi \left( W_{2}\right) =\left(
0\right) $. Thus we get the following pairwise non-isomorphic algebras:

\begin{itemize}
\item $J_{4,8}:a\circ b=c,a\circ c=d.$

\item $J_{4,9}:a\circ b=c,a\circ c=d,b^{2}=d.$
\end{itemize}

\item Suppose now that $\alpha _{3}\alpha _{4}\neq 0$. Let $\phi $ be as
follows:%
\begin{equation*}
\phi =%
\begin{bmatrix}
1 & 0 & 0 \\ 
0 & \alpha _{3}\alpha _{4}^{-1} & 0 \\ 
-\frac{1}{2}\alpha _{1}\alpha _{3}^{-1} & -\frac{1}{2}\alpha _{2}\alpha
_{3}\alpha _{4}^{-2} & \alpha _{3}\alpha _{4}^{-1}%
\end{bmatrix}%
.
\end{equation*}%
Then $\phi W=\left\langle \left[ \delta _{a,c}\right] +\left[ \delta _{b,c}%
\right] \right\rangle $. So we get the algebra $J_{4,10}:a\circ b=c,a\circ
c=d,b\circ c=d$.
\end{enumerate}

\subsection{1-dimensional annihilator extensions of\textbf{\ }$J_{3,4}$}

Here $H^{2}(J_{3,4},\mathbb{F})$ is spanned by $\left[ \delta _{a,c}\right] +%
\left[ \delta _{b,b}\right] $ and hence $T_{1}\left( J_{3,4}\right) =\left\{
\left\langle \left[ \delta _{a,c}\right] +\left[ \delta _{b,b}\right]
\right\rangle \right\} $. Therefore we get only one algebra, namely $%
J_{4,11}:a^{2}=b,a\circ b=c,a\circ c=d,b^{2}=d$.

\subsection{2-dimensional annihilator extensions of $J_{2,1}$}

Here $Aut\left( J_{2,1}\right) =GL\left( J_{2,1}\right) $ and $H^{2}(J_{2,1},%
\mathbb{F})=Sym\left( J_{2,1},\mathbb{F}\right) $. Consider any arbitrary
subspace $W\in T_{2}\left( J_{2,1}\right) $. Then $\Psi \left( W\right) \neq
\left( 0,0\right) $. Indeed, if $W=\left\langle \theta _{1},\theta
_{2}\right\rangle $ with $\dim \theta _{1}^{\perp }=\dim \theta _{2}^{\bot
}=0$, then there is a $\lambda \in \mathbb{F}^{\ast }$ such that $\dim
\left( \theta _{1}+\lambda \theta _{2}\right) ^{\perp }=1$. Thus $\Psi
\left( W\right) \in \left\{ \left( 1,0\right) ,\left( 1,1\right) \right\} $.
Let us consider two cases.

\begin{enumerate}
\item Let $\Psi \left( W\right) =\left( 1,0\right) $. Then there is a basis $%
\theta _{1},\theta _{2}$ of $W$ such that $\dim \theta _{1}^{\perp }=1$ and $%
\dim \theta _{2}^{\bot }=0$. Therefore, $\theta _{2}$ is a nondegenerate
symmetric bilinear form on $J_{2,1}$. Up to equivalence there is only one
nondegenerate symmetric bilinear form on $J_{2,1}$. So there is a $\phi \in
GL\left( J_{2,1}\right) $ such that $\phi W=\left\langle \phi \theta
_{1},\delta _{a,b}\right\rangle $. Let $\phi \theta _{1}=\alpha _{1}^{\prime
}\delta _{a,a}+\alpha _{2}^{\prime }\delta _{a,b}+\alpha _{3}^{\prime
}\delta _{b,b}$. Then $\alpha _{2}^{\prime 2}=\alpha _{1}^{\prime }\alpha
_{3}^{\prime }=0$; otherwise, $\Psi \left( W\right) =\left( 1,1\right) $.
Consequently, $\phi W=\left\langle \delta _{a,a},\delta _{a,b}\right\rangle $
if $\alpha _{1}^{\prime }\neq 0$, while $\phi W=\left\langle \delta
_{b,b},\delta _{a,b}\right\rangle $ otherwise. So we have two
representatives, namely $\left\langle \delta _{a,a},\delta
_{a,b}\right\rangle $ and $\left\langle \delta _{b,b},\delta
_{a,b}\right\rangle $. Furthermore, $\mbox{Orb}\left( \left\langle \delta
_{a,a},\delta _{a,b}\right\rangle \right) =\mbox{Orb}\left( \left\langle
\delta _{b,b},\delta _{a,b}\right\rangle \right) $. To illustrate this, let $%
\phi $ be as follows:%
\begin{equation*}
\phi =%
\begin{bmatrix}
0 & 1 \\ 
1 & 0%
\end{bmatrix}%
.
\end{equation*}%
Then $\phi \left\langle \delta _{b,b},\delta _{a,b}\right\rangle
=\left\langle \delta _{a,a},\delta _{a,b}\right\rangle $. So we get the
algebra $J_{4,12}:a^{2}=c,a\circ b=d$.

\item Let $\Psi \left( W\right) =\left( 1,1\right) $. Then there is a basis $%
\theta _{1},\theta _{2}$ of $W$ such that $\dim \theta _{1}^{\perp }=\dim
\theta _{2}^{\bot }=1$ and $\theta _{1}^{\perp }\cap \theta _{2}^{\bot }=0$.
Let $\theta _{1}=\alpha _{1}\delta _{a,a}+\alpha _{2}\delta _{a,b}+\alpha
_{3}\delta _{b,b}$ and $\theta _{2}=\beta _{1}\delta _{a,a}+\beta _{2}\delta
_{a,b}+\beta _{3}\delta _{b,b}$. As $\dim \theta _{1}^{\perp }=\dim \theta
_{2}^{\bot }=1$ and $\theta _{1}^{\perp }\cap \theta _{2}^{\bot }=0$ we have 
$\alpha _{2}^{2}=\alpha _{1}\alpha _{3},\beta _{2}^{2}=\beta _{1}\beta _{3}$
and $\alpha _{1}\beta _{3}-2\alpha _{2}\beta _{2}+\beta _{1}\alpha _{3}\neq
0 $. Then $\alpha _{2}=\epsilon \sqrt{\alpha _{1}}\sqrt{\alpha _{3}},\beta
_{2}=\varepsilon \sqrt{\beta _{1}}\sqrt{\beta _{3}}$ and $\epsilon \sqrt{%
\beta _{1}}\sqrt{\alpha _{3}}-\varepsilon \sqrt{\alpha _{1}}\sqrt{\beta _{3}}%
\neq 0$ where $\epsilon ^{2}=\varepsilon ^{2}=1$. Let $\phi $ be the
following automorphism:%
\begin{equation*}
\phi =\left( \epsilon \sqrt{\beta _{1}}\sqrt{\alpha _{3}}-\varepsilon \sqrt{%
\alpha _{1}}\sqrt{\beta _{3}}\right) ^{-1}%
\begin{bmatrix}
-\varepsilon \sqrt{\beta _{3}} & \epsilon \sqrt{\alpha _{3}} \\ 
\sqrt{\beta _{1}} & -\sqrt{\alpha _{1}}%
\end{bmatrix}%
.
\end{equation*}%
Then $\phi W=\left\langle \delta _{a,a},\delta _{b,b}\right\rangle $. Thus
we get only one algebra, namely $J_{4,13}:a^{2}=c,b^{2}=d$.
\end{enumerate}

\subsection{2-dimensional annihilator extensions of $J_{2,2}$}

Since $H^{2}(J_{2,2},\mathbb{F})=\left\langle \left[ \delta _{b,b}\right]
\right\rangle $, it follows that $T_{2}\left( J_{2,2}\right) =\emptyset $.

\begin{thm}
Any nilpotent Jordan algebra of dimension four is isomorphic to one of the
following\ pairwise non-isomorphic algebras:%
\begin{table}[H] \centering%
$%
\begin{tabular}[t]{|c|c|c|}
\hline
Algebra & Multiplication table & Comments \\ \hline
$J_{4,1}$ & \multicolumn{1}{|l|}{$\cdots $} & associative \\ \hline
$J_{4,2}$ & \multicolumn{1}{|l|}{$a^{2}=b$} & associative \\ \hline
$J_{4,3}$ & \multicolumn{1}{|l|}{$a\circ b=c$} & associative \\ \hline
$J_{4,4}$ & \multicolumn{1}{|l|}{$a^{2}=b,a\circ b=c$} & associative \\ 
\hline
$J_{4,5}$ & \multicolumn{1}{|l|}{$a\circ b=d,c^{2}=d$} & associative \\ 
\hline
$J_{4,6}$ & \multicolumn{1}{|l|}{$a^{2}=b,b\circ c=d$} & non-associative \\ 
\hline
$J_{4,7}$ & \multicolumn{1}{|l|}{$a^{2}=b,a\circ b=d,c^{2}=d$} & associative
\\ \hline
$J_{4,8}$ & \multicolumn{1}{|l|}{$a\circ b=c,a\circ c=d$} & non-associative
\\ \hline
$J_{4,9}$ & \multicolumn{1}{|l|}{$a\circ b=c,a\circ c=d,b^{2}=d$} & 
non-associative \\ \hline
$J_{4,10}$ & \multicolumn{1}{|l|}{$a\circ b=c,a\circ c=d,b\circ c=d$} & 
non-associative \\ \hline
$J_{4,11}$ & \multicolumn{1}{|l|}{$a^{2}=b,a\circ b=c,a\circ c=d,b^{2}=d$} & 
associative \\ \hline
$J_{4,12}$ & \multicolumn{1}{|l|}{$a^{2}=c,a\circ b=d$} & associative \\ 
\hline
$J_{4,13}$ & \multicolumn{1}{|l|}{$a^{2}=c,b^{2}=d$} & associative \\ \hline
\end{tabular}%
$%
\caption{Nilpotent Jordan algebras of dimension four.}%
\label{tab2}%
\end{table}%
\end{thm}

\begin{rem}
The classification of nilpotent Jordan algebras of dimension up to four over 
$%
\mathbb{C}
$ given in \cite{Ancochea} is the same as the one given here over an
algebraically closed field $\mathbb{F}$ of characteristic $\neq 2$.
\end{rem}

\section{Non-associative nilpotent Jordan algebras of dimension 5}

\label{dim5}In this section we only classify non-associative nilpotent
Jordan algebras of dimension $5$ since all associative nilpotent Jordan
algebras of dimension $5$ have been classified in \cite{Ponnen}. For this,
we need to describe $\mathcal{H}^{2}(J,\mathbb{F}),H^{2}(J,\mathbb{F})$ and $%
Ann\left( J\right) $ for all nilpotent Jordan algebras of dimension $\leq 4$%
. So we have the following table.%
\begin{table}[H] \centering%
$%
\begin{tabular}[t]{|c|c|c|c|c|}
\hline
$J$ & $\mathcal{H}^{2}(J,\mathbb{F})$ & $H^{2}(J,\mathbb{F})$ & $Ann\left(
J\right) $ & Comments \\ \hline
$J_{1,1}$ & $Sym\left( J_{1,1},\mathbb{F}\right) $ & $Sym\left( J_{1,1},%
\mathbb{F}\right) $ & $\left\langle a\right\rangle $ & $U_{4}\left(
J_{1,1}\right) =\emptyset $ \\ \hline
$J_{2,1}$ & $Sym\left( J_{2,1},\mathbb{F}\right) $ & $Sym\left( J_{2,1},%
\mathbb{F}\right) $ & $\left\langle a,b\right\rangle $ & $U_{3}\left(
J_{2,1}\right) =\emptyset $ \\ \hline
$J_{2,2}$ & $\left\langle \left[ \delta _{a,b}\right] \right\rangle $ & $%
\mathcal{H}^{2}(J_{2,2},\mathbb{F})$ & $\left\langle b\right\rangle $ & $%
U_{3}\left( J_{2,2}\right) =\emptyset $ \\ \hline
$J_{3,1}$ & $Sym\left( J_{3,1},\mathbb{F}\right) $ & $Sym\left( J_{3,1},%
\mathbb{F}\right) $ & $\left\langle a,b,c\right\rangle $ & $U_{2}\left(
J_{3,1}\right) =\emptyset $ \\ \hline
$J_{3,2}$ & $\left\langle \left[ \delta _{a,b}\right] ,\left[ \delta _{a,c}%
\right] ,\left[ \delta _{c,c}\right] \right\rangle $ & $\mathcal{H}%
^{2}(J_{3,2},\mathbb{F})\oplus \left\langle \left[ \delta _{b,c}\right]
\right\rangle $ & $\left\langle b,c\right\rangle $ & $U_{2}\left(
J_{3,1}\right) \neq \emptyset $ \\ \hline
$J_{3,3}$ & $\left\langle \left[ \delta _{a,a}\right] ,\left[ \delta _{b,b}%
\right] \right\rangle $ & $\mathcal{H}^{2}(J_{3,3},\mathbb{F})\oplus
\left\langle \left[ \delta _{a,c}\right] ,\left[ \delta _{b,c}\right]
\right\rangle $ & $\left\langle c\right\rangle $ & $U_{2}\left(
J_{3,2}\right) \neq \emptyset $ \\ \hline
$J_{3,4}$ & $\left\langle \left[ \delta _{a,c}\right] +\left[ \delta _{b,b}%
\right] \right\rangle $ & $\mathcal{H}^{2}(J_{3,4},\mathbb{F})$ & $%
\left\langle c\right\rangle $ & $U_{2}\left( J_{3,4}\right) =\emptyset $ \\ 
\hline
$J_{4,1}$ & $Sym\left( J_{4,1},\mathbb{F}\right) $ & $Sym\left( J_{4,1},%
\mathbb{F}\right) $ & $\left\langle a,b,c,d\right\rangle $ & $U_{1}\left(
J_{4,1}\right) =\emptyset $ \\ \hline
$J_{4,2}$ & $\left\langle \left[ \delta _{c,c}\right] ,\left[ \delta _{d,d}%
\right] ,\left[ \delta _{a,b}\right] ,\left[ \delta _{a,c}\right] ,\left[
\delta _{a,d}\right] ,\left[ \delta _{c,d}\right] \right\rangle $ & $%
\mathcal{H}^{2}\left( J_{4,2},\mathbb{F}\right) \oplus \left\langle \left[
\delta _{b,c}\right] ,\left[ \delta _{b,d}\right] \right\rangle $ & $%
\left\langle b,c,d\right\rangle $ & $U_{1}\left( J_{4,2}\right) \neq
\emptyset $ \\ \hline
$J_{4,3}$ & $\left\langle \left[ \delta _{a,a}\right] ,\left[ \delta _{b,b}%
\right] ,\left[ \delta _{d,d}\right] ,\left[ \delta _{a,d}\right] ,\left[
\delta _{b,d}\right] \right\rangle $ & $\mathcal{H}^{2}\left( J_{4,3},%
\mathbb{F}\right) \oplus \left\langle \left[ \delta _{a,c}\right] ,\left[
\delta _{b,c}\right] ,\left[ \delta _{c,d}\right] \right\rangle $ & $%
\left\langle c,d\right\rangle $ & $U_{1}\left( J_{4,3}\right) \neq \emptyset 
$ \\ \hline
$J_{4,4}$ & $\left\langle \left[ \delta _{a,c}\right] +\left[ \delta _{b,b}%
\right] ,\left[ \delta _{a,d}\right] ,\left[ \delta _{d,d}\right]
\right\rangle $ & $\mathcal{H}^{2}\left( J_{4,4},\mathbb{F}\right) \oplus
\left\langle \left[ \delta _{b,d}\right] \right\rangle $ & $\left\langle
c,d\right\rangle $ & $U_{1}\left( J_{4,4}\right) \neq \emptyset $ \\ \hline
$J_{4,5}$ & $\left\langle \left[ \delta _{a,a}\right] ,\left[ \delta _{b,b}%
\right] ,\left[ \delta _{a,b}\right] ,\left[ \delta _{a,c}\right] ,\left[
\delta _{b,c}\right] \right\rangle $ & $\mathcal{H}^{2}\left( J_{4,5},%
\mathbb{F}\right) \oplus \left\langle \left[ \delta _{a,d}\right] ,\left[
\delta _{b,d}\right] ,\left[ \delta _{c,d}\right] \right\rangle $ & $%
\left\langle d\right\rangle $ & $U_{1}\left( J_{4,5}\right) \neq \emptyset $
\\ \hline
$J_{4,6}$ & non-associative & $\left\langle \left[ \delta _{a,b}\right] ,%
\left[ \delta _{a,c}\right] ,\left[ \delta _{c,c}\right] \right\rangle $ & $%
\left\langle d\right\rangle $ & $T_{1}\left( J_{4,6}\right) =\emptyset $ \\ 
\hline
$J_{4,7}$ & $\left\langle \left[ \delta _{a,c}\right] ,\left[ \delta _{a,b}%
\right] \right\rangle $ & $\left\langle \left[ \delta _{a,c}\right] ,\left[
\delta _{a,b}\right] ,\left[ \delta _{b,b}\right] +\left[ \delta _{a,d}%
\right] \right\rangle $ & $\left\langle d\right\rangle $ & $U_{1}\left(
J_{4,7}\right) \neq \emptyset $ \\ \hline
$J_{4,8}$ & non-associative & $\left\langle \left[ \delta _{a,a}\right] ,%
\left[ \delta _{b,b}\right] ,\left[ \delta _{b,c}\right] \right\rangle $ & $%
\left\langle d\right\rangle $ & $T_{1}\left( J_{4,8}\right) =\emptyset $ \\ 
\hline
$J_{4,9}$ & non-associative & $\left\langle \left[ \delta _{a,a}\right] ,%
\left[ \delta _{a,c}\right] ,\left[ \delta _{b,c}\right] \right\rangle $ & $%
\left\langle d\right\rangle $ & $T_{1}\left( J_{4,9}\right) =\emptyset $ \\ 
\hline
$J_{4,10}$ & non-associative & $\left\langle \left[ \delta _{a,a}\right] ,%
\left[ \delta _{b,b}\right] ,\left[ \delta _{a,c}\right] \right\rangle $ & $%
\left\langle d\right\rangle $ & $T_{1}\left( J_{4,10}\right) =\emptyset $ \\ 
\hline
$J_{4,11}$ & $\left\langle \left[ \delta _{a,d}\right] +\left[ \delta _{b,b}%
\right] \right\rangle $ & $\mathcal{H}^{2}(J_{4,11},\mathbb{F})$ & $%
\left\langle d\right\rangle $ & $U_{1}\left( J_{4,11}\right) \neq \emptyset $
\\ \hline
$J_{4,12}$ & $\left\langle \left[ \delta _{a,c}\right] ,\left[ \delta _{b,b}%
\right] ,\left[ \delta _{a,d}\right] +\left[ \delta _{b,c}\right]
\right\rangle $ & $\left\langle \left[ \delta _{a,c}\right] ,\left[ \delta
_{b,b}\right] ,\left[ \delta _{a,d}\right] ,\left[ \delta _{b,c}\right] ,%
\left[ \delta _{b,d}\right] \right\rangle $ & $\left\langle c,d\right\rangle 
$ & $U_{1}\left( J_{4,12}\right) \neq \emptyset $ \\ \hline
$J_{4,13}$ & $\left\langle \left[ \delta _{a,b}\right] ,\left[ \delta _{a,c}%
\right] ,\left[ \delta _{b,d}\right] \right\rangle $ & $\mathcal{H}%
^{2}(J_{4,11},\mathbb{F})\oplus \left\langle \left[ \delta _{b,c}\right] ,%
\left[ \delta _{a,d}\right] \right\rangle $ & $\left\langle c,d\right\rangle 
$ & $U_{1}\left( J_{4,13}\right) \neq \emptyset $ \\ \hline
\end{tabular}%
$%
\caption{Description of $\mathcal{H}^{2}\left( J,\mathbb{F}\right), H^{2}\left( J,\mathbb{F}\right) $ and $Ann\left( J\right) $.}%
\label{tabx}%
\end{table}%
From this table we have the following lemma.

\begin{lem}
Any $5$-dimensional nilpotent non-associative Jordan algebra without
annihilator components is one of the following:

\begin{enumerate}
\item 1-dimensional annihilator extension of $J_{4,2}$.

\item 1-dimensional annihilator extension of $J_{4,3}$.

\item 1-dimensional annihilator extension of $J_{4,4}$.

\item 1-dimensional annihilator extension of $J_{4,5}$.

\item 1-dimensional annihilator extension of $J_{4,7}$.

\item 1-dimensional annihilator extension of $J_{4,12}$.

\item 1-dimensional annihilator extension of $J_{4,13}$.

\item 2-dimensional annihilator extension of $J_{3,2}$.

\item 2-dimensional annihilator extension of $J_{3,3}$.
\end{enumerate}
\end{lem}

Let us now classify $5$-dimensional nilpotent non-associative Jordan
algebras. For this, we have the following subsections.

\subsection{Nilpotent Jordan algebras with annihilator components}

Let $J$ be a nilpotent non-associative Jordan algebra with annihilator
component. Then by Definition \ref{Annihilator component}, $J=I\oplus
J_{1,1} $. Moreover, $J$ is non-associative if and only if $I$ is
non-associative. So, from Table \ref{tab2}, $I$ is one of $%
J_{4,6},J_{4,8},J_{4,9},J_{4,10}$. Therefore we get the algebras:

\begin{itemize}
\item $J_{5,1}=J_{4,6}\oplus J_{1,1}:a^{2}=b$ $,$ $b\circ c=d.$

\item $J_{5,2}=J_{4,8}\oplus J_{1,1}:a\circ b=c,a\circ c=d.$

\item $J_{5,3}=J_{4,9}\oplus J_{1,1}:a\circ b=c,a\circ c=d,b^{2}=d.$

\item $J_{5,4}=J_{4,10}\oplus J_{1,1}:a\circ b=c,a\circ c=d,b\circ c=d.$
\end{itemize}

\subsection{1-dimensional annihilator extensions of $J_{4,2}$}

The automorphism group $Aut(J_{4,2})$ consists of

\begin{equation*}
\phi =\left[ 
\begin{array}{cccc}
a_{11} & 0 & 0 & 0 \\ 
a_{21} & a_{11}^{2} & a_{23} & a_{24} \\ 
a_{31} & 0 & a_{33} & a_{34} \\ 
a_{41} & 0 & a_{43} & a_{44}%
\end{array}%
\right] :a_{11}\left( a_{33}a_{44}-a_{34}a_{43}\right) \neq 0.
\end{equation*}%
Choose an arbitrary subspace $W\in U_{1}\left( J_{4,2}\right) $. From Table %
\ref{tabx}, such a subspace is spanned by $\left[ \theta \right] =\alpha _{1}%
\left[ \delta _{c,c}\right] +\alpha _{2}\left[ \delta _{d,d}\right] +\alpha
_{3}\left[ \delta _{a,b}\right] +\alpha _{4}\left[ \delta _{a,c}\right]
+\alpha _{5}\left[ \delta _{a,d}\right] +\alpha _{6}\left[ \delta _{c,d}%
\right] +\alpha _{7}\left[ \delta _{b,c}\right] +\alpha _{8}\left[ \delta
_{b,d}\right] $ such that $\left( \alpha _{7},\alpha _{8}\right) \neq \left(
0,0\right) $ and $\theta ^{\bot }\cap \left\langle b,c,d\right\rangle =0$.
Further, $\theta ^{\bot }\cap \left\langle b,c,d\right\rangle =0$ if and
only if the matrix%
\begin{equation*}
A=%
\begin{bmatrix}
\alpha _{3} & \alpha _{4} & \alpha _{5} \\ 
0 & \alpha _{7} & \alpha _{8} \\ 
\alpha _{7} & \alpha _{1} & \alpha _{6} \\ 
\alpha _{8} & \alpha _{6} & \alpha _{2}%
\end{bmatrix}%
\end{equation*}%
has rank $3$. Hence $\theta ^{\bot }\cap \left\langle b,c,d\right\rangle =0$
if and only if $\left( \epsilon _{1},\epsilon _{2},\epsilon _{3},\epsilon
_{4}\right) \neq \left( 0,0,0,0\right) $ where%
\begin{eqnarray*}
\epsilon _{1} &=&\alpha _{2}\alpha _{7}^{2}-2\alpha _{6}\alpha _{7}\alpha
_{8}+\alpha _{1}\alpha _{8}^{2}, \\
\epsilon _{2} &=&\alpha _{5}\alpha _{7}^{2}+\alpha _{1}\alpha _{3}\alpha
_{8}-\alpha _{3}\alpha _{6}\alpha _{7}-\alpha _{4}\alpha _{7}\alpha _{8}, \\
\epsilon _{3} &=&\alpha _{3}\alpha _{6}\alpha _{8}-\alpha _{4}\alpha
_{8}^{2}+\alpha _{5}\alpha _{7}\alpha _{8}-\alpha _{2}\alpha _{3}\alpha _{7},
\\
\epsilon _{4} &=&\alpha _{1}\alpha _{2}\alpha _{3}-\alpha _{3}\alpha
_{6}^{2}-\alpha _{2}\alpha _{4}\alpha _{7}-\alpha _{1}\alpha _{5}\alpha
_{8}+\alpha _{4}\alpha _{6}\allowbreak \alpha _{8}+\alpha _{5}\alpha
_{6}\alpha _{7}.
\end{eqnarray*}%
Let $\phi =\big(a_{ij}\big)\in $ $Aut\left( J_{4,2}\right) $. Write $\left[
\phi \theta \right] =\alpha _{1}^{\prime }\left[ \delta _{c,c}\right]
+\alpha _{2}^{\prime }\left[ \delta _{d,d}\right] +\alpha _{3}^{\prime }%
\left[ \delta _{a,b}\right] +\alpha _{4}^{\prime }\left[ \delta _{a,c}\right]
+\alpha _{5}^{\prime }\left[ \delta _{a,d}\right] +\alpha _{6}^{\prime }%
\left[ \delta _{c,d}\right] +\alpha _{7}^{\prime }\left[ \delta _{b,c}\right]
+\alpha _{8}^{\prime }\left[ \delta _{b,d}\right] $. Then%
\begin{eqnarray*}
\alpha _{1}^{\prime } &=&a_{33}\left( \alpha _{1}a_{33}+\alpha
_{7}a_{23}+\alpha _{6}a_{43}\right) +a_{43}\left( \alpha _{2}a_{43}+\alpha
_{6}a_{33}+\alpha _{8}a_{23}\right) +a_{23}\left( \alpha _{7}a_{33}+\alpha
_{8}a_{43}\right) , \\
\alpha _{2}^{\prime } &=&a_{34}\left( \alpha _{1}a_{34}+\alpha
_{7}a_{24}+\alpha _{6}a_{44}\right) +a_{44}\left( \alpha _{2}a_{44}+\alpha
_{6}a_{34}+\alpha _{8}a_{24}\right) +a_{24}\left( \alpha _{7}a_{34}+\alpha
_{8}a_{44}\right) , \\
\alpha _{3}^{\prime } &=&a_{11}^{2}\left( \alpha _{3}a_{11}+\alpha
_{7}a_{31}+\alpha _{8}a_{41}\right) , \\
\alpha _{4}^{\prime } &=&a_{33}\left( \alpha _{1}a_{31}+\alpha
_{4}a_{11}+\alpha _{7}a_{21}+\alpha _{6}a_{41}\right) +a_{43}\left( \alpha
_{2}a_{41}+\alpha _{5}a_{11}+\alpha _{6}a_{31}+\alpha _{8}a_{21}\right)
+a_{23}\left( \alpha _{3}a_{11}+\alpha _{7}a_{31}+\alpha _{8}a_{41}\right) ,
\\
\alpha _{5}^{\prime } &=&a_{34}\left( \alpha _{1}a_{31}+\alpha
_{4}a_{11}+\alpha _{7}a_{21}+\alpha _{6}a_{41}\right) +a_{44}\left( \alpha
_{2}a_{41}+\alpha _{5}a_{11}+\alpha _{6}a_{31}+\alpha _{8}a_{21}\right)
+a_{24}\left( \alpha _{3}a_{11}+\alpha _{7}a_{31}+\alpha _{8}a_{41}\right) ,
\\
\alpha _{6}^{\prime } &=&a_{34}\left( \alpha _{1}a_{33}+\alpha
_{7}a_{23}+\alpha _{6}a_{43}\right) +a_{44}\left( \alpha _{2}a_{43}+\alpha
_{6}a_{33}+\alpha _{8}a_{23}\right) +a_{24}\left( \alpha _{7}a_{33}+\alpha
_{8}a_{43}\right) , \\
\alpha _{7}^{\prime } &=&a_{11}^{2}\left( \alpha _{7}a_{33}+\alpha
_{8}a_{43}\right) , \\
\alpha _{8}^{\prime } &=&a_{11}^{2}\left( \alpha _{7}a_{34}+\alpha
_{8}a_{44}\right) .
\end{eqnarray*}%
Let us now consider two cases:

\begin{enumerate}
\item Suppose first that $\alpha _{7}=0$. This then implies that $\alpha
_{8}\neq 0$. Moreover,%
\begin{equation*}
\epsilon _{1}=\alpha _{1}\alpha _{8}^{2},\epsilon _{2}=\alpha _{3}\alpha
_{8}^{-1}\epsilon _{1},\epsilon _{3}=\alpha _{8}\left( \alpha _{3}\alpha
_{6}-\alpha _{4}\alpha _{8}\right) ,\epsilon _{4}=\alpha _{8}^{-2}\left(
\alpha _{2}\alpha _{3}-\alpha _{5}\alpha _{8}\right) \epsilon _{1}-\alpha
_{6}\alpha _{8}^{-1}\epsilon _{3}.
\end{equation*}%
This shows that the matrix $A$ has rank $3$ if and only if $\left( \epsilon
_{1},\epsilon _{3}\right) \neq \left( 0,0\right) $ or, equivalently, $\left(
\alpha _{1},\alpha _{3}\alpha _{6}-\alpha _{4}\alpha _{8}\right) \neq \left(
0,0\right) $.

\begin{enumerate}
\item Suppose that $\alpha _{1}\neq 0$. Let $\phi $\ be the following
automorphism:%
\begin{equation*}
\phi =\left[ 
\begin{array}{cccc}
1 & 0 & 0 & 0 \\ 
\alpha _{1}^{-1}\alpha _{8}^{-2}\epsilon _{4} & 1 & -\frac{1}{2}\alpha
_{2}\alpha _{8}^{-2} & -\alpha _{1}^{-\frac{1}{2}}\alpha _{6}\alpha _{8}^{-1}
\\ 
\alpha _{1}^{-1}\alpha _{8}^{-1}\left( \alpha _{3}\alpha _{6}-\alpha
_{4}\alpha _{8}\right) & 0 & 0 & \alpha _{1}^{-\frac{1}{2}} \\ 
-\alpha _{3}\alpha _{8}^{-1} & 0 & \alpha _{8}^{-1} & 0%
\end{array}%
\right] .
\end{equation*}%
Then $\phi W=\left\langle \left[ \delta _{d,d}\right] +\left[ \delta _{b,c}%
\right] \right\rangle $. Hence we get a representative $W_{1}=\left\langle %
\left[ \delta _{d,d}\right] +\left[ \delta _{b,c}\right] \right\rangle $.

\item Suppose that $\alpha _{1}=0$. So $\alpha _{3}\alpha _{6}-\alpha
_{4}\alpha _{8}\neq 0$. Let $\phi $\ be the following automorphism:%
\begin{equation*}
\phi =\left[ 
\begin{array}{cccc}
1 & 0 & 0 & 0 \\ 
\alpha _{8}^{-2}\left( \alpha _{2}\alpha _{3}-\alpha _{5}\alpha _{8}\right)
& 1 & -\frac{1}{2}\alpha _{2}\alpha _{8}^{-2} & \alpha _{6}\left( \alpha
_{3}\alpha _{6}-\alpha _{4}\alpha _{8}\right) ^{-1} \\ 
0 & 0 & 0 & -\alpha _{8}\left( \alpha _{3}\alpha _{6}-\alpha _{4}\alpha
_{8}\right) ^{-1} \\ 
-\alpha _{3}\alpha _{8}^{-1} & 0 & \alpha _{8}^{-1} & 0%
\end{array}%
\right] .
\end{equation*}%
Then $\phi W=\left\langle \left[ \delta _{a,d}\right] +\left[ \delta _{b,c}%
\right] \right\rangle $. Hence we get a representative $W_{2}=\left\langle %
\left[ \delta _{a,d}\right] +\left[ \delta _{b,c}\right] \right\rangle $.
\end{enumerate}

\item Assume now that $\alpha _{7}\neq 0$. Then%
\begin{eqnarray*}
\epsilon _{1} &=&\alpha _{2}\alpha _{7}^{2}-2\alpha _{6}\alpha _{7}\alpha
_{8}+\alpha _{1}\alpha _{8}^{2}, \\
\epsilon _{2} &=&\alpha _{5}\alpha _{7}^{2}+\alpha _{1}\alpha _{3}\alpha
_{8}-\alpha _{3}\alpha _{6}\alpha _{7}-\alpha _{4}\alpha _{7}\alpha _{8}, \\
\epsilon _{3} &=&\alpha _{7}^{-1}\left( \alpha _{8}\epsilon _{2}-\alpha
_{3}\epsilon _{1}\right) , \\
\epsilon _{4} &=&\alpha _{7}^{-2}\left( \alpha _{1}\alpha _{3}\epsilon
_{1}-\alpha _{1}\epsilon _{2}\alpha _{8}-\alpha _{4}\epsilon _{1}\alpha
_{7}+\alpha _{6}\epsilon _{2}\alpha _{7}\right) .
\end{eqnarray*}%
So the matrix $A$ has rank $3$ if and only if $\left( \epsilon _{1},\epsilon
_{2}\right) \neq \left( 0,0\right) $.

\begin{enumerate}
\item Suppose first that $\epsilon _{1}\neq 0$. Let $\phi $\ be the
following automorphism:%
\begin{equation*}
\phi =\left[ 
\begin{array}{cccc}
1 & 0 & 0 & 0 \\ 
-\alpha _{7}^{-2}\epsilon _{1}^{-1}\epsilon _{4} & 1 & -\frac{1}{2}\alpha
_{1}\alpha _{7}^{-2} & \epsilon _{1}^{-\frac{1}{2}}\alpha _{7}^{-1}\left(
\alpha _{1}\alpha _{8}-\alpha _{6}\alpha _{7}\right) \\ 
\epsilon _{1}^{-1}\alpha _{7}^{-1}\left( \alpha _{8}\epsilon _{2}-\alpha
_{3}\epsilon _{1}\right) & 0 & \alpha _{7}^{-1} & -\epsilon _{1}^{-\frac{1}{2%
}}\alpha _{8} \\ 
-\epsilon _{1}^{-1}\epsilon _{2} & 0 & 0 & \epsilon _{1}^{-\frac{1}{2}%
}\alpha _{7}%
\end{array}%
\right] .
\end{equation*}%
Then $\phi W=W_{1}$.

\item Assume now that $\epsilon _{1}=0$. This then implies $\epsilon
_{2}\neq 0$. Let $\phi $\ be the following automorphism:%
\begin{equation*}
\phi =\left[ 
\begin{array}{cccc}
1 & 0 & 0 & 0 \\ 
\alpha _{7}^{-2}\left( \alpha _{1}\alpha _{3}-\alpha _{4}\alpha _{7}\right)
& 1 & -\frac{1}{2}\alpha _{1}\alpha _{7}^{-2} & \epsilon _{2}^{-1}\left(
\alpha _{1}\alpha _{8}-\alpha _{6}\alpha _{7}\right) \\ 
-\alpha _{7}^{-1}\alpha _{3} & 0 & \alpha _{7}^{-1} & -\epsilon
_{2}^{-1}\alpha _{7}\alpha _{8} \\ 
0 & 0 & 0 & \epsilon _{2}^{-1}\alpha _{7}^{2}%
\end{array}%
\right] .
\end{equation*}%
Then $\phi W=W_{2}$.
\end{enumerate}
\end{enumerate}

As shown we have two representatives, namely $W_{1}=\left\langle \left[
\delta _{d,d}\right] +\left[ \delta _{b,c}\right] \right\rangle $ and $%
W_{2}=\left\langle \left[ \delta _{a,d}\right] +\left[ \delta _{b,c}\right]
\right\rangle $. Since $\Psi \left( W_{1}\right) =\left( 1\right) $ and $%
\Psi \left( W_{2}\right) =\left( 0\right) $, $W_{1},W_{2}$ are not in the
same orbit. Thus we get the following non-isomorphic algebras:

\begin{itemize}
\item $J_{5,5}:a^{2}=b,d^{2}=e,b\circ c=e.$

\item $J_{5,6}:a^{2}=b,a\circ d=e,b\circ c=e.$
\end{itemize}

\subsection{1-dimensional annihilator extensions of $J_{4,3}$}

The automorphism group $Aut\left( J_{4,3}\right) $ consists of

\begin{equation}
\phi =\left[ 
\begin{array}{cccc}
a_{11} & a_{12} & 0 & 0 \\ 
a_{21} & a_{22} & 0 & 0 \\ 
a_{31} & a_{32} & a_{11}a_{22}+a_{12}a_{21} & a_{34} \\ 
a_{41} & a_{42} & 0 & a_{44}%
\end{array}%
\right] :a_{44}\left( a_{11}^{2}a_{22}^{2}-a_{12}^{2}a_{21}^{2}\right) \neq 0%
\mbox{ and }a_{11}a_{21}=a_{12}a_{22}=0.  \label{AutJ4,3}
\end{equation}%
Choose an arbitrary subspace $W\in U_{1}\left( J_{4,3}\right) $. From Table %
\ref{tabx}, such a subspace is spanned by%
\begin{equation*}
\left[ \theta \right] =\alpha _{1}\left[ \delta _{a,a}\right] +\alpha _{2}%
\left[ \delta _{b,b}\right] +\alpha _{3}\left[ \delta _{d,d}\right] +\alpha
_{4}\left[ \delta _{a,c}\right] +\alpha _{5}\left[ \delta _{a,d}\right]
+\alpha _{6}\left[ \delta _{b,c}\right] +\alpha _{7}\left[ \delta _{b,d}%
\right] +\alpha _{8}\left[ \delta _{c,d}\right]
\end{equation*}%
such that $\left( \alpha _{4},\alpha _{6},\alpha _{8}\right) \neq \left(
0,0,0\right) $ and $\theta ^{\bot }\cap \left\langle c,d\right\rangle =0$.
Further, $\theta ^{\bot }\cap \left\langle c,d\right\rangle =0$ if and only
if the matrix%
\begin{equation*}
A=%
\begin{bmatrix}
\alpha _{4} & \alpha _{5} \\ 
\alpha _{6} & \alpha _{7} \\ 
0 & \alpha _{8} \\ 
\alpha _{8} & \alpha _{3}%
\end{bmatrix}%
\end{equation*}%
has rank $2$. Easy computation shows that $A$ has rank $2$ if and only if $%
\left( \alpha _{8},\alpha _{3}\alpha _{4},\alpha _{3}\alpha _{6},\alpha
_{4}\alpha _{7}-\alpha _{5}\alpha _{6}\right) \neq \left( 0,0,0,0\right) $.
Let $\phi =\big(a_{ij}\big)\in $ $Aut\left( J_{4,3}\right) $. Write $\left[
\phi \theta \right] =\alpha _{1}^{\prime }\left[ \delta _{a,a}\right]
+\alpha _{2}^{\prime }\left[ \delta _{b,b}\right] +\alpha _{3}^{\prime }%
\left[ \delta _{d,d}\right] +\alpha _{4}^{\prime }\left[ \delta _{a,c}\right]
+\alpha _{5}^{\prime }\left[ \delta _{a,d}\right] +\alpha _{6}^{\prime }%
\left[ \delta _{b,c}\right] +\alpha _{7}^{\prime }\left[ \delta _{b,d}\right]
+\alpha _{8}^{\prime }\left[ \delta _{c,d}\right] $. Then%
\begin{eqnarray*}
\alpha _{1}^{\prime } &=&\alpha _{1}a_{11}^{2}+2\alpha
_{5}a_{11}a_{41}+2\alpha _{4}a_{31}a_{11}+\alpha _{2}a_{21}^{2}+2\alpha
_{7}\allowbreak a_{21}a_{41}+2\alpha _{6}a_{31}a_{21}+\alpha
_{3}a_{41}^{2}+2\alpha _{8}a_{31}a_{41}, \\
\alpha _{2}^{\prime } &=&\alpha _{1}a_{12}^{2}+2\alpha
_{5}a_{12}a_{42}+2\alpha _{4}a_{32}a_{12}+\alpha _{2}a_{22}^{2}+2\alpha
_{7}\allowbreak a_{22}a_{42}+2\alpha _{6}a_{32}a_{22}+\alpha
_{3}a_{42}^{2}+2\alpha _{8}a_{32}a_{42}, \\
\alpha _{3}^{\prime } &=&a_{44}\left( \alpha _{3}a_{44}+2\alpha
_{8}a_{34}\right) , \\
\alpha _{4}^{\prime } &=&\left( a_{11}a_{22}+a_{12}a_{21}\right) \left(
\alpha _{4}a_{11}+\alpha _{6}a_{21}+\alpha _{8}a_{41}\right) , \\
\alpha _{5}^{\prime } &=&a_{11}\left( \alpha _{4}a_{34}+\alpha
_{5}a_{44}\right) +a_{21}\left( \alpha _{6}a_{34}+\alpha _{7}a_{44}\right)
+a_{41}\left( \alpha _{3}a_{44}+\alpha _{8}a_{34}\right) +\alpha
_{8}a_{31}a_{44}, \\
\alpha _{6}^{\prime } &=&\left( a_{11}a_{22}+a_{12}a_{21}\right) \left(
\alpha _{4}a_{12}+\alpha _{6}a_{22}+\alpha _{8}a_{42}\right) , \\
\alpha _{7}^{\prime } &=&a_{12}\left( \alpha _{4}a_{34}+\alpha
_{5}a_{44}\right) +a_{22}\left( \alpha _{6}a_{34}+\alpha _{7}a_{44}\right)
+a_{42}\left( \alpha _{3}a_{44}+\alpha _{8}a_{34}\right) +\alpha
_{8}a_{32}a_{44}, \\
\alpha _{8}^{\prime } &=&\alpha _{8}a_{44}\left(
a_{11}a_{22}+a_{12}a_{21}\right) .
\end{eqnarray*}%
As $a_{44}\left( a_{11}a_{22}+a_{12}a_{21}\right) \neq 0$, the equation for $%
\alpha _{8}^{\prime }$ implies that $\alpha _{8}^{\prime }\neq 0$ if and
only if $\alpha _{8}\neq 0$. Thus $\mbox{Orb}\left( \left\langle \left[
\theta \right] :\alpha _{8}\neq 0\right\rangle \right) \cap \mbox{Orb}\left(
\left\langle \left[ \theta \right] :\alpha _{8}=0\right\rangle \right)
=\emptyset $ and hence $Aut(J_{4,3})$ has at least two orbits on $%
U_{1}\left( J_{4,3}\right) $. Therefore we distinguish two cases.

\textsc{Case 1. }$\alpha _{8}\neq 0$\textsc{.}

Set $\alpha =\alpha _{8}^{-2}\left( \alpha _{3}\alpha _{4}^{2}-2\alpha
_{5}\alpha _{4}\alpha _{8}+\alpha _{1}\alpha _{8}^{2}\right) $ and $\beta
=\alpha _{8}^{-2}\left( \alpha _{3}\alpha _{6}^{2}-2\alpha _{6}\alpha
_{7}\alpha _{8}+\alpha _{2}\alpha _{8}^{2}\right) $. We next consider some
case distinctions.

\begin{enumerate}
\item Suppose first that $\alpha =\beta =0$. Let $\phi $ be the following
automorphism: 
\begin{equation*}
\phi =%
\begin{bmatrix}
1 & 0 & 0 & 0 \\ 
0 & 1 & 0 & 0 \\ 
\alpha _{8}^{-2}\left( \alpha _{3}\alpha _{4}-\alpha _{5}\alpha _{8}\right)
& \alpha _{8}^{-2}\left( \alpha _{3}\alpha _{6}-\alpha _{7}\alpha _{8}\right)
& 1 & -\frac{1}{2}\alpha _{3}\alpha _{8}^{-2} \\ 
-\alpha _{4}\alpha _{8}^{-1} & -\alpha _{6}\alpha _{8}^{-1} & 0 & \alpha
_{8}^{-1}%
\end{bmatrix}%
.
\end{equation*}%
Then $\phi W=\left\langle \left[ \delta _{c,d}\right] \right\rangle $. So we
get a representative $W_{1}=\left\langle \left[ \delta _{c,d}\right]
\right\rangle $.

\item Suppose that $\alpha \neq 0$ and $\beta =0$. Let $\phi $ be the
following automorphism: 
\begin{equation*}
\phi =%
\begin{bmatrix}
\alpha ^{-\frac{1}{2}} & 0 & 0 & 0 \\ 
0 & \alpha ^{\frac{1}{2}} & 0 & 0 \\ 
\alpha ^{-\frac{1}{2}}\alpha _{8}^{-2}\left( \alpha _{3}\alpha _{4}-\alpha
_{5}\alpha _{8}\right) & \alpha ^{\frac{1}{2}}\alpha _{8}^{-2}\left( \alpha
_{3}\alpha _{6}-\alpha _{7}\alpha _{8}\right) & 1 & -\frac{1}{2}\alpha
_{3}\alpha _{8}^{-2} \\ 
-\alpha ^{-\frac{1}{2}}\alpha _{4}\alpha _{8}^{-1} & -\alpha ^{\frac{1}{2}%
}\alpha _{6}\alpha _{8}^{-1} & 0 & \alpha _{8}^{-1}%
\end{bmatrix}%
.
\end{equation*}%
Then $\phi W=\left\langle \left[ \delta _{a,a}\right] +\left[ \delta _{c,d}%
\right] \right\rangle $. So we get a representative $W_{2}=\left\langle %
\left[ \delta _{a,a}\right] +\left[ \delta _{c,d}\right] \right\rangle $.

\item Suppose that $\alpha =0$ and $\beta \neq 0$. Let $\phi $ be the
following automorphism:%
\begin{equation*}
\phi =%
\begin{bmatrix}
0 & \beta ^{\frac{1}{2}} & 0 & 0 \\ 
\beta ^{-\frac{1}{2}} & 0 & 0 & 0 \\ 
\beta ^{-\frac{1}{2}}\alpha _{8}^{-2}\left( \alpha _{3}\alpha _{6}-\alpha
_{7}\alpha _{8}\right) & \beta ^{\frac{1}{2}}\alpha _{8}^{-2}\left( \alpha
_{3}\alpha _{4}-\alpha _{5}\alpha _{8}\right) & 1 & -\frac{1}{2}\alpha
_{3}\alpha _{8}^{-2} \\ 
-\beta ^{-\frac{1}{2}}\alpha _{6}\alpha _{8}^{-1} & -\beta ^{\frac{1}{2}%
}\alpha _{4}\alpha _{8}^{-1} & 0 & \alpha _{8}^{-1}%
\end{bmatrix}%
.
\end{equation*}%
Then $\phi W=$ $W_{2}$.

\item Assume now that $\alpha \beta \neq 0$. Let $\phi $ be the following
automorphism:%
\begin{equation*}
\phi =%
\begin{bmatrix}
\alpha ^{-\frac{1}{2}} & 0 & 0 & 0 \\ 
0 & \beta ^{-\frac{1}{2}} & 0 & 0 \\ 
\alpha ^{-\frac{1}{2}}\alpha _{8}^{-2}\left( \alpha _{3}\alpha _{4}-\alpha
_{5}\alpha _{8}\right) & \beta ^{-\frac{1}{2}}\alpha _{8}^{-2}\left( \alpha
_{3}\alpha _{6}-\alpha _{7}\alpha _{8}\right) & \alpha ^{-\frac{1}{2}}\beta
^{-\frac{1}{2}} & -\frac{1}{2}\alpha ^{\frac{1}{2}}\beta ^{\frac{1}{2}%
}\alpha _{3}\alpha _{8}^{-2} \\ 
-\alpha ^{-\frac{1}{2}}\alpha _{4}\alpha _{8}^{-1} & -\beta ^{-\frac{1}{2}%
}\alpha _{6}\alpha _{8}^{-1} & 0 & \alpha ^{\frac{1}{2}}\beta ^{\frac{1}{2}%
}\alpha _{8}^{-1}%
\end{bmatrix}%
.
\end{equation*}%
Then $\phi W=\left\langle \left[ \delta _{a,a}\right] +\left[ \delta _{b,b}%
\right] +\left[ \delta _{c,d}\right] \right\rangle $. Hence we get a
representative $W_{3}=\left\langle \left[ \delta _{a,a}\right] +\left[
\delta _{b,b}\right] +\left[ \delta _{c,d}\right] \right\rangle $.
\end{enumerate}

As previously described we have three representatives, namely $%
W_{1}=\left\langle \left[ \delta _{c,d}\right] \right\rangle
,W_{2}=\left\langle \left[ \delta _{a,a}\right] +\left[ \delta _{c,d}\right]
\right\rangle $ and $W_{3}=\left\langle \left[ \delta _{a,a}\right] +\left[
\delta _{b,b}\right] +\left[ \delta _{c,d}\right] \right\rangle $. Moreover,
we have $\Psi \left( W_{1}\right) =\left( 2\right) $ and $\Psi \left(
W_{2}\right) =\Psi \left( W_{3}\right) =\left( 1\right) $. Therefore, $%
\mbox{Orb}\left( W_{1}\right) \cap \mbox{Orb}\left( W_{i}\right) =\emptyset $
for $i=2,3$. Further we claim that $W_{2},W_{3}$ are not in the same orbit.
Suppose, on the contrary, that $\phi W_{2}=W_{3}$ for some $\phi =\big(a_{ij}%
\big)\in $ $Aut\left( J_{4,3}\right) $. Then there is a $\lambda \in \mathbb{%
F}^{\ast }$ such that%
\begin{eqnarray}
a_{11}^{2}+2a_{31}a_{41} &=&\lambda ,  \label{v1} \\
a_{12}^{2}+2a_{32}a_{42} &=&\lambda ,  \label{v2} \\
2a_{44}a_{34} &=&0,  \notag \\
\left( a_{11}a_{22}+a_{12}a_{21}\right) a_{41} &=&0,  \label{v4} \\
a_{41}a_{34}+a_{31}a_{44} &=&0,  \notag \\
\left( a_{11}a_{22}+a_{12}a_{21}\right) a_{42} &=&0,  \label{v6} \\
a_{42}a_{34}+a_{32}a_{44} &=&0,  \notag \\
a_{44}\left( a_{11}a_{22}+a_{12}a_{21}\right) &=&\lambda .  \notag
\end{eqnarray}%
As $a_{11}a_{22}+a_{12}a_{21}\neq 0$, we obtain from Eqs. $\left( \ref{v4}%
\right) $ and $\left( \ref{v6}\right) $ that $a_{41}=a_{42}=0$. So, from Eq. 
$\left( \ref{v1}\right) $ and Eq. $\left( \ref{v2}\right) $ we obtain $%
a_{11}^{2}a_{12}^{2}=\lambda ^{2}\neq 0$. This then contradicts our
assumption that $\phi =\big(a_{ij}\big)\in $ $Aut\left( J_{4,3}\right) $
because according to\ $\left( \ref{AutJ4,3}\right) $ we either have $%
a_{11}=0 $ or $a_{12}=0$. Therefore any two of $W_{1},W_{2},W_{3}$ are not
in the same orbit. So we get the following pairwise non-isomorphic algebras:

\begin{itemize}
\item $J_{5,7}:a\circ b=c,c\circ d=e.$

\item $J_{5,8}:a\circ b=c,c\circ d=e,a^{2}=e.$

\item $J_{5,9}:a\circ b=c,c\circ d=e,a^{2}=e,b^{2}=e.$
\end{itemize}

\textsc{Case 2. }$\alpha _{8}=0$\textsc{.}

The equation for $\alpha _{3}^{\prime }$ then amounts to $\alpha
_{3}^{\prime }=a_{44}^{2}\alpha _{3}$. Since $a_{44}^{2}\neq 0$, $\alpha
_{3}^{\prime }\neq 0$ if and only if $\alpha _{3}\neq 0$. Thus $\mbox{Orb}%
\left( \left\langle \left[ \theta \right] :\alpha _{3}\neq 0,\alpha
_{8}=0\right\rangle \right) \cap \mbox{Orb}\left( \left\langle \left[ \theta %
\right] :\alpha _{3}=\alpha _{8}=0\right\rangle \right) =\emptyset $. So we
distinguish two cases.

\textsc{Case 2.1. }$\alpha _{3}\neq 0$\textsc{.}

As $\alpha _{8}=0$, we have $\left( \alpha _{4},\alpha _{6}\right) \neq
\left( 0,0\right) $. So $\theta ^{\bot }\cap \left\langle c,d\right\rangle
=0 $. Let us now consider some case distinctions.

\begin{enumerate}
\item Suppose first that $\alpha _{4}\alpha _{6}\neq 0$. Let $\phi $ be the
following automorphism:%
\begin{equation*}
\phi =\left[ 
\begin{array}{cccc}
\alpha _{4}^{-\frac{2}{3}}\alpha _{6}^{\frac{1}{3}} & 0 & 0 & 0 \\ 
0 & \alpha _{4}^{\frac{1}{3}}\alpha _{6}^{-\frac{2}{3}} & 0 & 0 \\ 
\frac{1}{2}\alpha _{3}^{-1}\alpha _{4}^{-\frac{5}{3}}\alpha _{6}^{\frac{1}{3}%
}\left( \alpha _{5}^{2}-\alpha _{1}\alpha _{3}\right) & \frac{1}{2}\alpha
_{3}^{-1}\alpha _{4}^{\frac{1}{3}}\alpha _{6}^{-\frac{5}{3}}\left( \alpha
_{7}^{2}-\alpha _{2}\alpha _{3}\right) & \alpha _{4}^{-\frac{1}{3}}\alpha
_{6}^{-\frac{1}{3}} & 0 \\ 
-\alpha _{3}^{-1}\alpha _{4}^{-\frac{2}{3}}\alpha _{5}\alpha _{6}^{\frac{1}{3%
}} & -\alpha _{3}^{-1}\alpha _{4}^{\frac{1}{3}}\alpha _{6}^{-\frac{2}{3}%
}\alpha _{7} & 0 & \alpha _{3}^{-\frac{1}{2}}%
\end{array}%
\right] .
\end{equation*}%
Then $\phi W=\left\langle \left[ \delta _{d,d}\right] +\left[ \delta _{a,c}%
\right] +\left[ \delta _{b,c}\right] \right\rangle $. Hence we get a
representative $W_{4}=\left\langle \left[ \delta _{d,d}\right] +\left[
\delta _{a,c}\right] +\left[ \delta _{b,c}\right] \right\rangle $.$%
\allowbreak $

\item Suppose that $\alpha _{4}\neq 0$ and $\alpha _{6}=0$. Let $\phi $ be
the following automorphism:%
\begin{equation*}
\phi =\left[ 
\begin{array}{cccc}
\epsilon ^{\frac{1}{4}}\alpha _{4}^{-\frac{1}{2}} & 0 & 0 & 0 \\ 
0 & \epsilon ^{-\frac{1}{2}} & 0 & 0 \\ 
\frac{1}{2}\epsilon ^{\frac{1}{4}}\alpha _{3}^{-1}\left( \alpha
_{5}^{2}-\alpha _{1}\alpha _{3}\right) \alpha _{4}^{-\frac{3}{2}} & 0 & 
\epsilon ^{-\frac{1}{4}}\alpha _{4}^{-\frac{1}{2}} & 0 \\ 
-\epsilon ^{\frac{1}{4}}\alpha _{3}^{-1}\alpha _{4}^{-\frac{1}{2}}\alpha _{5}
& -\epsilon ^{-\frac{1}{2}}\alpha _{3}^{-1}\alpha _{7} & 0 & \alpha _{3}^{-%
\frac{1}{2}}%
\end{array}%
\right]
\end{equation*}%
where $\epsilon =1$ if $\alpha _{7}^{2}-\alpha _{2}\alpha _{3}=0$;
otherwise, $\epsilon =-\alpha _{3}^{-1}\left( \alpha _{7}^{2}-\alpha
_{2}\alpha _{3}\right) $. Then $\phi W=\left\langle \left[ \delta _{d,d}%
\right] +\left[ \delta _{a,c}\right] \right\rangle $ if $\alpha
_{7}^{2}-\alpha _{2}\alpha _{3}=0$, while $\phi W=\left\langle \left[ \delta
_{b,b}\right] +\left[ \delta _{d,d}\right] +\left[ \delta _{a,c}\right]
\right\rangle $ otherwise. Hence we get the representatives $%
W_{5}=\left\langle \left[ \delta _{d,d}\right] +\left[ \delta _{a,c}\right]
\right\rangle ,W_{6}=\left\langle \left[ \delta _{b,b}\right] +\left[ \delta
_{d,d}\right] +\left[ \delta _{a,c}\right] \right\rangle $.

\item Assume now that $\alpha _{4}=0$ and $\alpha _{6}\neq 0$. Let $\phi $
be the following automorphism:%
\begin{equation*}
\phi =\left[ 
\begin{array}{cccc}
0 & \epsilon ^{-\frac{1}{2}} & 0 & 0 \\ 
\epsilon ^{\frac{1}{4}}\alpha _{6}^{-\frac{1}{2}} & 0 & 0 & 0 \\ 
\frac{1}{2}\epsilon ^{\frac{1}{4}}\alpha _{3}^{-1}\left( \alpha
_{7}^{2}-\alpha _{2}\alpha _{3}\right) \alpha _{6}^{-\frac{3}{2}} & 0 & 
\epsilon ^{-\frac{1}{4}}\alpha _{6}^{-\frac{1}{2}} & 0 \\ 
-\epsilon ^{\frac{1}{4}}\alpha _{3}^{-1}\alpha _{6}^{-\frac{1}{2}}\alpha _{7}
& -\epsilon ^{-\frac{1}{2}}\alpha _{3}^{-1}\alpha _{5} & 0 & \alpha _{3}^{-%
\frac{1}{2}}%
\end{array}%
\right]
\end{equation*}%
where $\epsilon =1$ if $\alpha _{5}^{2}-\alpha _{1}\alpha _{3}=0$;
otherwise, $\epsilon =-\alpha _{3}^{-1}\left( \alpha _{5}^{2}-\alpha
_{1}\alpha _{3}\right) $. Then $\phi W=W_{5}$ if $\alpha _{5}^{2}-\alpha
_{1}\alpha _{3}=0$, while $\phi W=W_{6}$ otherwise.
\end{enumerate}

From the above we have three representatives, namely $W_{4}=\left\langle %
\left[ \delta _{d,d}\right] +\left[ \delta _{a,c}\right] +\left[ \delta
_{b,c}\right] \right\rangle ,W_{5}=\left\langle \left[ \delta _{d,d}\right] +%
\left[ \delta _{a,c}\right] \right\rangle $ and $W_{6}=\left\langle \left[
\delta _{b,b}\right] +\left[ \delta _{d,d}\right] +\left[ \delta _{a,c}%
\right] \right\rangle $. Moreover, we have $\Psi \left( W_{4}\right) =\Psi
\left( W_{5}\right) =\left( 1\right) $ and $\Psi \left( W_{6}\right) =\left(
0\right) $. Therefore, $\mbox{Orb}\left( W_{6}\right) \cap \mbox{Orb}\left(
W_{i}\right) =\emptyset $ for $i=4,5$. Further $W_{4},W_{5}$ are not in the
same orbit since the algebras ($J_{5,10},J_{5,11}$) corresponding to these
representatives are non-isomorphic (see Example \ref{non-isom.}). Thus we
get the following pairwise non-isomorphic algebras:

\begin{itemize}
\item $J_{5,10}:a\circ b=c,a\circ c=e,d^{2}=e,b\circ c=e.$

\item $J_{5,11}:a\circ b=c,a\circ c=e,d^{2}=e.$

\item $J_{5,12}:a\circ b=c,a\circ c=e,d^{2}=e,b^{2}=e.$
\end{itemize}

\textsc{Case 2.2. }$\alpha _{3}=0$

In this case $\theta ^{\bot }\cap \left\langle c,d\right\rangle =0$\ if and
only if $\alpha _{4}\alpha _{7}-\alpha _{5}\alpha _{6}\neq 0$. We examine
some cases as follows:

\begin{enumerate}
\item Suppose first that $\alpha _{4}\alpha _{6}\neq 0$. Let $\phi $ be the
following automorphism:%
\begin{equation*}
\phi =\left[ 
\begin{array}{cccc}
\alpha _{4}^{-\frac{2}{3}}\alpha _{6}^{\frac{1}{3}} & 0 & 0 & 0 \\ 
0 & \alpha _{4}^{\frac{1}{3}}\alpha _{6}^{-\frac{2}{3}} & 0 & 0 \\ 
-\frac{1}{2}\alpha _{1}\alpha _{4}^{-\frac{5}{3}}\alpha _{6}^{\frac{1}{3}} & 
-\frac{1}{2}\alpha _{2}\alpha _{4}^{\frac{1}{3}}\alpha _{6}^{-\frac{5}{3}} & 
\alpha _{4}^{-\frac{1}{3}}\alpha _{6}^{-\frac{1}{3}} & \alpha _{4}^{\frac{2}{%
3}}\alpha _{6}^{-\frac{1}{3}}\alpha _{7}\left( \alpha _{4}\alpha _{7}-\alpha
_{5}\alpha _{6}\right) ^{-1} \\ 
0 & 0 & 0 & -\alpha _{4}^{\frac{2}{3}}\alpha _{6}^{\frac{2}{3}}\left( \alpha
_{4}\alpha _{7}-\alpha _{5}\alpha _{6}\right) ^{-1}%
\end{array}%
\right] .
\end{equation*}%
Then $\phi W=\left\langle \left[ \delta _{a,c}\right] +\left[ \delta _{a,d}%
\right] +\left[ \delta _{b,c}\right] \right\rangle $. Hence we get a
representative $W_{7}=\left\langle \left[ \delta _{a,c}\right] +\left[
\delta _{a,d}\right] +\left[ \delta _{b,c}\right] \right\rangle $.

\item Suppose that $\alpha _{4}\neq 0$ and $\alpha _{6}=0$. As $\alpha
_{4}\alpha _{7}-\alpha _{5}\alpha _{6}\neq 0$ we have $\alpha _{7}\neq 0$.
Let $\phi $ be the following automorphism:%
\begin{equation*}
\phi =%
\begin{bmatrix}
\alpha _{4}^{-\frac{1}{2}} & 0 & 0 & 0 \\ 
0 & 1 & 0 & 0 \\ 
-\frac{1}{2}\alpha _{1}\alpha _{4}^{-\frac{3}{2}} & 0 & \alpha _{4}^{-\frac{1%
}{2}} & -\alpha _{4}^{-1}\alpha _{5}\alpha _{7}^{-1} \\ 
0 & -\frac{1}{2}\alpha _{2}\alpha _{7}^{-1} & 0 & \alpha _{7}^{-1}%
\end{bmatrix}%
.
\end{equation*}%
Then $\phi W=\left\langle \left[ \delta _{a,c}\right] +\left[ \delta _{b,d}%
\right] \right\rangle $. Hence we get a representative $W_{8}=\left\langle %
\left[ \delta _{a,c}\right] +\left[ \delta _{b,d}\right] \right\rangle $.

\item Assume now that $\alpha _{4}=0$ and $\alpha _{6}\neq 0$. As $\alpha
_{4}\alpha _{7}-\alpha _{5}\alpha _{6}\neq 0$, we thus have $\alpha _{5}\neq
0$. Let $\phi $ be the following automorphism:%
\begin{equation*}
\phi =%
\begin{bmatrix}
0 & 1 & 0 & 0 \\ 
\alpha _{6}^{-\frac{1}{2}} & 0 & 0 & 0 \\ 
-\frac{1}{2}\alpha _{2}\alpha _{6}^{-\frac{3}{2}} & 0 & \alpha _{6}^{-\frac{1%
}{2}} & -\alpha _{6}^{-1}\alpha _{7}\alpha _{5}^{-1} \\ 
0 & -\frac{1}{2}\alpha _{1}\alpha _{5}^{-1} & 0 & \alpha _{5}^{-1}%
\end{bmatrix}%
.
\end{equation*}%
Then $\phi W=W_{8}$. $\allowbreak $
\end{enumerate}

As seen above we have two representatives, namely $W_{7}=\left\langle \left[
\delta _{a,c}\right] +\left[ \delta _{a,d}\right] +\left[ \delta _{b,c}%
\right] \right\rangle $ and $W_{8}=\left\langle \left[ \delta _{a,c}\right] +%
\left[ \delta _{b,d}\right] \right\rangle $. The algebras ($%
J_{5,13},J_{5,14} $) corresponding to these representatives are
non-isomorphic (see Example \ref{non-isom.}). So we get the following
non-isomorphic algebras:

\begin{itemize}
\item $J_{5,13}:a\circ b=c,a\circ c=e,a\circ d=e,b\circ c=e.$

\item $J_{5,14}:a\circ b=c,a\circ c=e,b\circ d=e.$
\end{itemize}

\subsection{1-dimensional annihilator extensions of $J_{4,4}$}

The automorphism group$\ Aut(J_{4,4})$ consists of:

\begin{equation*}
\phi =\left[ 
\begin{array}{cccc}
a_{11} & 0 & 0 & 0 \\ 
a_{21} & a_{11}^{2} & 0 & 0 \\ 
a_{31} & 2a_{11}a_{21} & a_{11}^{3} & a_{34} \\ 
a_{41} & 0 & 0 & a_{44}%
\end{array}%
\right] :a_{11}a_{44}\neq 0.
\end{equation*}%
Choose an arbitrary subspace $W\in U_{1}\left( J_{4,4}\right) $. From Table %
\ref{tabx}, such a subspace is spanned by $\left[ \theta \right] =\alpha
_{1}\left( \left[ \delta _{a,c}\right] +\left[ \delta _{b,b}\right] \right)
+\alpha _{2}\left[ \delta _{a,d}\right] +\alpha _{3}\left[ \delta _{b,d}%
\right] +\alpha _{4}\left[ \delta _{d,d}\right] $ such that $\alpha
_{1}\alpha _{3}\neq 0$. Let $\phi =\big(a_{ij}\big)\in $ $Aut\left(
J_{4,4}\right) $. Then%
\begin{equation*}
\phi W=\left\langle \alpha _{1}a_{11}^{4}\left( \left[ \delta _{a,c}\right] +%
\left[ \delta _{b,b}\right] \right) +\left( a_{44}\left( \alpha
_{2}a_{11}+\alpha _{3}a_{21}+\alpha _{4}a_{41}\right) +\alpha
_{1}a_{11}a_{34}\right) \left[ \delta _{a,d}\right] +\alpha
_{3}a_{11}^{2}a_{44}\left[ \delta _{b,d}\right] +\alpha _{4}a_{44}^{2}\left[
\delta _{d,d}\right] \right\rangle .
\end{equation*}%
Now let $\phi $ be the following automorphism:%
\begin{equation*}
\phi =\left[ 
\begin{array}{cccc}
\alpha _{1}^{-\frac{1}{4}} & 0 & 0 & 0 \\ 
0 & \alpha _{1}^{-\frac{1}{2}} & 0 & 0 \\ 
0 & 0 & \alpha _{1}^{-\frac{3}{4}} & -\alpha _{1}^{-\frac{1}{2}}\alpha
_{2}\alpha _{3}^{-1} \\ 
0 & 0 & 0 & \alpha _{1}^{\frac{1}{2}}\alpha _{3}^{-1}%
\end{array}%
\right] .
\end{equation*}%
Then $\phi W=\left\langle \left[ \delta _{a,c}\right] +\left[ \delta _{b,b}%
\right] +\left[ \delta _{b,d}\right] +\alpha _{1}\alpha _{3}^{-2}\alpha _{4}%
\left[ \delta _{d,d}\right] \right\rangle $. Set $\alpha =\alpha _{1}\alpha
_{3}^{-2}\alpha _{4}$. Then we get the representatives $W^{\alpha \in 
\mathbb{F}}=\left\langle \left[ \delta _{a,c}\right] +\left[ \delta _{b,b}%
\right] +\left[ \delta _{b,d}\right] +\alpha \left[ \delta _{d,d}\right]
\right\rangle $. We claim that $W^{\alpha },W^{\beta }$ are in the same
orbit if and only if $\alpha =\beta $. To prove this, assume $\phi W^{\alpha
}=W^{\beta }$ for some $\phi =\big(a_{ij}\big)\in Aut(J_{4,4})$. Then there
is a $\lambda \in \mathbb{F}^{\ast }$ such that%
\begin{equation*}
\left[ a_{11}^{4}\left( \delta _{a,c}+\delta _{b,b}\right) +\left(
a_{44}a_{21}+\alpha a_{44}a_{41}+a_{11}a_{34}\right) \delta
_{a,d}+a_{11}^{2}a_{44}\delta _{b,d}+\alpha a_{44}^{2}\delta _{d,d}\right]
=\lambda \left[ \delta _{a,c}+\delta _{b,b}+\delta _{b,d}+\beta \delta _{d,d}%
\right] \text{.}
\end{equation*}%
This then amounts to the following polynomial equations:%
\begin{eqnarray}
a_{11}^{4} &=&\lambda ,  \label{Eq.1} \\
a_{44}a_{21}+\alpha a_{44}a_{41}+a_{11}a_{34} &=&0,  \notag \\
a_{11}^{2}a_{44} &=&\lambda ,  \label{Eq.3} \\
a_{44}^{2}\alpha  &=&\lambda \beta .  \label{Eq.4}
\end{eqnarray}%
Using Eq. $\left( \ref{Eq.1}\right) $\ and Eq. $\left( \ref{Eq.3}\right) $,
we obtain $a_{44}^{2}=\lambda ^{2}a_{11}^{-4}=\lambda $. Thus, from Eq. $%
\left( \ref{Eq.4}\right) $ we obtain $\alpha =\beta $, as claimed. Therefore
we get the algebras:

\begin{itemize}
\item $J_{5,15}^{\alpha \in \mathbb{F}}:a^{2}=b,a\circ b=c,a\circ
c=e,b^{2}=e,b\circ d=e,d^{2}=\alpha e$. Isomorphism: $J_{5,15}^{\alpha
}\cong J_{5,15}^{\beta }$ if and only if $\alpha =\beta $.
\end{itemize}

\subsection{1-dimensional annihilator extensions of $J_{4,5}$}

\begin{lem}
\label{lemma}Let $\alpha =\left[ \alpha _{1},\alpha _{2},\alpha _{3}\right]
\neq 0$. Then there exists an invertible matrix 
\begin{equation}
A=\left[ 
\begin{array}{ccc}
a_{11} & a_{12} & a_{13} \\ 
a_{21} & a_{22} & a_{23} \\ 
a_{31} & a_{32} & a_{33}%
\end{array}%
\right] \mbox{ with }\left[ 
\begin{array}{ccc}
a_{23} & a_{13} & a_{33} \\ 
a_{22} & a_{12} & a_{32} \\ 
a_{21} & a_{11} & a_{31}%
\end{array}%
\right] \left[ 
\begin{array}{ccc}
a_{11} & a_{12} & a_{13} \\ 
a_{21} & a_{22} & a_{23} \\ 
a_{31} & a_{32} & a_{33}%
\end{array}%
\right] =\left[ 
\begin{array}{ccc}
0 & 0 & c \\ 
c & 0 & 0 \\ 
0 & c & 0%
\end{array}%
\right]  \label{form}
\end{equation}%
such that $\alpha A\in \left\{ \left[ 1,0,0\right] ,\left[ 0,0,1\right]
\right\} $. Further, if $\alpha =\left[ 1,0,0\right] $, then there is no
invertible matrix of the form $\left( \ref{form}\right) $\ such that $\left[
1,0,0\right] A=\left[ 0,0,1\right] $.$\allowbreak $
\end{lem}

\begin{proof}
Let $\alpha =\left[ \alpha _{1},\alpha _{2},\alpha _{3}\right] \neq 0$.
Suppose first that $\alpha _{3}=0$. Next consider the following cases:

\begin{enumerate}
\item If $\alpha _{1}\neq 0$ and $\alpha _{2}=0$, we take%
\begin{equation*}
A=\left[ 
\begin{array}{ccc}
\alpha _{1}^{-1} & 0 & 0 \\ 
0 & \alpha _{1} & 0 \\ 
0 & 0 & 1%
\end{array}%
\right] .
\end{equation*}%
Then $\alpha A=\left[ 1,0,0\right] $.

\item If $\alpha _{1}=0$ and $\alpha _{2}\neq 0$, we take%
\begin{equation*}
A=\left[ 
\begin{array}{ccc}
0 & \alpha _{2} & 0 \\ 
\alpha _{2}^{-1} & 0 & 0 \\ 
0 & 0 & 1%
\end{array}%
\right] .
\end{equation*}%
Then $\alpha A=\left[ 1,0,0\right] $.

\item If $\alpha _{1}\alpha _{2}\neq 0$, we take%
\begin{equation*}
A=\left[ 
\begin{array}{ccc}
-\sqrt{\frac{1}{8}\alpha _{1}^{-3}\alpha _{2}} & \sqrt{\frac{1}{8}\alpha
_{1}^{-1}\alpha _{2}^{-1}} & \frac{1}{2}\alpha _{1}^{-1} \\ 
\sqrt{\frac{1}{8}\alpha _{1}^{-1}\alpha _{2}^{-1}} & -\sqrt{\frac{1}{8}%
\alpha _{1}\alpha _{2}^{-3}} & \frac{1}{2}\alpha _{2}^{-1} \\ 
\frac{1}{2}\alpha _{1}^{-1} & \frac{1}{2}\alpha _{2}^{-1} & 0%
\end{array}%
\right] .
\end{equation*}%
Then $\alpha A=\left[ 0,0,1\right] $.
\end{enumerate}

Assume now that $\alpha _{3}\neq 0$. Next consider the following cases:

\begin{enumerate}
\item If $\alpha _{1}\neq 0$ and $\alpha _{2}=0$, we take%
\begin{equation*}
A=\left[ 
\begin{array}{ccc}
0 & -\alpha _{1}^{-1}\alpha _{3} & 0 \\ 
-\alpha _{1}\alpha _{3}^{-3} & \frac{1}{2}\alpha _{1}\alpha _{3}^{-1} & 
\alpha _{1}\alpha _{3}^{-2} \\ 
0 & 1 & \alpha _{3}^{-1}%
\end{array}%
\right] .
\end{equation*}%
Then $\alpha A=\left[ 0,0,1\right] $.

\item If $\alpha _{1}=0$ and $\alpha _{2}\neq 0$, we take%
\begin{equation*}
A=\left[ 
\begin{array}{ccc}
-\alpha _{2}\alpha _{3}^{-3} & \frac{1}{2}\alpha _{2}\alpha _{3}^{-1} & 
\alpha _{2}\alpha _{3}^{-2} \\ 
0 & -\alpha _{2}^{-1}\alpha _{3} & 0 \\ 
0 & 1 & \alpha _{3}^{-1}%
\end{array}%
\right] .
\end{equation*}%
Then $\alpha A=\left[ 0,0,1\right] $.

\item If $\alpha _{1}=\alpha _{2}=0$, we take%
\begin{equation*}
A=\alpha _{3}^{-1}\left[ 
\begin{array}{ccc}
1 & 0 & 0 \\ 
0 & 1 & 0 \\ 
0 & 0 & 1%
\end{array}%
\right] .
\end{equation*}%
Then $\alpha A=\left[ 0,0,1\right] $.

\item If $\alpha _{1}\alpha _{2}\neq 0$ and $\epsilon =2\alpha _{1}\alpha
_{2}+\alpha _{3}^{2}\neq 0$, we take%
\begin{equation*}
A=\epsilon ^{-1}\left[ 
\begin{array}{ccc}
\frac{1}{2}\alpha _{1}^{-1}\sqrt{\alpha _{1}\alpha _{2}}\left( -\sqrt{%
\epsilon }-\alpha _{3}\right) & \frac{1}{2}\alpha _{1}^{-1}\sqrt{\alpha
_{1}\alpha _{2}}\left( \sqrt{\epsilon }-\alpha _{3}\right) & \alpha _{2} \\ 
\frac{1}{2}\alpha _{2}^{-1}\sqrt{\alpha _{1}\alpha _{2}}\left( \sqrt{%
\epsilon }-\alpha _{3}\right) & \frac{1}{2}\alpha _{2}^{-1}\sqrt{\alpha
_{1}\alpha _{2}}\left( -\sqrt{\epsilon }-\alpha _{3}\right) & \alpha _{1} \\ 
\sqrt{\alpha _{1}\alpha _{2}} & \sqrt{\alpha _{1}\alpha _{2}} & \alpha _{3}%
\end{array}%
\right] .
\end{equation*}%
Then $\alpha A=\left[ 0,0,1\right] $.

\item If $\alpha _{1}\alpha _{2}\neq 0$ and $2\alpha _{1}\alpha _{2}+\alpha
_{3}^{2}=0$, we take%
\begin{equation*}
A=\left[ 
\begin{array}{ccc}
\frac{1}{4}\alpha _{1}^{-1} & -\frac{1}{2}\alpha _{1}^{-1}\alpha _{3}^{2} & -%
\frac{1}{2}\alpha _{1}^{-1}\alpha _{3} \\ 
\frac{1}{4}\alpha _{2}^{-1} & -\frac{1}{2}\alpha _{2}^{-1}\alpha _{3}^{2} & 
\frac{1}{2}\alpha _{2}^{-1}\alpha _{3} \\ 
\frac{1}{2}\alpha _{3}^{-1} & \alpha _{3} & 0%
\end{array}%
\right] .
\end{equation*}%
Then $\alpha A=\left[ 1,0,0\right] $.
\end{enumerate}

Hence there exists an invertible matrix of the form $\left( \ref{form}%
\right) $ such that $\alpha A\in \left\{ \left[ 1,0,0\right] ,\left[ 0,0,1%
\right] \right\} $. Further, let $A$ be of the form $\left( \ref{form}%
\right) $ such that $\left[ 1,0,0\right] A=\left[ 0,0,1\right] $.$%
\allowbreak $ Then $a_{11}=a_{12}=0$ and $a_{13}=1$. Moreover,%
\begin{equation*}
\left[ 
\begin{array}{ccc}
a_{21}+a_{31}a_{33} & a_{22}+a_{32}a_{33} & a_{33}^{2}+2a_{23} \\ 
a_{31}a_{32} & a_{32}^{2} & a_{22}+a_{32}a_{33} \\ 
a_{31}^{2} & a_{31}a_{32} & a_{21}+a_{31}a_{33}%
\end{array}%
\right] =\left[ 
\begin{array}{ccc}
0 & 0 & c \\ 
c & 0 & 0 \\ 
0 & c & 0%
\end{array}%
\right] .
\end{equation*}%
This then implies $a_{31}=a_{32}=a_{21}=a_{22}=c=0$, which is impossible as $%
\det A\neq 0$. This completes the proof.
\end{proof}

Let us now find a set of representatives of the orbits of $U_{1}\left(
J_{4,5}\right) $ under the action of $Aut(J_{4,5})$. Choose an arbitrary
subspace $W\in U_{1}\left( J_{4,5}\right) $. From Table \ref{tabx}, such a
subspace is spanned by $\left[ \theta \right] =\alpha _{1}\left[ \delta
_{a,d}\right] +\alpha _{2}\left[ \delta _{b,d}\right] +\alpha _{3}\left[
\delta _{c,d}\right] +\alpha _{4}\left[ \delta _{a,a}\right] +\alpha _{5}%
\left[ \delta _{b,b}\right] +\alpha _{6}\left[ \delta _{a,b}\right] +\alpha
_{7}\left[ \delta _{a,c}\right] +\alpha _{8}\left[ \delta _{b,c}\right] $
such that $\left( \alpha _{1},\alpha _{2},\alpha _{3}\right) \neq \left(
0,0,0\right) $. Furthermore the automorphism group $Aut(J_{4,5})$ consists of%
\begin{equation}
\phi =\left[ 
\begin{array}{cccc}
a_{11} & a_{12} & a_{13} & 0 \\ 
a_{21} & a_{22} & a_{23} & 0 \\ 
a_{31} & a_{32} & a_{33} & 0 \\ 
a_{41} & a_{42} & a_{43} & a_{44}%
\end{array}%
\right] :\det \phi \neq 0\mbox{ and }\left[ 
\begin{array}{ccc}
a_{23} & a_{13} & a_{33} \\ 
a_{22} & a_{12} & a_{32} \\ 
a_{21} & a_{11} & a_{31}%
\end{array}%
\right] \left[ 
\begin{array}{ccc}
a_{11} & a_{12} & a_{13} \\ 
a_{21} & a_{22} & a_{23} \\ 
a_{31} & a_{32} & a_{33}%
\end{array}%
\right] =\left[ 
\begin{array}{ccc}
0 & 0 & a_{44} \\ 
a_{44} & 0 & 0 \\ 
0 & a_{44} & 0%
\end{array}%
\right] .  \label{AutJ4,5}
\end{equation}%
Let $\phi =\big(a_{ij}\big)$ and write $\left[ \phi \theta \right] =\alpha
_{1}^{\prime }\left[ \delta _{a,d}\right] +\alpha _{2}^{\prime }\left[
\delta _{b,d}\right] +\alpha _{3}^{\prime }\left[ \delta _{c,d}\right]
+\alpha _{4}^{\prime }\left[ \delta _{a,a}\right] +\alpha _{5}^{\prime }%
\left[ \delta _{b,b}\right] +\alpha _{6}^{\prime }\left[ \delta _{a,b}\right]
+\alpha _{7}^{\prime }\left[ \delta _{a,c}\right] +\alpha _{8}^{\prime }%
\left[ \delta _{b,c}\right] $. Then%
\begin{eqnarray*}
\alpha _{1}^{\prime } &=&a_{44}\left( \alpha _{1}a_{11}+\alpha
_{2}a_{21}+\alpha _{3}a_{31}\right) , \\
\alpha _{2}^{\prime } &=&a_{44}\left( \alpha _{1}a_{12}+\alpha
_{2}a_{22}+\alpha _{3}a_{32}\right) , \\
\alpha _{3}^{\prime } &=&a_{44}\left( \alpha _{1}a_{13}+\alpha
_{2}a_{23}+\alpha _{3}a_{33}\right) , \\
\alpha _{4}^{\prime } &=&\alpha _{4}a_{11}^{2}+\alpha _{5}a_{21}^{2}+2\alpha
_{1}a_{11}a_{41}+2\alpha _{2}a_{21}a_{41}+2\alpha _{6}a_{11}a_{21}+2\alpha
_{3}a_{31}a_{41}+2\alpha _{7}a_{11}a_{31}+2\alpha _{8}a_{21}a_{31}, \\
\alpha _{5}^{\prime } &=&\alpha _{4}a_{12}^{2}+\alpha _{5}a_{22}^{2}+2\alpha
_{1}a_{12}a_{42}+2\alpha _{2}a_{22}a_{42}+2\alpha _{6}a_{12}a_{22}+2\alpha
_{3}a_{32}a_{42}+2\alpha _{7}a_{12}a_{32}+2\alpha _{8}a_{22}a_{32}, \\
\alpha _{6}^{\prime } &=&a_{11}\left( \alpha _{1}a_{42}+\alpha
_{4}a_{12}+\alpha _{6}a_{22}+\alpha _{7}a_{32}\right) +a_{21}\left( \alpha
_{2}a_{42}+\alpha _{5}a_{22}+\alpha _{6}a_{12}+\alpha _{8}a_{32}\right)
+a_{41}\left( \alpha _{1}a_{12}+\alpha _{2}a_{22}+\alpha _{3}a_{32}\right) \\
&&+a_{31}\left( \alpha _{3}a_{42}+\alpha _{7}a_{12}+\alpha _{8}a_{22}\right)
-\alpha _{4}a_{13}^{2}-\alpha _{5}a_{23}^{2}-2\alpha
_{1}a_{13}a_{43}-2\alpha _{2}a_{23}a_{43}-2\alpha _{6}a_{13}a_{23}-2\alpha
_{3}a_{33}a_{43} \\
&&-2\alpha _{7}a_{13}a_{33}-2\alpha _{8}a_{23}a_{33}, \\
\alpha _{7}^{\prime } &=&a_{11}\left( \alpha _{1}a_{43}+\alpha
_{4}a_{13}+\alpha _{6}a_{23}+\alpha _{7}a_{33}\right) +a_{21}\left( \alpha
_{2}a_{43}+\alpha _{5}a_{23}+\alpha _{6}a_{13}+\alpha _{8}a_{33}\right)
+a_{41}\left( \alpha _{1}a_{13}+\alpha _{2}a_{23}+\alpha _{3}a_{33}\right) \\
&&+a_{31}\left( \alpha _{3}a_{43}+\alpha _{7}a_{13}+\alpha _{8}a_{23}\right)
, \\
\alpha _{8}^{\prime } &=&a_{12}\left( \alpha _{1}a_{43}+\alpha
_{4}a_{13}+\alpha _{6}a_{23}+\alpha _{7}a_{33}\right) +a_{22}\left( \alpha
_{2}a_{43}+\alpha _{5}a_{23}+\alpha _{6}a_{13}+\alpha _{8}a_{33}\right)
+a_{42}\left( \alpha _{1}a_{13}+\alpha _{2}a_{23}+\alpha _{3}a_{33}\right) \\
&&+a_{32}\left( \alpha _{3}a_{43}+\alpha _{7}a_{13}+\alpha _{8}a_{23}\right)
.
\end{eqnarray*}%
According to Lemma \ref{lemma}, we may now assume without loss of generality
that $\left( \alpha _{1},\alpha _{2},\alpha _{3}\right) \in \left\{ \left(
1,0,0\right) ,\left( 0,0,1\right) \right\} $. Moreover, $\mbox{Orb}\left(
\left\langle \left[ \theta \right] :\alpha _{1}=1,\alpha _{2}=\alpha
_{3}=0\right\rangle \right) \cap \mbox{Orb}\left( \left\langle \left[ \theta %
\right] :\alpha _{1}=\alpha _{2}=0,\alpha _{3}=1\right\rangle \right)
=\emptyset \ $and hence $Aut(J_{4,5})$ has at least two orbits on $%
U_{1}\left( J_{4,5}\right) $. So we distinguish two cases.

\textsc{Case 1. }Suppose first that\textsc{\ }$\left( \alpha _{1},\alpha
_{2},\alpha _{3}\right) =\left( 1,0,0\right) $.

Let us consider a few cases:

\begin{enumerate}
\item Suppose that $\alpha _{5}=\alpha _{8}=0$. Let $\phi $ be the following
automorphism: 
\begin{equation*}
\phi =\left[ 
\begin{array}{cccc}
1 & 0 & 0 & 0 \\ 
0 & 1 & 0 & 0 \\ 
0 & 0 & 1 & 0 \\ 
-\frac{1}{2}\alpha _{4} & -\alpha _{6} & -\alpha _{7} & 1%
\end{array}%
\right] .
\end{equation*}%
Then $\left[ \phi \theta \right] =\left[ \delta _{a,d}\right] $. Hence we
get a representative $W_{1}=\left\langle \left[ \delta _{a,d}\right]
\right\rangle $.

\item Suppose that $\alpha _{5}=0$ and $\alpha _{8}\neq 0$. Let $\phi $ be
the following automorphism:%
\begin{equation*}
\phi =\left[ 
\begin{array}{cccc}
\alpha _{8}^{\frac{2}{5}} & 0 & 0 & 0 \\ 
0 & \alpha _{8}^{-\frac{4}{5}} & 0 & 0 \\ 
0 & 0 & \alpha _{8}^{-\frac{1}{5}} & 0 \\ 
-\frac{1}{2}\alpha _{4}\alpha _{8}^{\frac{2}{5}} & -\alpha _{6}\alpha _{8}^{-%
\frac{4}{5}} & -\alpha _{7}\alpha _{8}^{-\frac{1}{5}} & \alpha _{8}^{-\frac{2%
}{5}}%
\end{array}%
\right] .
\end{equation*}%
Then $\left[ \phi \theta \right] =\left[ \delta _{a,d}\right] +\left[ \delta
_{b,c}\right] $. Hence we get a representative $W_{2}=\left\langle \left[
\delta _{a,d}\right] +\left[ \delta _{b,c}\right] \right\rangle $.

\item Assume now that\textsc{\ }$\alpha _{5}\neq 0$. Let $\phi $ be the
following automorphism:%
\begin{equation*}
\phi =\left[ 
\begin{array}{cccc}
\alpha _{5}^{\frac{1}{4}} & 0 & 0 & 0 \\ 
0 & \alpha _{5}^{-\frac{1}{2}} & -\alpha _{5}^{-\frac{9}{8}}\alpha _{8} & 0
\\ 
0 & 0 & \alpha _{5}^{-\frac{1}{8}} & 0 \\ 
-\frac{1}{2}\left( \alpha _{5}^{-\frac{6}{4}}\alpha _{8}^{2}+\alpha
_{4}\alpha _{5}^{\frac{1}{4}}\right) & -\alpha _{5}^{-\frac{1}{2}}\alpha _{6}
& \left( \alpha _{6}\alpha _{8}-\alpha _{5}\alpha _{7}\right) \alpha _{5}^{-%
\frac{9}{8}} & \alpha _{5}^{-\frac{1}{4}}%
\end{array}%
\right] .
\end{equation*}%
Then $\left[ \phi \theta \right] =\left[ \delta _{a,d}\right] +\left[ \delta
_{b,b}\right] $. So we get a representative $W_{3}=\left\langle \left[
\delta _{a,d}\right] +\left[ \delta _{b,b}\right] \right\rangle $.$%
\allowbreak $
\end{enumerate}

As shown we have three representatives, namely $W_{1}=\left\langle \left[
\delta _{a,d}\right] \right\rangle ,W_{2}=\left\langle \left[ \delta _{a,d}%
\right] +\left[ \delta _{b,c}\right] \right\rangle ,W_{3}=\left\langle \left[
\delta _{a,d}\right] +\left[ \delta _{b,b}\right] \right\rangle $. Any two
of them are not in the same orbit since $\Psi \left( W_{1}\right) =\left(
2\right) ,\Psi \left( W_{2}\right) =\left( 0\right) $ and $\Psi \left(
W_{3}\right) =\left( 1\right) $. So we get the following pairwise
non-isomorphic algebras:

\begin{itemize}
\item $J_{5,16}:a\circ b=d,c^{2}=d,a\circ d=e.$

\item $J_{5,17}:a\circ b=d,c^{2}=d,a\circ d=e,b\circ c=e.$

\item $J_{5,18}:a\circ b=d,c^{2}=d,a\circ d=e,b^{2}=e.$
\end{itemize}

\textsc{Case 2. }Assume now that $\left( \alpha _{1},\alpha _{2},\alpha
_{3}\right) =\left( 0,0,1\right) $.

Let us consider a few cases:

\begin{enumerate}
\item Suppose that $\alpha _{4}=\alpha _{5}=0$. Let $\phi $ be the following
automorphism:%
\begin{equation*}
\phi =%
\begin{bmatrix}
1 & 0 & 0 & 0 \\ 
0 & 1 & 0 & 0 \\ 
0 & 0 & 1 & 0 \\ 
-\alpha _{7} & -\alpha _{8} & \frac{1}{2}\alpha _{6} & 1%
\end{bmatrix}%
.
\end{equation*}%
Then $\left[ \phi \theta \right] =\left[ \delta _{c,d}\right] $. Hence we
get a representative $W_{4}=\left\langle \left[ \delta _{c,d}\right]
\right\rangle $.

\item Suppose that $\alpha _{4}\neq 0$ and $\alpha _{5}=0$. Let $\phi $ be
the following automorphism:%
\begin{equation*}
\phi =%
\begin{bmatrix}
\alpha _{4}^{-\frac{1}{2}} & 0 & 0 & 0 \\ 
0 & \alpha _{4}^{\frac{1}{2}} & 0 & 0 \\ 
0 & 0 & 1 & 0 \\ 
-\alpha _{4}^{-\frac{1}{2}}\alpha _{7} & -\alpha _{4}^{\frac{1}{2}}\alpha
_{8} & \frac{1}{2}\alpha _{6} & 1%
\end{bmatrix}%
.
\end{equation*}%
Then $\left[ \phi \theta \right] =\left[ \delta _{c,d}\right] +\left[ \delta
_{a,a}\right] $. Hence we get a representative $W_{5}=\left\langle \left[
\delta _{c,d}\right] +\left[ \delta _{a,a}\right] \right\rangle $.

\item Suppose that $\alpha _{4}=0$ and $\alpha _{5}\neq 0$. Let $\phi $ be
the following automorphism:%
\begin{equation*}
\phi =%
\begin{bmatrix}
0 & \alpha _{5}^{\frac{1}{2}} & 0 & 0 \\ 
\alpha _{5}^{-\frac{1}{2}} & 0 & 0 & 0 \\ 
0 & 0 & 1 & 0 \\ 
-\alpha _{5}^{-\frac{1}{2}}\alpha _{8} & -\alpha _{5}^{\frac{1}{2}}\alpha
_{7} & \frac{1}{2}\alpha _{6} & 1%
\end{bmatrix}%
.
\end{equation*}%
Then $\left[ \phi \theta \right] =\left[ \delta _{c,d}\right] +\left[ \delta
_{a,a}\right] $. So we get again a representative $W_{5}=\left\langle \left[
\delta _{c,d}\right] +\left[ \delta _{a,a}\right] \right\rangle $.

\item Assume now that $\alpha _{4}\alpha _{5}\neq 0$. Let $\phi $ be the
following automorphism:%
\begin{equation*}
\phi =%
\begin{bmatrix}
\alpha _{4}^{\frac{1}{4}}\alpha _{5}^{\frac{3}{4}} & 0 & 0 & 0 \\ 
0 & \alpha _{4}^{\frac{3}{4}}\alpha _{5}^{\frac{1}{4}} & 0 & 0 \\ 
0 & 0 & \alpha _{4}^{\frac{1}{2}}\alpha _{5}^{\frac{1}{2}} & 0 \\ 
-\alpha _{4}^{\frac{1}{4}}\alpha _{5}^{\frac{3}{4}}\alpha _{7} & -\alpha
_{4}^{\frac{3}{4}}\alpha _{5}^{\frac{1}{4}}\alpha _{8} & \frac{1}{2}\alpha
_{4}^{\frac{1}{2}}\alpha _{5}^{\frac{1}{2}}\alpha _{6} & \alpha _{4}\alpha
_{5}%
\end{bmatrix}%
.
\end{equation*}%
Then $\left[ \phi \theta \right] =\alpha _{4}^{\frac{3}{2}}\alpha _{5}^{%
\frac{3}{2}}\left( \left[ \delta _{c,d}\right] +\left[ \delta _{a,a}\right] +%
\left[ \delta _{b,b}\right] \right) $. Consequently, we have a
representative $W_{6}=\left\langle \left[ \delta _{c,d}\right] +\left[
\delta _{a,a}\right] +\left[ \delta _{b,b}\right] \right\rangle $.
\end{enumerate}

As shown we have three representatives, namely $W_{4}=\left\langle \delta
_{c,d}\right\rangle ,W_{5}=\left\langle \left[ \delta _{c,d}\right] +\left[
\delta _{a,a}\right] \right\rangle $ and $W_{6}=\left\langle \left[ \delta
_{c,d}\right] +\left[ \delta _{a,a}\right] +\left[ \delta _{b,b}\right]
\right\rangle $. Moreover, we have $\Psi \left( W_{4}\right) =\left(
2\right) $ and $\Psi \left( W_{5}\right) =\Psi \left( W_{6}\right) =\left(
1\right) $. Therefore, for $i=5,6$, $\mbox{Orb}\left( W_{4}\right) \cap %
\mbox{Orb}\left( W_{i}\right) =\emptyset $. Further we claim that $%
W_{5},W_{6}$ are not in the same orbit. To see this, assume $\phi
W_{5}=W_{6} $ for some $\phi =\big(a_{ij}\big)\in Aut(J_{4,5})$. Then there
is a $\lambda \in \mathbb{F}^{\ast }$ such that $\phi \left( \left[ \delta
_{c,d}\right] +\left[ \delta _{a,a}\right] \right) =\lambda \left( \left[
\delta _{c,d}\right] +\left[ \delta _{a,a}\right] +\left[ \delta _{b,b}%
\right] \right) $. By comparing the coefficients of $\left[ \delta _{a,d}%
\right] ,\left[ \delta _{b,d}\right] ,\left[ \delta _{a,a}\right] $ and $%
\left[ \delta _{b,b}\right] $ on the two sides, we obtain%
\begin{eqnarray}
a_{44}a_{31} &=&0,  \label{F8} \\
a_{44}a_{32} &=&0,  \label{F9} \\
a_{12}^{2}+2a_{32}a_{42} &=&\lambda ,  \label{F10} \\
a_{11}^{2}+2a_{31}a_{41} &=&\lambda .  \label{F11}
\end{eqnarray}%
Moreover, from $\left( \ref{AutJ4,5}\right) $ we have 
\begin{eqnarray}
a_{44} &\neq &0,  \label{F4} \\
a_{31}^{2}+2a_{11}a_{21} &=&0,  \label{F5} \\
a_{32}^{2}+2a_{12}a_{22} &=&0,  \label{F6} \\
a_{11}a_{22}+a_{12}a_{21}+a_{31}a_{32} &=&a_{44}.  \label{F7}
\end{eqnarray}%
From $\left( \ref{F8}\right) $ and $\left( \ref{F9}\right) $, we obtain $%
a_{31}=a_{32}=0$. So, from Eq. $\left( \ref{F10}\right) $ and Eq. $\left( %
\ref{F11}\right) \ $we get $a_{11}a_{12}\neq 0$. Therefore, from $\left( \ref%
{F5}\right) $ and $\left( \ref{F6}\right) $, we get $a_{21}=a_{22}=0$. Thus,
Eq. $\left( \ref{F7}\right) $ gives $a_{44}=0$, which is a contradiction. So
any two of $W_{4},W_{5},W_{6}$ are not in the same orbit. Therefore we get
the following pairwise non-isomorphic algebras:

\begin{itemize}
\item $J_{5,19}:a\circ b=d,c^{2}=d,c\circ d=e.$

\item $J_{5,20}:a\circ b=d,c^{2}=d,c\circ d=e,a^{2}=e.$

\item $J_{5,21}:a\circ b=d,c^{2}=d,c\circ d=e,a^{2}=e,b^{2}=e.$
\end{itemize}

\subsection{1-dimensional annihilator extensions of $J_{4,7}$}

The automorphism group $Aut\left( J_{4,7}\right) $ consists of

\begin{equation*}
\phi =\left[ 
\begin{array}{cccc}
a_{11} & 0 & 0 & 0 \\ 
a_{21} & a_{11}^{2} & 0 & 0 \\ 
0 & 0 & a_{11}^{\frac{3}{2}} & 0 \\ 
a_{41} & 2a_{11}a_{21} & a_{43} & a_{11}^{3}%
\end{array}%
\right] :a_{11}\neq 0.
\end{equation*}%
Choose an arbitrary subspace $W\in U_{1}\left( J_{4,7}\right) $. From Table %
\ref{tabx}, such a subspace is spanned by $\left[ \theta \right] =\alpha
_{1}\left( \left[ \delta _{b,b}\right] +\left[ \delta _{a,d}\right] \right)
+\alpha _{2}\left[ \delta _{a,c}\right] +\alpha _{3}\left[ \delta _{a,b}%
\right] \ $such that $\alpha _{1}\neq 0$.$\allowbreak $ Let $\phi =\big(%
a_{ij}\big)\in Aut(J_{4,7})$. Then $\left[ \phi \theta \right] =\alpha
_{1}a_{11}^{4}\left( \left[ \delta _{b,b}\right] +\left[ \delta _{a,d}\right]
\right) +a_{11}\left( \alpha _{1}a_{43}+\alpha _{2}a_{11}^{\frac{3}{2}%
}\right) \left[ \delta _{a,c}\right] +a_{11}^{2}\left( 3\alpha
_{1}a_{21}+\alpha _{3}a_{11}\right) \left[ \delta _{a,b}\right] $. This
shows that the coefficient of $\delta _{a,b}$ depends on whether
characteristic $\mathbb{F=}3$ or not. So we distinguish two cases.

\begin{enumerate}
\item Suppose first that characteristic $\mathbb{F}\neq 3$\textbf{.} Let $%
\phi $ be the following automorphism:%
\begin{equation*}
\phi =%
\begin{bmatrix}
1 & 0 & 0 & 0 \\ 
-\frac{1}{3}\alpha _{1}^{-1}\alpha _{3} & 1 & 0 & 0 \\ 
0 & 0 & 1 & 0 \\ 
0 & 0 & -\alpha _{1}^{-1}\alpha _{2} & 1%
\end{bmatrix}%
.
\end{equation*}%
Then $\left[ \phi \theta \right] =\alpha _{1}\left( \left[ \delta _{b,b}%
\right] +\left[ \delta _{a,d}\right] \right) $. So we get a representative $%
W_{1}=\left\langle \left[ \delta _{b,b}\right] +\left[ \delta _{a,d}\right]
\right\rangle $. Therefore, if characteristic $\mathbb{F}\neq 3$\textbf{,} $%
Aut(J_{4,7})$ has only one orbit on $U_{1}\left( J_{4,7}\right) $.

\item Assume now that characteristic $\mathbb{F}=3$. Then $\left[ \phi
\theta \right] =\alpha _{1}a_{11}^{4}\left( \left[ \delta _{b,b}\right] +%
\left[ \delta _{a,d}\right] \right) +a_{11}\left( \alpha _{1}a_{43}+\alpha
_{2}a_{11}^{\frac{3}{2}}\right) \left[ \delta _{a,c}\right] +\alpha
_{3}a_{11}^{3}\left[ \delta _{a,b}\right] $. So $\mbox{Orb}\left(
\left\langle \left[ \theta \right] :\alpha _{3}=0\right\rangle \right) \cap %
\mbox{Orb}\left( \left\langle \left[ \theta \right] :\alpha _{3}\neq
0\right\rangle \right) =\emptyset $ and hence $Aut(J_{4,7})$ has at least
two orbits on $U_{1}\left( J_{4,7}\right) $. Suppose first that $\alpha
_{3}=0$. Let $\phi $ be the following automorphism:%
\begin{equation*}
\phi =%
\begin{bmatrix}
1 & 0 & 0 & 0 \\ 
0 & 1 & 0 & 0 \\ 
0 & 0 & 1 & 0 \\ 
0 & 0 & -\alpha _{1}^{-1}\alpha _{2} & 1%
\end{bmatrix}%
.
\end{equation*}%
Then $\left[ \phi \theta \right] =\alpha _{1}\left( \left[ \delta _{b,b}%
\right] +\left[ \delta _{a,d}\right] \right) $. So we get the same
representative as above. Assume now that $\alpha _{3}\neq 0$. Let $\phi $ be
the following automorphism:%
\begin{equation*}
\phi =%
\begin{bmatrix}
\alpha _{1}^{-1}\alpha _{3} & 0 & 0 & 0 \\ 
0 & \alpha _{1}^{-2}\alpha _{3}^{2} & 0 & 0 \\ 
0 & 0 & \alpha _{1}^{-\frac{3}{2}}\alpha _{3}^{\frac{3}{2}} & 0 \\ 
0 & 0 & \alpha _{1}^{-\frac{5}{2}}\alpha _{2}\alpha _{3}^{\frac{3}{2}} & 
\alpha _{1}^{-3}\alpha _{3}^{3}%
\end{bmatrix}%
.
\end{equation*}%
Then $\left[ \phi \theta \right] =\alpha _{1}^{-3}\alpha _{3}^{4}\left( %
\left[ \delta _{b,b}\right] +\left[ \delta _{a,d}\right] +\left[ \delta
_{a,b}\right] \right) $. So we get a representative $W_{2}=\left\langle %
\left[ \delta _{b,b}\right] +\left[ \delta _{a,d}\right] +\left[ \delta
_{a,b}\right] \right\rangle $. Therefore, if characteristic $\mathbb{F}=3$, $%
Aut(J_{4,7})$ has exactly two orbits on $U_{1}\left( J_{4,7}\right) $. So we
get the algebra:
\end{enumerate}

\begin{itemize}
\item $J_{5,22}:a^{2}=b,a\circ b=d,c^{2}=d,b^{2}=e,a\circ d=e.$

\item $M_{5,1}:a^{2}=b,a\circ b=d+e,c^{2}=d,b^{2}=e,a\circ d=e$.
Isomorphism: $M_{5,1}\cong J_{5,22}$ if and only if characteristic $\mathbb{F%
}\neq 3$. To describe the isomorphism, let $\varphi \left( a\right) =a+\frac{%
1}{3}b,\varphi \left( b\right) =b+\frac{2}{3}d+\frac{1}{9}e,\varphi \left(
c\right) =c,\varphi \left( d\right) =d,\varphi \left( e\right) =e$. Then $%
\varphi :$\ $M_{5,1}\longrightarrow J_{5,22}$ is an isomorphism.
\end{itemize}

\subsection{1-dimensional annihilator extensions of $J_{4,12}$}

The automorphism group $Aut\left( J_{4,12}\right) $ consists of

\begin{equation*}
\phi =\left[ 
\begin{array}{cccc}
a_{11} & 0 & 0 & 0 \\ 
a_{21} & a_{22} & 0 & 0 \\ 
a_{31} & a_{32} & a_{11}^{2} & 0 \\ 
a_{41} & a_{42} & 2a_{11}a_{21} & a_{11}a_{22}%
\end{array}%
\right] :a_{11}a_{22}\neq 0.
\end{equation*}%
Choose an arbitrary subspace $W\in U_{1}\left( J_{4,12}\right) $. From Table %
\ref{tabx}, such a subspace is spanned by $\left[ \theta \right] =\alpha _{1}%
\left[ \delta _{b,b}\right] +\alpha _{2}\left[ \delta _{a,c}\right] +\alpha
_{3}\left[ \delta _{b,c}\right] +\alpha _{4}\left[ \delta _{a,d}\right]
+\alpha _{5}\left[ \delta _{b,d}\right] $ such that $\left( \alpha
_{5},\alpha _{3}-\alpha _{4}\right) \neq \left( 0,0\right) $ and $\theta
^{\bot }\cap \left\langle c,d\right\rangle =0$. Moreover, $\theta ^{\bot
}\cap \left\langle c,d\right\rangle =0$ if and only if $\epsilon =\alpha
_{2}\alpha _{5}-\alpha _{3}\alpha _{4}\neq 0$. Let $\phi =\big(a_{ij}\big)%
\in Aut\left( J_{4,12}\right) $, and write $\left[ \phi \theta \right]
=\alpha _{1}^{\prime }\left[ \delta _{b,b}\right] +\alpha _{2}^{\prime }%
\left[ \delta _{a,c}\right] +\alpha _{3}^{\prime }\left[ \delta _{b,c}\right]
+\alpha _{4}^{\prime }\left[ \delta _{a,d}\right] +\alpha _{5}^{\prime }%
\left[ \delta _{b,d}\right] $. Then%
\begin{eqnarray*}
\alpha _{1}^{\prime } &=&a_{22}\left( \alpha _{1}a_{22}+2\alpha
_{3}a_{32}+2\alpha _{5}a_{42}\right) , \\
\alpha _{2}^{\prime } &=&a_{11}^{2}\left( \alpha _{2}a_{11}+\alpha
_{3}a_{21}\right) +2a_{11}a_{21}\left( \alpha _{4}a_{11}+\alpha
_{5}a_{21}\right) , \\
\alpha _{3}^{\prime } &=&a_{22}\left( \alpha _{3}a_{11}^{2}+2\alpha
_{5}a_{21}a_{11}\right) , \\
\alpha _{4}^{\prime } &=&a_{11}a_{22}\left( \alpha _{4}a_{11}+\alpha
_{5}a_{21}\right) , \\
\alpha _{5}^{\prime } &=&\alpha _{5}a_{11}a_{22}^{2}.
\end{eqnarray*}%
As $a_{11}a_{22}^{2}\neq 0$, the equation for $\alpha _{5}^{\prime }$
implies that $\alpha _{5}^{\prime }\neq 0$ if and only if $\alpha _{5}\neq 0$%
. Thus $\mbox{Orb}\left( \left\langle \left[ \theta \right] :\alpha _{5}\neq
0\right\rangle \right) \cap \mbox{Orb}\left( \left\langle \left[ \theta %
\right] :\alpha _{5}=0\right\rangle \right) =\emptyset $. Therefore we
distinguish two cases.

\textsc{Case 1. }$\alpha _{5}\neq 0$.

Let $\phi $ be the following automorphism:%
\begin{equation*}
\phi =\allowbreak 
\begin{bmatrix}
\epsilon ^{-\frac{1}{3}} & 0 & 0 & 0 \\ 
-\epsilon ^{-\frac{1}{3}}\alpha _{4}\alpha _{5}^{-1} & \epsilon ^{\frac{1}{6}%
}\alpha _{5}^{-\frac{1}{2}} & 0 & 0 \\ 
0 & 0 & \epsilon ^{-\frac{2}{3}} & 0 \\ 
0 & -\frac{1}{2}\epsilon ^{\frac{1}{6}}\alpha _{1}\alpha _{5}^{-\frac{3}{2}}
& -2\epsilon ^{-\frac{2}{3}}\alpha _{4}\alpha _{5}^{-1} & \epsilon ^{-\frac{1%
}{6}}\alpha _{5}^{-\frac{1}{2}}%
\end{bmatrix}%
.
\end{equation*}%
Then $\left[ \phi \theta \right] =\left[ \delta _{a,c}\right] +\epsilon ^{-%
\frac{1}{2}}\alpha _{5}^{-\frac{1}{2}}\left( \alpha _{3}-2\alpha _{4}\right) %
\left[ \delta _{b,c}\right] +\left[ \delta _{b,d}\right] $. Set $\alpha
=\epsilon ^{-\frac{1}{2}}\alpha _{5}^{-\frac{1}{2}}\left( \alpha
_{3}-2\alpha _{4}\right) $. Then we get the representatives $W_{1}^{\alpha
\in \mathbb{F}}=\left\langle \left[ \delta _{a,c}\right] +\alpha \left[
\delta _{b,c}\right] +\left[ \delta _{b,d}\right] \right\rangle $. Let us
now determine the possible orbits among the representatives $W_{1}^{\alpha
\in \mathbb{F}}$. We claim that the subspaces $W_{1}^{\alpha },W_{1}^{\beta
} $ are in the same orbit if and only if $\alpha ^{2}=\beta ^{2}$. To prove
this, assume first that $\phi W_{1}^{\alpha }=W_{1}^{\beta }$ for some $\phi
=\big(a_{ij}\big)\in Aut\left( J_{4,12}\right) $. Then there is a $\lambda
\in \mathbb{F}^{\ast }$ such that 
\begin{equation*}
\left( 
\begin{array}{c}
2a_{22}\left( a_{42}+\alpha a_{32}\right) \left[ \delta _{b,b}\right]
+a_{11}\left( a_{11}^{2}+\alpha a_{11}a_{21}+2a_{21}^{2}\right) \left[
\delta _{a,c}\right] \\ 
+a_{11}a_{22}\left( 2a_{21}+\alpha a_{11}\right) \left[ \delta _{b,c}\right]
+a_{11}a_{21}a_{22}\left[ \delta _{a,d}\right] +a_{11}a_{22}^{2}\left[
\delta _{b,d}\right]%
\end{array}%
\right) =\lambda \left( \left[ \delta _{a,c}\right] +\beta \left[ \delta
_{b,c}\right] +\left[ \delta _{b,d}\right] \right) .
\end{equation*}%
This amounts to the following polynomial equations:%
\begin{eqnarray}
2a_{22}\left( a_{42}+\alpha a_{32}\right) &=&0,  \notag \\
a_{11}\left( a_{11}^{2}+\alpha a_{11}a_{21}+2a_{21}^{2}\right) &=&\lambda ,
\label{H2} \\
a_{11}a_{22}\left( 2a_{21}+\alpha a_{11}\right) &=&\lambda \beta ,
\label{H3} \\
a_{11}a_{21}a_{22} &=&0,  \label{H4} \\
a_{11}a_{22}^{2} &=&\lambda .  \label{H5}
\end{eqnarray}%
Since $a_{11}a_{22}\neq 0$, Eq. $\left( \ref{H4}\right) $ gives $a_{21}=0$.
This allows us to rewrite Eqs. $\left( \ref{H2}\right) $ and $\left( \ref{H3}%
\right) $ as follows: 
\begin{eqnarray}
a_{11}^{3} &=&\lambda ,  \label{H6} \\
a_{11}^{2}a_{22}\alpha &=&\lambda \beta .  \label{H7}
\end{eqnarray}%
Using Eqs.$\ \left( \ref{H6}\right) ,\left( \ref{H7}\right) $ and $\left( %
\ref{H5}\right) $, one obtains that%
\begin{equation*}
\lambda ^{2}\beta ^{2}=\left( a_{22}a_{11}^{2}\alpha \right) ^{2}=\left(
a_{11}^{3}\right) \left( a_{11}a_{22}^{2}\right) \alpha ^{2}=\lambda
^{2}\alpha ^{2}.
\end{equation*}%
Consequently, $\alpha ^{2}=\beta ^{2}$. Conversely, suppose that $\alpha
^{2}=\beta ^{2}$. When $\alpha =\beta $, there is nothing to prove. Assume
now that $\alpha =-\beta $. Let $\phi $ be the diagonal matrix with the
entries $(1,-1,1,-1)$ in the diagonal. Then $\phi W_{1}^{\alpha
}=W_{1}^{-\alpha }$, this completes the proof of the claim. So we get the
following algebras:

\begin{itemize}
\item $J_{5,23}^{\alpha \in \mathbb{F}}:a^{2}=c,a\circ b=d,a\circ c=e,b\circ
d=e,b\circ c=\alpha e$. Isomorphism: $J_{5,23}^{\alpha }\cong
J_{5,23}^{\beta }$ if and only if $\beta =$ $\alpha $ or $\beta =-\alpha $.
In the latter case, $\varphi (a)=a,\varphi (b)=-b,\varphi (c)=c,\varphi
(d)=-d,\varphi (e)=e$, defines an isomorphism $\varphi :J_{5,23}^{\alpha
}\longrightarrow J_{5,23}^{-\alpha }$.
\end{itemize}

\textsc{Case 2. }$\alpha _{5}=0.$

In this case $W=\left\langle \left[ \theta \right] :\alpha
_{5}=0\right\rangle \in U_{1}\left( J_{4,12}\right) $ if and only if $\alpha
_{3}\alpha _{4}\left( \alpha _{3}-\alpha _{4}\right) \neq 0$. Let us now
consider a few cases:

\begin{enumerate}
\item Suppose first that $\alpha _{3}\neq -2\alpha _{4}$. Let $\phi $ be the
following automorphism:%
\begin{equation*}
\phi =%
\begin{bmatrix}
1 & 0 & 0 & 0 \\ 
-\alpha _{2}\left( \alpha _{3}+2\alpha _{4}\right) ^{-1} & \alpha _{3}^{-1}
& 0 & 0 \\ 
0 & -\frac{1}{2}\alpha _{1}\alpha _{3}^{-2} & 1 & 0 \\ 
0 & 0 & -2\alpha _{2}\left( \alpha _{3}+2\alpha _{4}\right) ^{-1} & \alpha
_{3}^{-1}%
\end{bmatrix}%
.
\end{equation*}%
Then $\left[ \phi \theta \right] =\left[ \delta _{b,c}\right] +\alpha
_{3}^{-1}\alpha _{4}\left[ \delta _{a,d}\right] $. Set $\alpha =\alpha
_{3}^{-1}\alpha _{4}$, then $\alpha \notin \left\{ 0,1,-\frac{1}{2}\right\} $%
. So we get the representatives $W_{2}^{\alpha \in \mathbb{F}^{\ast }\mathbb{%
-}\left\{ 1,-\frac{1}{2}\right\} }=\left\langle \left[ \delta _{b,c}\right]
+\alpha \left[ \delta _{a,d}\right] \right\rangle $.

\item Suppose now that $\alpha _{3}=-2\alpha _{4}$. We distinguish two cases:

\begin{enumerate}
\item $\alpha _{2}=0$. Let $\phi $ be the following automorphism:%
\begin{equation*}
\phi =%
\begin{bmatrix}
1 & 0 & 0 & 0 \\ 
0 & \alpha _{3}^{-1} & 0 & 0 \\ 
0 & -\frac{1}{2}\alpha _{1}\alpha _{3}^{-2} & 1 & 0 \\ 
0 & 0 & 0 & \alpha _{3}^{-1}%
\end{bmatrix}%
.
\end{equation*}%
Then $\left[ \phi \theta \right] =\left[ \delta _{b,c}\right] -\frac{1}{2}%
\left[ \delta _{a,d}\right] $. So we get a representative $W_{2}^{\alpha =-%
\frac{1}{2}}$.

\item $\alpha _{2}\neq 0$. Let $\phi $ be the following automorphism:%
\begin{equation*}
\phi =%
\begin{bmatrix}
\alpha _{2}^{-\frac{1}{3}} & 0 & 0 & 0 \\ 
0 & \alpha _{2}^{\frac{2}{3}}\alpha _{4}^{-1} & 0 & 0 \\ 
0 & \frac{1}{4}\alpha _{1}\alpha _{2}^{\frac{2}{3}}\alpha _{4}^{-2} & \alpha
_{2}^{-\frac{2}{3}} & 0 \\ 
0 & 0 & 0 & \alpha _{2}^{\frac{1}{3}}\alpha _{4}^{-1}%
\end{bmatrix}%
.
\end{equation*}%
Then $\left[ \phi \theta \right] =\left[ \delta _{a,c}\right] -2\left[
\delta _{b,c}\right] +\left[ \delta _{a,d}\right] $. Therefore we have a
representative $W_{3}=\left\langle \left[ \delta _{a,c}\right] -2\left[
\delta _{b,c}\right] +\left[ \delta _{a,d}\right] \right\rangle $.
\end{enumerate}
\end{enumerate}

As shown we have the representatives $W_{2}^{\alpha \in \mathbb{F}^{\ast }%
\mathbb{-}\left\{ 1\right\} }=\left\langle \delta _{b,c}+\alpha \delta
_{a,d}\right\rangle ,W_{3}=\left\langle \delta _{a,c}-2\delta _{b,c}+\delta
_{a,d}\right\rangle $. Let us determine the possible orbits among such
representatives. First we claim that $\mbox{Orb}\left( W_{2}^{\alpha
}\right) \cap \mbox{Orb}\left( W_{3}\right) =\emptyset $ for all $\alpha \in 
\mathbb{F}^{\ast }\mathbb{-}\left\{ 1\right\} $. To prove this, suppose on
the contrary that $\phi W_{2}^{\alpha }=W_{3}$ for some $\phi =\big(a_{ij}%
\big)\in Aut\left( J_{4,12}\right) $. Then there is a $\lambda \in \mathbb{F}%
^{\ast }$ such that%
\begin{equation*}
2a_{22}a_{32}\left[ \delta _{b,b}\right] +a_{11}^{2}a_{21}\left( 2\alpha
+1\right) \left[ \delta _{a,c}\right] +a_{11}^{2}a_{22}\left[ \delta _{b,c}%
\right] +\alpha a_{11}^{2}a_{22}\left[ \delta _{a,d}\right] =\lambda \left( %
\left[ \delta _{a,c}\right] -2\left[ \delta _{b,c}\right] +\left[ \delta
_{a,d}\right] \right) .
\end{equation*}%
This amounts to the following polynomial equations:%
\begin{eqnarray}
2a_{22}a_{32} &=&0,  \notag \\
a_{11}^{2}a_{21}\left( 2\alpha +1\right) &=&\lambda ,  \label{H8} \\
a_{11}^{2}a_{22} &=&-2\lambda ,  \label{H9} \\
\alpha a_{11}^{2}a_{22} &=&\lambda .  \label{H10}
\end{eqnarray}%
Substitute Eq. $\left( \ref{H9}\right) $\ into Eq. $\left( \ref{H10}\right) $%
, we get $\alpha =-\frac{1}{2}$. So, from Eq. $\left( \ref{H8}\right) $ we
get $\lambda =0$, which is a contradiction. Therefore, for all $\alpha \in 
\mathbb{F}^{\ast }\mathbb{-}\left\{ 1\right\} $, $\mbox{Orb}\left(
W_{2}^{\alpha }\right) \cap \mbox{Orb}\left( W_{3}\right) =\emptyset $. It
remains to determine the possible orbits among the representatives $%
W_{2}^{\alpha \in \mathbb{F}^{\ast }\mathbb{-}\left\{ 1\right\} }$. We claim
that $\mbox{Orb}\left( W_{2}^{\alpha }\right) =\mbox{Orb}\left( W_{2}^{\beta
}\right) $ if and only if $\alpha =\beta $. To prove this, assume that $\phi
W_{2}^{\alpha }=W_{2}^{\beta }$ for some $\phi =\big(a_{ij}\big)\in
Aut\left( J_{4,12}\right) $. Then there is a $\lambda \in \mathbb{F}^{\ast }$
such that%
\begin{equation*}
2a_{22}a_{32}\left[ \delta _{b,b}\right] +a_{11}^{2}a_{21}\left( 2\alpha
+1\right) \left[ \delta _{a,c}\right] +a_{11}^{2}a_{22}\left[ \delta _{b,c}%
\right] +\alpha a_{11}^{2}a_{22}\left[ \delta _{a,d}\right] =\lambda \left( %
\left[ \delta _{b,c}\right] +\beta \left[ \delta _{a,d}\right] \right) .
\end{equation*}%
This amounts to the following polynomial equations:%
\begin{eqnarray}
2a_{22}a_{32} &=&0,  \notag \\
a_{11}^{2}a_{21}\left( 2\alpha +1\right) &=&0,  \notag \\
a_{11}^{2}a_{22} &=&\lambda ,  \label{H11} \\
\alpha a_{11}^{2}a_{22} &=&\lambda \beta .  \label{H12}
\end{eqnarray}%
Substitute Eq. $\left( \ref{H11}\right) $ into Eq. $\left( \ref{H12}\right) $
to obtain $\alpha =\beta $, as claimed. So we get the following algebras:

\begin{itemize}
\item $J_{5,24}^{\alpha \in \mathbb{F-}\left\{ 0,1\right\} }:a^{2}=c,a\circ
b=d,b\circ c=e,a\circ d=\alpha e$. Isomorphism: $J_{5,24}^{\alpha }\cong
J_{5,24}^{\beta }$ if and only if $\alpha =\beta $.

\item $J_{5,25}:a^{2}=c,a\circ b=d,a\circ c=e,b\circ c=-2e,a\circ
d=e.\allowbreak $
\end{itemize}

\subsection{1-dimensional annihilator extensions of $J_{4,13}$}

The automorphism group $Aut\left( J_{4,13}\right) $ consists of

\begin{equation}
\phi =\left[ 
\begin{array}{cccc}
a_{11} & a_{12} & 0 & 0 \\ 
a_{21} & a_{22} & 0 & 0 \\ 
a_{31} & a_{32} & a_{11}^{2} & a_{12}^{2} \\ 
a_{41} & a_{42} & a_{21}^{2} & a_{22}^{2}%
\end{array}%
\right] :a_{11}a_{22}\neq 0,a_{12}=a_{21}=0\mbox{ or }a_{12}a_{21}\neq
0,a_{11}=a_{22}=0.  \label{AutJ4,13}
\end{equation}%
Choose an arbitrary subspace $W\in U_{1}\left( J_{4,13}\right) $. From Table %
\ref{tabx}, such a subspace is spanned by $\left[ \theta \right] =\alpha _{1}%
\left[ \delta _{a,b}\right] +\alpha _{2}\left[ \delta _{a,c}\right] +\alpha
_{3}\left[ \delta _{b,d}\right] +\alpha _{4}\left[ \delta _{b,c}\right]
+\alpha _{5}\left[ \delta _{a,d}\right] $ such that $\left( \alpha
_{4},\alpha _{5}\right) \neq \left( 0,0\right) $ and $\theta ^{\bot }\cap
\left\langle c,d\right\rangle =0$. Moreover, $\theta ^{\bot }\cap
\left\langle c,d\right\rangle =0$ if and only if $\alpha _{2}\alpha
_{3}-\alpha _{4}\alpha _{5}\neq 0$. Let $\phi =\big(a_{ij}\big)\in Aut\left(
J_{4,13}\right) $, and write $\left[ \phi \theta \right] =\alpha
_{1}^{\prime }\left[ \delta _{a,b}\right] +\alpha _{2}^{\prime }\left[
\delta _{a,c}\right] +\alpha _{3}^{\prime }\left[ \delta _{b,d}\right]
+\alpha _{4}^{\prime }\left[ \delta _{b,c}\right] +\alpha _{5}^{\prime }%
\left[ \delta _{a,d}\right] $. Then%
\begin{eqnarray*}
\alpha _{1}^{\prime } &=&a_{11}\left( \alpha _{1}a_{22}+\alpha
_{2}a_{32}+\alpha _{5}a_{42}\right) +a_{21}\left( \alpha _{1}a_{12}+\alpha
_{3}a_{42}+\alpha _{4}a_{32}\right) \\
&&+a_{31}\left( \alpha _{2}a_{12}+\alpha _{4}a_{22}\right) +a_{41}\left(
\alpha _{3}a_{22}+\alpha _{5}a_{12}\right) , \\
\alpha _{2}^{\prime } &=&\alpha _{2}a_{11}^{2}a_{11}+\alpha
_{3}a_{21}^{2}a_{21}, \\
\alpha _{3}^{\prime } &=&\alpha _{2}a_{12}^{2}a_{12}+\alpha
_{3}a_{22}^{2}a_{22}, \\
\alpha _{4}^{\prime } &=&\alpha _{4}a_{11}^{2}a_{22}+\alpha
_{5}a_{21}^{2}a_{12}, \\
\alpha _{5}^{\prime } &=&\alpha _{4}a_{12}^{2}\allowbreak a_{21}+\alpha
_{5}a_{22}^{2}a_{11}\allowbreak .
\end{eqnarray*}%
Let us now consider some case distinctions.

\begin{enumerate}
\item If $\alpha _{4}=0$ and $\alpha _{5}\neq 0$, then $\alpha _{2}\alpha
_{3}\neq 0$. Let $\phi $ be the following automorphism:%
\begin{equation*}
\phi =\left[ 
\begin{array}{cccc}
\alpha _{2}^{-\frac{1}{3}} & 0 & 0 & 0 \\ 
0 & \alpha _{2}^{\frac{1}{6}}\alpha _{5}^{-\frac{1}{2}} & 0 & 0 \\ 
0 & 0 & \alpha _{2}^{-\frac{2}{3}} & 0 \\ 
0 & -\alpha _{1}\alpha _{2}^{\frac{1}{6}}\alpha _{5}^{-\frac{3}{2}} & 0 & 
\alpha _{2}^{\frac{1}{3}}\alpha _{5}^{-1}%
\end{array}%
\right] .
\end{equation*}%
Then $\left[ \phi \theta \right] =\left[ \delta _{a,c}\right] +\alpha _{2}^{%
\frac{1}{2}}\alpha _{3}\alpha _{5}^{-\frac{3}{2}}\left[ \delta _{b,d}\right]
+\left[ \delta _{a,d}\right] $. Set $\alpha =\alpha _{2}^{\frac{1}{2}}\alpha
_{3}\alpha _{5}^{-\frac{3}{2}}$, then $\alpha \in \mathbb{F}^{\ast }$. So we
get the representatives $W_{1}^{\alpha \in \mathbb{F}^{\ast }}=\left\langle %
\left[ \delta _{a,c}\right] +\alpha \left[ \delta _{b,d}\right] +\left[
\delta _{a,d}\right] \right\rangle $.

\item If $\alpha _{4}\neq 0$ and $\alpha _{5}=0$, then $\alpha _{2}\alpha
_{3}\neq 0$. Let $\phi $ be the following automorphism:%
\begin{equation*}
\phi =\left[ 
\begin{array}{cccc}
0 & \alpha _{3}^{\frac{1}{6}}\alpha _{4}^{-\frac{1}{2}} & 0 & 0 \\ 
\alpha _{3}^{-\frac{1}{3}} & 0 & 0 & 0 \\ 
0 & -\alpha _{1}\alpha _{3}^{\frac{1}{6}}\alpha _{4}^{-\frac{3}{2}} & 0 & 
\alpha _{3}^{\frac{1}{3}}\alpha _{4}^{-1} \\ 
0 & 0 & \alpha _{3}^{-\frac{2}{3}} & 0%
\end{array}%
\right] .
\end{equation*}%
Then $\left[ \phi \theta \right] =\left[ \delta _{a,c}\right] +\alpha
_{2}\alpha _{3}^{\frac{1}{2}}\alpha _{4}^{-\frac{3}{2}}\left[ \delta _{b,d}%
\right] +\left[ \delta _{a,d}\right] $. So $\phi W=W_{1}^{\alpha }$, where $%
\alpha =\alpha _{2}\alpha _{3}^{\frac{1}{2}}\alpha _{4}^{-\frac{3}{2}}\in 
\mathbb{F}^{\ast }$.

\item If $\alpha _{4}\alpha _{5}\neq 0$, then $\alpha _{2}\alpha _{3}\neq
\alpha _{4}\alpha _{5}$. Let $\phi $ be as follows:%
\begin{equation*}
\phi =\left[ 
\begin{array}{cccc}
\alpha _{4}^{-\frac{2}{3}}\alpha _{5}^{\frac{1}{3}} & 0 & 0 & 0 \\ 
0 & \alpha _{4}^{\frac{1}{3}}\alpha _{5}^{-\frac{2}{3}} & 0 & 0 \\ 
0 & 0 & \alpha _{4}^{-\frac{4}{3}}\alpha _{5}^{\frac{2}{3}} & 0 \\ 
0 & -\alpha _{1}\alpha _{4}^{\frac{1}{3}}\alpha _{5}^{-\frac{5}{3}} & 0 & 
\alpha _{4}^{\frac{2}{3}}\alpha _{5}^{-\frac{4}{3}}%
\end{array}%
\right] .
\end{equation*}%
Then $\left[ \phi \theta \right] =\alpha _{2}\alpha _{4}^{-2}\alpha _{5}%
\left[ \delta _{a,c}\right] +\alpha _{3}\alpha _{4}\alpha _{5}^{-2}\left[
\delta _{b,d}\right] +\left[ \delta _{b,c}\right] +\left[ \delta _{a,d}%
\right] $. Set $\alpha =\alpha _{2}\alpha _{4}^{-2}\alpha _{5},\beta =\alpha
_{3}\alpha _{4}\alpha _{5}^{-2}$, then $\alpha ,\beta \in \mathbb{F}$ and $%
\alpha \beta \neq 1$. So we get the representatives $W_{2}^{\alpha ,\beta
\in \mathbb{F}}=\left\langle \alpha \left[ \delta _{a,c}\right] +\beta \left[
\delta _{b,d}\right] +\left[ \delta _{b,c}\right] +\left[ \delta _{a,d}%
\right] :\alpha \beta \neq 1\right\rangle $.
\end{enumerate}

From the above we have the representatives $W_{1}^{\alpha \in \mathbb{F}%
^{\ast }}=\left\langle \left[ \delta _{a,c}\right] +\alpha \left[ \delta
_{b,d}\right] +\left[ \delta _{a,d}\right] \right\rangle ,W_{2}^{\alpha
,\beta \in \mathbb{F}}=\left\langle \alpha \left[ \delta _{a,c}\right]
+\beta \left[ \delta _{b,d}\right] +\left[ \delta _{b,c}\right] +\left[
\delta _{a,d}\right] :\alpha \beta \neq 1\right\rangle $. Let us now
determine the possible orbits among such representatives. First we claim
that $\mbox{Orb}\left( W^{\alpha ^{\prime }}\right) \cap \mbox{Orb}\left(
W_{2}^{\alpha ,\beta }\right) =\emptyset $ for all $\alpha ^{\prime },\alpha
,\beta \in \mathbb{F}$. To prove this, suppose on the contrary that $\phi
W^{\alpha ^{\prime }}=W_{2}^{\alpha ,\beta }$ for some $\phi =\big(a_{ij}%
\big)\in Aut(J_{4,13})$. Then there is a $\lambda \in \mathbb{F}^{\ast }$
such that:%
\begin{equation}
\phi \left( \left[ \delta _{a,c}\right] +\alpha ^{\prime }\left[ \delta
_{b,d}\right] +\left[ \delta _{a,d}\right] \right) =\lambda \left( \alpha %
\left[ \delta _{a,c}\right] +\beta \left[ \delta _{b,d}\right] +\left[
\delta _{b,c}\right] +\left[ \delta _{a,d}\right] \right) .  \label{CJ4,13}
\end{equation}%
Using the equations for $\alpha _{4}^{\prime }$ and $\alpha _{5}^{\prime }$
to compare coefficients of $\delta _{b,c}$ and $\delta _{a,d}$ on the two
sides of $\left( \ref{CJ4,13}\right) $, we get%
\begin{equation*}
a_{12}a_{21}^{2}=\lambda ,a_{11}a_{22}^{2}=\lambda .
\end{equation*}%
From which we get $a_{12}a_{21}^{2}=a_{11}a_{22}^{2}\neq 0$, which
contradicts with $\left( \ref{AutJ4,13}\right) $. Hence $\mbox{Orb}\left(
W^{\alpha ^{\prime }}\right) \cap \mbox{Orb}\left( W_{2}^{\alpha ,\beta
}\right) =\emptyset $ for all $\alpha ^{\prime },\alpha ,\beta \in \mathbb{F}
$, as claimed. Next, we claim that $\mbox{Orb}\left( W_{1}^{\alpha }\right) =%
\mbox{Orb}\left( W_{1}^{\beta }\right) $ if and only if $\alpha ^{2}=\beta
^{2}$. To prove this, suppose first that $\phi W_{1}^{\alpha }=W_{1}^{\beta
} $ for some $\phi =\big(a_{ij}\big)\in Aut\left( J_{4,13}\right) $. Then
there is a $\lambda \in \mathbb{F}^{\ast }$ such that $\phi \left( \left[
\delta _{a,c}\right] +\alpha \left[ \delta _{b,d}\right] +\left[ \delta
_{a,d}\right] \right) =\lambda \left( \left[ \delta _{a,c}\right] +\beta %
\left[ \delta _{b,d}\right] +\left[ \delta _{a,d}\right] \right) $. This
then amounts to the following polynomial equations:%
\begin{eqnarray}
a_{11}\left( a_{32}+a_{42}\right) +\alpha
a_{21}a_{42}+a_{31}a_{12}+a_{41}\left( \alpha a_{22}+a_{12}\right) &=&0, 
\notag \\
a_{11}^{3}+\alpha a_{21}^{3} &=&\lambda ,  \label{eq2} \\
a_{12}^{3}+\alpha a_{22}^{3} &=&\lambda \beta ,  \label{eq3} \\
a_{21}^{3} &=&0,  \label{eq4} \\
a_{11}a_{22}^{2} &=&\lambda .  \label{eq5}
\end{eqnarray}%
From Eq. $\left( \ref{eq4}\right) $ we get $a_{21}=0$. So, from $\left( \ref%
{CJ4,13}\right) $ we get $a_{12}=0$. Consequently, we can rewrite Eqs. $%
\left( \ref{eq2}\right) ,\left( \ref{eq3}\right) $ and $\left( \ref{eq5}%
\right) $ as follows:%
\begin{eqnarray}
a_{11}^{3} &=&\lambda ,  \label{q1} \\
a_{22}^{3}\alpha &=&\lambda \beta ,  \label{q2} \\
a_{22}^{2}a_{11}\allowbreak &=&\lambda .  \label{q3}
\end{eqnarray}%
From $\left( \ref{q1}\right) $ and $\left( \ref{q3}\right) $, we get $%
\lambda ^{2}=a_{22}^{6}$. So, from Eq. $\left( \ref{q2}\right) $ we get $%
\alpha ^{2}=\beta ^{2}$. Conversely, suppose that $\alpha ^{2}=\beta ^{2}$.
If $\alpha =\beta $, we have nothing to prove. Suppose that $\alpha =-\beta $%
. Choose $\phi $ to be the diagonal matrix with the entries $(1,-1,1,1)$ in
the diagonal. Then $\phi W_{1}^{\alpha }=W_{1}^{-\alpha }$, this completes
the proof of the claim. Finally, we claim that $\mbox{Orb}\left(
W_{2}^{\alpha ,\beta }\right) =\mbox{Orb}\left( W_{2}^{\alpha ^{\prime
},\beta ^{\prime }}\right) $ if and only if $\left( \alpha ,\beta \right)
=\left( \alpha ^{\prime },\beta ^{\prime }\right) $ or $\left( \alpha ,\beta
\right) =\left( \beta ^{\prime },\alpha ^{\prime }\right) $. To prove this,
suppose first that $\phi W_{2}^{\alpha ,\beta }=W_{2}^{\alpha ^{\prime
},\beta ^{\prime }}$ for some $\phi =\big(a_{ij}\big)\in Aut(J_{4,13})$.
Then there is a $\lambda \in \mathbb{F}^{\ast }$ such that 
\begin{equation}
\phi \left( \alpha \left[ \delta _{a,c}\right] +\beta \left[ \delta _{b,d}%
\right] +\left[ \delta _{b,c}\right] +\left[ \delta _{a,d}\right] \right)
=\lambda \left( \alpha ^{\prime }\left[ \delta _{a,c}\right] +\beta ^{\prime
}\left[ \delta _{b,d}\right] +\left[ \delta _{b,c}\right] +\left[ \delta
_{a,d}\right] \right) .  \label{m}
\end{equation}%
From $\left( \ref{AutJ4,13}\right) $, we have either $a_{12}=a_{21}=0$ or $%
a_{11}=a_{22}=0$. Suppose first that $a_{12}=a_{21}=0$. Then $\left( \ref{m}%
\right) $ amounts to the following polynomial equations:%
\begin{eqnarray}
a_{11}\left( \alpha a_{32}+a_{42}\right) +a_{31}a_{22}+\beta a_{41}a_{22}
&=&0,  \notag \\
a_{11}^{3}\alpha &=&\lambda \alpha ^{\prime },  \label{m2} \\
a_{22}^{3}\beta &=&\lambda \beta ^{\prime },  \label{m3} \\
a_{11}^{2}a_{22} &=&\lambda ,  \label{m4} \\
a_{22}^{2}a_{11} &=&\lambda .  \label{m5}
\end{eqnarray}%
Equations $\left( \ref{m4}\right) $ and $\left( \ref{m5}\right) $ give that $%
a_{11}^{3}=a_{22}^{3}=\lambda $. Whence, Eq. $\left( \ref{m2}\right) $\
gives $\alpha =\alpha ^{\prime }$, while Eq. $\left( \ref{m3}\right) $ gives 
$\beta =\beta ^{\prime }$. Suppose now that $a_{11}=a_{22}=0$. Then $\left( %
\ref{m}\right) $ amounts to the following polynomial equations:%
\begin{eqnarray}
a_{21}\left( \beta a_{42}+a_{32}\right) +\alpha a_{31}a_{12}+a_{41}a_{12}
&=&0,  \notag \\
a_{21}^{3}\beta &=&\lambda \alpha ^{\prime },  \label{n2} \\
a_{12}^{3}\alpha &=&\lambda \beta ^{\prime },  \label{n3} \\
a_{21}^{2}a_{12} &=&\lambda ,  \label{n4} \\
a_{12}^{2}a_{21} &=&\lambda .  \label{n5}
\end{eqnarray}%
Equations $\left( \ref{n4}\right) $ and $\left( \ref{n5}\right) $ give that $%
a_{12}^{3}=a_{21}^{3}=\lambda $. Whence, Eq. $\left( \ref{n2}\right) $\
gives $\beta =\alpha ^{\prime }$, while Eq. $\left( \ref{n3}\right) $ gives $%
\alpha =\beta ^{\prime }$. Conversely, suppose that $\left( \alpha ,\beta
\right) =\left( \beta ,\alpha \right) $. Let $\phi $ be as follows:%
\begin{equation*}
\phi =\left[ 
\begin{array}{cccc}
0 & 1 & 0 & 0 \\ 
1 & 0 & 0 & 0 \\ 
0 & 0 & 0 & 1 \\ 
0 & 0 & 1 & 0%
\end{array}%
\right] .
\end{equation*}%
Then $\phi W_{2}^{\alpha ,\beta }=W_{2}^{\beta ,\alpha }$, this completes
the proof of the claim. Therefore we get the following algebras:

\begin{itemize}
\item $J_{5,26}^{\alpha \in \mathbb{F}^{\ast }}:a^{2}=c,b^{2}=d,a\circ
c=e,a\circ d=e,b\circ d=\alpha e$. Isomorphism: $J_{5,26}^{\alpha }\cong
J_{5,26}^{\beta }$ if and only if $\beta =$ $\alpha $ or $\beta =-\alpha $.
To describe the isomorphism in the latter case, let $\varphi (a)=a,\varphi
(b)=-b,\varphi (c)=c,\varphi (d)=d,\varphi (e)=e$. Then $\varphi
:J_{5,26}^{\alpha }\longrightarrow J_{5,26}^{-\alpha }$ is an isomorphism.

\item $J_{5,27}^{\alpha ,\beta \in \mathbb{F}}\left( \alpha \beta \neq
1\right) :a^{2}=c,b^{2}=d,b\circ c=e,a\circ d=e,b\circ d=\alpha e,a\circ
c=\beta e$. Isomorphism: $J_{5,27}^{\alpha ,\beta }\cong J_{5,27}^{\alpha
^{\prime },\beta ^{\prime }}$ if and only if $\left( \alpha ,\beta \right)
=\left( \alpha ^{\prime },\beta ^{\prime }\right) $ or $\left( \alpha ,\beta
\right) =\left( \beta ^{\prime },\alpha ^{\prime }\right) $. To describe the
isomorphism in the latter case, let $\varphi (a)=b,\varphi (b)=a,\varphi
(c)=d,\varphi (d)=c,\varphi (e)=e$. Then $\varphi :J_{5,27}^{\alpha ,\beta
}\longrightarrow J_{5,27}^{\beta ,\alpha }$ is an isomorphism.
\end{itemize}

\subsection{2-dimensional annihilator extensions of $J_{3,2}$}

Choose an arbitrary subspace $W\in U_{2}\left( J_{3,2}\right) $. From Table %
\ref{tabx}, such a subspace is spanned by $\left[ \theta _{1}\right] =\alpha
_{1}\left[ \delta _{a,b}\right] +\alpha _{2}\left[ \delta _{a,c}\right]
+\alpha _{3}\left[ \delta _{b,c}\right] +\alpha _{4}\left[ \delta _{c,c}%
\right] $ and $\left[ \theta _{2}\right] =\beta _{1}\left[ \delta _{a,b}%
\right] +\beta _{2}\left[ \delta _{a,c}\right] +\beta _{3}\left[ \delta
_{b,c}\right] +\beta _{4}\left[ \delta _{c,c}\right] $ such that $\left(
\alpha _{3},\beta _{3}\right) \neq \left( 0,0\right) $. Without loss of
generality we may assume that $\left( \alpha _{3},\beta _{3}\right) =\left(
1,0\right) $. Assume now that $W=\left\langle \alpha _{1}\left[ \delta _{a,b}%
\right] +\alpha _{2}\left[ \delta _{a,c}\right] +\left[ \delta _{b,c}\right]
+\alpha _{4}\left[ \delta _{c,c}\right] ,\beta _{1}\left[ \delta _{a,b}%
\right] +\beta _{2}\left[ \delta _{a,c}\right] +\beta _{4}\left[ \delta
_{c,c}\right] \right\rangle $. Let $\phi =\big(a_{ij}\big)\in $ $Aut\left(
J_{3,2}\right) $. Then $\phi W=\left\langle \alpha _{1}^{\prime }\left[
\delta _{a,b}\right] +\alpha _{2}^{\prime }\left[ \delta _{a,c}\right]
+\alpha _{3}^{\prime }\left[ \delta _{b,c}\right] +\alpha _{4}^{\prime }%
\left[ \delta _{c,c}\right] ,\beta _{1}^{\prime }\left[ \delta _{a,b}\right]
+\beta _{2}^{\prime }\left[ \delta _{a,c}\right] +\beta _{4}^{\prime }\left[
\delta _{c,c}\right] \right\rangle $ where%
\begin{eqnarray*}
\alpha _{1}^{\prime } &=&a_{11}^{3}\alpha _{1}+a_{31}a_{11}^{2}, \\
\alpha _{2}^{\prime } &=&a_{11}a_{23}\alpha _{1}+a_{11}a_{33}\alpha
_{2}+a_{21}a_{33}+a_{31}a_{23}+a_{31}a_{33}\alpha _{4}, \\
\alpha _{3}^{\prime } &=&a_{11}^{2}a_{33}, \\
\alpha _{4}^{\prime } &=&2a_{23}a_{33}+a_{33}^{2}\alpha _{4}, \\
\beta _{1}^{\prime } &=&a_{11}^{3}\beta _{1}, \\
\beta _{2}^{\prime } &=&a_{11}a_{23}\beta _{1}+a_{11}a_{33}\beta
_{2}+a_{31}a_{33}\beta _{4}, \\
\beta _{4}^{\prime } &=&a_{33}^{2}\beta _{4}.
\end{eqnarray*}%
We claim that there is a $\phi \in Aut\left( J_{3,2}\right) $ such that $%
\phi W=\left\langle \left[ \delta _{b,c}\right] ,\beta _{1}^{\prime }\left[
\delta _{a,b}\right] +\beta _{2}^{\prime }\left[ \delta _{a,c}\right] +\beta
_{4}^{\prime }\left[ \delta _{c,c}\right] \right\rangle $. To see this, let $%
\phi $\ to be the following automorphism: 
\begin{equation*}
\phi =\left[ 
\begin{array}{ccc}
1 & 0 & 0 \\ 
\alpha _{1}\alpha _{4}-\alpha _{2} & 1 & -\frac{1}{2}\alpha _{4} \\ 
-\alpha _{1} & 0 & 1%
\end{array}%
\right] .
\end{equation*}%
Then $\phi W=\left\langle \left[ \delta _{b,c}\right] ,\beta _{1}^{\prime }%
\left[ \delta _{a,b}\right] +\beta _{2}^{\prime }\left[ \delta _{a,c}\right]
+\beta _{4}^{\prime }\left[ \delta _{c,c}\right] \right\rangle $. Thus we
may now assume without loss of generality that $W=\left\langle \left[ \delta
_{b,c}\right] ,\beta _{1}\left[ \delta _{a,b}\right] +\beta _{2}\left[
\delta _{a,c}\right] +\beta _{4}\left[ \delta _{c,c}\right] \right\rangle $.
Let us choose $\phi =\big(a_{ij}\big)\in Aut\left( J_{3,2}\right) $ with $%
a_{31}=a_{23}=a_{21}=0$. Then 
\begin{equation}
\phi W=\left\langle \left[ \delta _{b,c}\right] ,a_{11}^{3}\beta _{1}\left[
\delta _{a,b}\right] +a_{11}a_{33}\beta _{2}\left[ \delta _{a,c}\right]
+a_{33}^{2}\beta _{4}\left[ \delta _{c,c}\right] \right\rangle .
\label{J3,2}
\end{equation}%
We next consider the following cases:

\begin{enumerate}
\item If $\beta _{1}=\beta _{2}=0,\beta _{4}\neq 0$, we then have $\phi
W=\left\langle \left[ \delta _{b,c}\right] ,\left[ \delta _{c,c}\right]
\right\rangle $. Thus we get a representative $W_{1}=\left\langle \left[
\delta _{b,c}\right] ,\left[ \delta _{c,c}\right] \right\rangle $.

\item If $\beta _{1}=\beta _{4}=0,\beta _{2}\neq 0$, we then have $\phi
W=\left\langle \left[ \delta _{b,c}\right] ,\left[ \delta _{a,c}\right]
\right\rangle $. So we get a representative $W_{2}=\left\langle \left[
\delta _{b,c}\right] ,\left[ \delta _{a,c}\right] \right\rangle $.

\item If $\beta _{1}=0,\beta _{2}\beta _{4}\neq 0$, we put $a_{11}=1$ and $%
a_{33}=\beta _{2}\beta _{4}^{-1}$ into $\left( \ref{J3,2}\right) $. Then $%
\phi W=\left\langle \left[ \delta _{b,c}\right] ,\left[ \delta _{a,c}\right]
+\left[ \delta _{c,c}\right] \right\rangle $. Hence we get a representative $%
W_{3}=\left\langle \left[ \delta _{b,c}\right] ,\left[ \delta _{a,c}\right] +%
\left[ \delta _{c,c}\right] \right\rangle $.

\item If $\beta _{1}\neq 0,\beta _{2}=\beta _{4}=0$, we then have $\phi
W=\left\langle \left[ \delta _{b,c}\right] ,\left[ \delta _{a,b}\right]
\right\rangle $. Hence we get a representative $W_{4}=\left\langle \left[
\delta _{b,c}\right] ,\left[ \delta _{a,b}\right] \right\rangle .$

\item If $\beta _{1}\beta _{2}\neq 0,\beta _{4}=0$, we put $a_{11}=1$ and $%
a_{33}=\beta _{2}^{-1}\beta _{1}$ into $\left( \ref{J3,2}\right) $. Then $%
\phi W=\left\langle \left[ \delta _{b,c}\right] ,\left[ \delta _{a,b}\right]
+\left[ \delta _{a,c}\right] \right\rangle $. Thus we get a representative $%
W_{5}=\left\langle \left[ \delta _{b,c}\right] ,\left[ \delta _{a,b}\right] +%
\left[ \delta _{a,c}\right] \right\rangle $.

\item If $\beta _{1}\beta _{4}\neq 0,\beta _{2}=0$, we put $a_{11}=1$ and $%
a_{33}=\sqrt{\beta _{1}\beta _{4}^{-1}}$ into $\left( \ref{J3,2}\right) $.
Then $\phi W=\left\langle \left[ \delta _{b,c}\right] ,\left[ \delta _{a,b}%
\right] +\left[ \delta _{c,c}\right] \right\rangle $. So we get a
representative $W_{6}=\left\langle \left[ \delta _{b,c}\right] ,\left[
\delta _{a,b}\right] +\left[ \delta _{c,c}\right] \right\rangle $.

\item If $\beta _{1}\beta _{2}\beta _{4}\neq 0$, we put $a_{11}=\beta
_{1}^{-1}\beta _{2}^{2}\beta _{4}^{-1}$ and $a_{33}=\beta _{1}^{-1}\beta
_{2}^{3}\beta _{4}^{-2}$ into $\left( \ref{J3,2}\right) $. Then $\phi
W=\left\langle \left[ \delta _{b,c}\right] ,\left[ \delta _{a,b}\right] +%
\left[ \delta _{a,c}\right] +\left[ \delta _{c,c}\right] \right\rangle $.
Thus we get a representative $W_{7}=\left\langle \left[ \delta _{b,c}\right]
,\left[ \delta _{a,b}\right] +\left[ \delta _{a,c}\right] +\left[ \delta
_{c,c}\right] \right\rangle $.
\end{enumerate}

Let us now determine the possible orbits among the representatives $%
W_{1},\ldots ,W_{7}$. First note that the algebras corresponding to the
representatives $W_{1},W_{2},W_{3}$\ are of nilpotency type $\left(
2,2,1\right) $ while those corresponding to the representatives $%
W_{4},W_{5},W_{6},W_{7}$ are of nilpotency type $\left( 2,1,2\right) $.
Therefore $\mbox{Orb}\left( X\right) \cap \mbox{Orb}\left( Y\right)
=\emptyset $ for $X\in \left\{ W_{1},W_{2},W_{3}\right\} $ and $Y\in \left\{
W_{4},W_{5},W_{6},W_{7}\right\} $. Let us first determine the possible
orbits among the representatives $W_{1},W_{2},W_{3}$. Since $\Psi \left(
W_{1}\right) =\left( 2,1\right) $ and $\Psi \left( W_{2}\right) =\Psi \left(
W_{3}\right) =\left( 1,1\right) $, $\mbox{Orb}\left( W_{1}\right) \cap %
\mbox{Orb}\left( W_{i}\right) =\emptyset $ for $i=2,3$. Further we claim
that $\mbox{Orb}\left( W_{2}\right) \cap \mbox{Orb}\left( W_{3}\right)
=\emptyset $. To see this, consider any $\phi =\big(a_{ij}\big)\in Aut\left(
J_{3,2}\right) $. Then $\phi W_{2}=\left\langle \left[ \phi \delta _{b,c}%
\right] ,a_{11}a_{33}\left[ \delta _{a,c}\right] \right\rangle $. Since $%
\left[ \delta _{a,c}\right] \in \phi W_{2}$ and $\left[ \delta _{a,c}\right]
\notin W_{3}$, $\phi W_{2}\neq W_{3}$ as claimed. It remains to determine
the possible orbits among the representatives $W_{4},W_{5},W_{6},W_{7}$.
First we claim that $\mbox{Orb}\left( W_{4}\right) \cap \mbox{Orb}\left(
W_{i}\right) =\emptyset $ for $i=5,6,7$. To see this, let $\phi =\big(a_{ij}%
\big)\in Aut\left( J_{3,2}\right) $. Then $\phi W_{4}=\left\langle \left[
\phi \delta _{b,c}\right] ,\left[ \phi \delta _{a,b}\right] \right\rangle $
where%
\begin{eqnarray*}
\left[ \phi \delta _{b,c}\right] &=&a_{31}a_{11}^{2}\left[ \delta _{a,b}%
\right] +\left( a_{21}a_{33}+a_{31}a_{23}\right) \left[ \delta _{a,c}\right]
+a_{11}^{2}a_{33}\left[ \delta _{b,c}\right] +2a_{23}a_{33}\left[ \delta
_{c,c}\right] , \\
\left[ \phi \delta _{a,b}\right] &=&a_{11}^{3}\left[ \delta _{a,b}\right]
+a_{11}a_{23}\left[ \delta _{a,c}\right] .
\end{eqnarray*}%
Now if $a_{23}=0$ then $\left[ \phi \delta _{a,b}\right] \notin W_{i}$ for $%
i=5,6,7$. Further if $a_{23}\neq 0$ then $\left[ \phi \delta _{b,c}\right]
\notin W_{5}$ and $\left[ \phi \delta _{a,b}\right] \notin W_{i}$ for $i=6,7$%
. This proves our claim. Next we claim that $\mbox{Orb}\left( W_{5}\right)
\cap \mbox{Orb}\left( W_{i}\right) =\emptyset $ for $i=6,7$. Consider any $%
\phi =\big(a_{ij}\big)\in Aut\left( J_{3,2}\right) $. Then 
\begin{equation*}
\phi W_{5}=\left\langle \left[ \phi \delta _{b,c}\right] ,a_{11}^{3}\left[
\delta _{a,b}\right] +a_{11}\left( a_{23}+a_{33}\right) \left[ \delta _{a,c}%
\right] \right\rangle .
\end{equation*}%
This shows that $\phi W_{5}\neq W_{i}$ ($i=6,7$) since $a_{11}^{3}\delta
_{a,b}+a_{11}\left( a_{23}+a_{33}\right) \delta _{a,c}$ neither belongs to $%
W_{6}$ nor $W_{7}$. So $\mbox{Orb}\left( W_{5}\right) \cap \mbox{Orb}\left(
W_{i}\right) =\emptyset $ for $i=6,7$. Thus we only need to decide whether $%
W_{6},W_{7}$ are in the same orbit or not. We claim that $W_{6},W_{7}$ are
in the same orbit if and only if characteristic $\mathbb{F}\neq 3$. To prove
this, assume first that $\phi W_{6}=W_{7}$ for some $\phi =\big(a_{ij}\big)%
\in Aut\left( J_{3,2}\right) $. Then $\left\langle \left[ \phi \delta _{b,c}%
\right] ,\left[ \phi \delta _{a,b}\right] +\left[ \phi \delta _{c,c}\right]
\right\rangle =\left\langle \left[ \delta _{b,c}\right] ,\left[ \delta _{a,b}%
\right] +\left[ \delta _{a,c}\right] +\left[ \delta _{c,c}\right]
\right\rangle $. Moreover,%
\begin{eqnarray*}
\left[ \phi \delta _{b,c}\right] &=&a_{11}^{2}a_{31}\left[ \delta _{a,b}%
\right] +\left( a_{21}a_{33}+a_{31}a_{23}\right) \left[ \delta _{a,c}\right]
+a_{11}^{2}a_{33}\left[ \delta _{b,c}\right] +2a_{23}a_{33}\left[ \delta
_{c,c}\right] , \\
\left[ \phi \delta _{a,b}\right] +\left[ \phi \delta _{c,c}\right]
&=&a_{11}^{3}\left[ \delta _{a,b}\right] +\left(
a_{11}a_{23}+a_{31}a_{33}\right) \left[ \delta _{a,c}\right] +a_{33}^{2}%
\left[ \delta _{c,c}\right] .
\end{eqnarray*}%
So there is an invertible matrix $\big(b_{ij}\big)$ such that%
\begin{eqnarray*}
a_{11}^{2}a_{31}\left[ \delta _{a,b}\right] +\left(
a_{21}a_{33}+a_{31}a_{23}\right) \left[ \delta _{a,c}\right]
+a_{11}^{2}a_{33}\left[ \delta _{b,c}\right] +2a_{23}a_{33}\left[ \delta
_{c,c}\right] &=&b_{11}\left[ \delta _{b,c}\right] +b_{21}\left( \left[
\delta _{a,b}\right] +\left[ \delta _{a,c}\right] +\left[ \delta _{c,c}%
\right] \right) , \\
a_{11}^{3}\left[ \delta _{a,b}\right] +\left(
a_{11}a_{23}+a_{31}a_{33}\right) \left[ \delta _{a,c}\right] +a_{33}^{2}%
\left[ \delta _{c,c}\right] &=&b_{21}\left[ \delta _{b,c}\right]
+b_{22}\left( \left[ \delta _{a,b}\right] +\left[ \delta _{a,c}\right] +%
\left[ \delta _{c,c}\right] \right) .
\end{eqnarray*}%
By equating the coefficients of like terms on the two sides, we obtain%
\begin{equation*}
a_{11}^{2}a_{31}=\left( a_{21}a_{33}+a_{31}a_{23}\right)
=2a_{23}a_{33},a_{11}^{3}=\left( a_{11}a_{23}+a_{31}a_{33}\right)
=a_{33}^{2}.
\end{equation*}%
Consequently,%
\begin{equation*}
a_{11}^{5}=a_{11}^{2}\left( a_{11}a_{23}+a_{31}a_{33}\right) =\left(
a_{11}^{3}\right) a_{23}+\left( a_{11}^{2}a_{31}\right)
a_{33}=a_{33}^{2}a_{23}+2a_{23}a_{33}^{2}=3a_{23}a_{33}^{2}.
\end{equation*}%
Since $a_{11}\neq 0$, it follows that characteristic $\mathbb{F}\neq 3$.
Conversely, suppose that characteristic $\mathbb{F}\neq 3$. Let $\phi \in
Aut\left( J_{3,2}\right) $ be such that%
\begin{equation*}
\phi =\frac{1}{9}%
\begin{bmatrix}
9 & 0 & 0 \\ 
4 & 9 & 3 \\ 
6 & 0 & 9%
\end{bmatrix}%
.
\end{equation*}%
Then $\phi W_{6}=\left\langle \left[ \delta _{b,c}\right] +\frac{2}{3}\left[
\delta _{a,b}+\delta _{a,c}+\delta _{c,c}\right] ,\left[ \delta
_{a,b}+\delta _{a,c}+\delta _{c,c}\right] \right\rangle =$ $\left\langle %
\left[ \delta _{b,c}\right] ,\left[ \delta _{a,b}+\delta _{a,c}+\delta _{c,c}%
\right] \right\rangle =W_{7}$, as claimed. Therefore we get the following
algebras:

\begin{itemize}
\item $J_{5,28}:a^{2}=b,b\circ c=d,c^{2}=e.$

\item $J_{5,29}:a^{2}=b,b\circ c=d,a\circ c=e.$

\item $J_{5,30}:a^{2}=b,b\circ c=d,a\circ c=e,c^{2}=e.$

\item $J_{5,31}:a^{2}=b,b\circ c=d,a\circ b=e.$

\item $J_{5,32}:a^{2}=b,b\circ c=d,a\circ b=e,a\circ c=e.$

\item $J_{5,33}:a^{2}=b,b\circ c=d,a\circ b=e,c^{2}=e.$

\item $M_{5,2}:a^{2}=b,b\circ c=d,a\circ b=e,a\circ c=e,c^{2}=e$.
Isomorphism: $M_{5,2}\cong J_{5,33}$ if and only if characteristic $\mathbb{%
F\neq }3$. To describe the isomorphism, let $\varphi (a)=a+\frac{4}{9}b+%
\frac{2}{3}c,\varphi (b)=b+\frac{16}{27}d+\frac{4}{3}e,\varphi (c)=\frac{1}{3%
}b+c,\varphi (d)=d,\varphi (e)=\frac{2}{3}d+e$. Then $\varphi
:M_{5,2}\longrightarrow J_{5,33}$ is an isomorphism.
\end{itemize}

\subsection{2-dimensional annihilator extensions of $J_{3,3}$}

Choose an arbitrary subspace $W\in U_{2}\left( J_{3,2}\right) $. From Table %
\ref{tabx}, such a subspace is spanned by $\theta _{1}=\alpha _{1}\delta
_{a,a}+\alpha _{2}\delta _{b,b}+\alpha _{3}\delta _{a,c}+\alpha _{4}\delta
_{b,c}$ and $\theta _{2}=\beta _{1}\delta _{a,a}+\beta _{2}\delta
_{b,b}+\beta _{3}\delta _{a,c}+\beta _{4}\delta _{b,c}$ such that $\theta
_{1}^{\bot }\cap \theta _{2}^{\bot }\cap \left\langle c\right\rangle =0$.
Moreover, $\theta _{1}^{\bot }\cap \theta _{2}^{\bot }\cap \left\langle
c\right\rangle =0$ if and only if $\left( \left( \alpha _{3},\alpha
_{4}\right) ,\left( \beta _{3},\beta _{4}\right) \right) \neq \left( \left(
0,0\right) ,\left( 0,0\right) \right) $. By possibly swapping $\theta _{1}$
and $\theta _{2}$, we may assume without loss of generality that $\left(
\alpha _{3},\alpha _{4}\right) \neq \left( 0,0\right) $. Let $\phi =\big(%
a_{ij}\big)\in $ $Aut\left( J_{3,3}\right) $, then $\phi W=\left\langle
\alpha _{1}^{\prime }\left[ \delta _{a,a}\right] +\alpha _{2}^{\prime }\left[
\delta _{b,b}\right] +\alpha _{3}^{\prime }\left[ \delta _{a,c}\right]
+\alpha _{4}^{\prime }\left[ \delta _{b,c}\right] ,\beta _{1}^{\prime }\left[
\delta _{a,a}\right] +\beta _{2}^{\prime }\left[ \delta _{b,b}\right] +\beta
_{3}^{\prime }\left[ \delta _{a,c}\right] +\beta _{4}^{\prime }\left[ \delta
_{b,c}\right] \right\rangle $ where%
\begin{eqnarray*}
\alpha _{1}^{\prime } &=&a_{11}^{2}\alpha _{1}+a_{21}^{2}\alpha
_{2}+2a_{11}a_{31}\alpha _{3}+2a_{21}a_{31}\alpha _{4}, \\
\alpha _{2}^{\prime } &=&a_{12}^{2}\alpha _{1}+a_{22}^{2}\alpha
_{2}+2a_{12}a_{32}\alpha _{3}+2a_{22}a_{32}\alpha _{4}, \\
\alpha _{3}^{\prime } &=&a_{11}^{2}a_{22}\alpha _{3}+a_{21}^{2}a_{12}\alpha
_{4}, \\
\alpha _{4}^{\prime } &=&a_{11}a_{22}^{2}\alpha _{4}+a_{21}a_{12}^{2}\alpha
_{3}, \\
\beta _{1}^{\prime } &=&a_{11}^{2}\beta _{1}+a_{21}^{2}\beta
_{2}+2a_{11}a_{31}\beta _{3}+2a_{21}a_{31}\beta _{4}, \\
\beta _{2}^{\prime } &=&a_{12}^{2}\beta _{1}+a_{22}^{2}\beta
_{2}+2a_{12}a_{32}\beta _{3}+2a_{22}a_{32}\beta _{4}, \\
\beta _{3}^{\prime } &=&a_{11}^{2}a_{22}\beta _{3}+a_{21}^{2}a_{12}\beta
_{4}, \\
\beta _{4}^{\prime } &=&a_{11}a_{22}^{2}\beta _{4}+a_{21}a_{12}^{2}\beta
_{3}.
\end{eqnarray*}%
Since $\left( \alpha _{3},\alpha _{4}\right) \neq \left( 0,0\right) $, we
may choose the $a_{ij}$ in such a way that $\alpha _{3}^{\prime }=1$. Thus
no generality is lost by assuming $W=\left\langle \alpha _{1}\left[ \delta
_{a,a}\right] +\alpha _{2}\left[ \delta _{b,b}\right] +\left[ \delta _{a,c}%
\right] +\alpha _{4}\left[ \delta _{b,c}\right] ,\beta _{1}\left[ \delta
_{a,a}\right] +\beta _{2}\left[ \delta _{b,b}\right] +\beta _{4}\left[
\delta _{b,c}\right] \right\rangle $. Further we claim that there is a $\phi
\in Aut\left( J_{3,2}\right) $ such that $\phi W\in \left\{ \left\langle %
\left[ \delta _{a,c}\right] ,\left[ \theta \right] \right\rangle
,\left\langle \left[ \delta _{b,b}\right] +\left[ \delta _{a,c}\right] ,%
\left[ \theta \right] \right\rangle ,\left\langle \left[ \delta _{a,c}\right]
+\left[ \delta _{b,c}\right] ,\left[ \theta \right] \right\rangle \right\} $
where $\left[ \theta \right] =\beta _{1}^{\prime }\left[ \delta _{a,a}\right]
+\beta _{2}^{\prime }\left[ \delta _{b,b}\right] +\beta _{4}^{\prime }\left[
\delta _{b,c}\right] $. To prove this, suppose first that $\alpha _{4}=0$.
Let $\phi $\ denote the first of the following matrices if $\alpha _{2}=0$,
or the second if $\alpha _{2}\neq 0$:%
\begin{equation*}
\begin{bmatrix}
1 & 0 & 0 \\ 
0 & 1 & 0 \\ 
-\frac{1}{2}\alpha _{1} & 0 & 1%
\end{bmatrix}%
,%
\begin{bmatrix}
1 & 0 & 0 \\ 
0 & \alpha _{2}^{-1} & 0 \\ 
-\frac{1}{2}\alpha _{1} & 0 & \alpha _{2}^{-1}%
\end{bmatrix}%
.
\end{equation*}%
Then $\phi W=\left\langle \left[ \delta _{a,c}\right] ,\left[ \theta \right]
\right\rangle $ if $\alpha _{2}=0$, while $\phi W=\left\langle \left[ \delta
_{b,b}\right] +\left[ \delta _{a,c}\right] ,\left[ \theta \right]
\right\rangle $ otherwise. Assume now that $\alpha _{4}\neq 0$. Let $\phi $
be as follows: 
\begin{equation*}
\phi =\left[ 
\begin{array}{ccc}
\alpha _{4} & 0 & 0 \\ 
0 & 1 & 0 \\ 
-\frac{1}{2}\alpha _{1}\alpha _{4} & -\frac{1}{2}\alpha _{2}\alpha _{4}^{-1}
& \alpha _{4}%
\end{array}%
\right] .
\end{equation*}%
Then $\phi W=\left\langle \left[ \delta _{a,c}\right] +\left[ \delta _{b,c}%
\right] ,\left[ \theta \right] \right\rangle $. This completes the proof of
the claim. Let us hence suppose without loss of generality that $W\in
\left\{ \left\langle \left[ \delta _{a,c}\right] ,\left[ \theta \right]
\right\rangle ,\left\langle \left[ \delta _{b,b}\right] +\left[ \delta _{a,c}%
\right] ,\left[ \theta \right] \right\rangle ,\left\langle \left[ \delta
_{a,c}\right] +\left[ \delta _{b,c}\right] ,\left[ \theta \right]
\right\rangle \right\} $ where $\left[ \theta \right] =\beta _{1}\left[
\delta _{a,a}\right] +\beta _{2}\left[ \delta _{b,b}\right] +\beta _{4}\left[
\delta _{b,c}\right] $. Let us now consider the following cases:

\textsc{Case 1.} $W=\left\langle \left[ \delta _{a,c}\right] ,\beta _{1}%
\left[ \delta _{a,a}\right] +\beta _{2}\left[ \delta _{b,b}\right] +\beta
_{4}\left[ \delta _{b,c}\right] \right\rangle $.

\begin{enumerate}
\item If $\beta _{1}\neq 0,\beta _{2}=\beta _{4}=0$, then $W=\left\langle %
\left[ \delta _{a,c}\right] ,\left[ \delta _{a,a}\right] \right\rangle $. So
we have a representative $W_{1}=\left\langle \left[ \delta _{a,c}\right] ,%
\left[ \delta _{a,a}\right] \right\rangle $.

\item If $\beta _{2}\neq 0,\beta _{1}=\beta _{4}=0$, then $W=\left\langle %
\left[ \delta _{a,c}\right] ,\left[ \delta _{b,b}\right] \right\rangle $. So
we have a representative $W_{2}=\left\langle \left[ \delta _{a,c}\right] ,%
\left[ \delta _{b,b}\right] \right\rangle $.

\item If $\beta _{4}\neq 0,\beta _{1}=\beta _{2}=0$, then $W=\left\langle %
\left[ \delta _{a,c}\right] ,\left[ \delta _{b,c}\right] \right\rangle $. So
we have a representative $W_{3}=\left\langle \left[ \delta _{a,c}\right] ,%
\left[ \delta _{b,c}\right] \right\rangle $.

\item If $\beta _{1}\beta _{2}\neq 0,\beta _{4}=0$, we have%
\begin{equation*}
\phi =%
\begin{bmatrix}
\sqrt{\beta _{1}^{-1}} & 0 & 0 \\ 
0 & \sqrt{\beta _{2}^{-1}} & 0 \\ 
0 & 0 & \sqrt{\beta _{1}^{-1}\beta _{2}^{-1}}%
\end{bmatrix}%
\in Aut\left( J_{3,3}\right) \mbox{ and }\phi W=\left\langle \left[ \delta
_{a,c}\right] ,\left[ \delta _{a,a}\right] +\left[ \delta _{b,b}\right]
\right\rangle .
\end{equation*}%
So we get a representative $W_{4}=\left\langle \left[ \delta _{a,c}\right] ,%
\left[ \delta _{a,a}\right] +\left[ \delta _{b,b}\right] \right\rangle $.

\item If $\beta _{1}\beta _{4}\neq 0,\beta _{2}=0$, we have%
\begin{equation*}
\phi =%
\begin{bmatrix}
\beta _{1}^{-1}\beta _{4} & 0 & 0 \\ 
0 & 1 & 0 \\ 
0 & 0 & \beta _{1}^{-1}\beta _{4}%
\end{bmatrix}%
\in Aut\left( J_{3,3}\right) \mbox{ and }\phi W=\left\langle \left[ \delta
_{a,c}\right] ,\left[ \delta _{a,a}\right] +\left[ \delta _{b,c}\right]
\right\rangle .
\end{equation*}%
So we get a representative $W_{5}=\left\langle \left[ \delta _{a,c}\right] ,%
\left[ \delta _{a,a}\right] +\left[ \delta _{b,c}\right] \right\rangle $.

\item If $\beta _{2}\beta _{4}\neq 0,\beta _{1}=0$, we have%
\begin{equation*}
\phi =%
\begin{bmatrix}
\beta _{4}^{-1} & 0 & 0 \\ 
0 & 1 & 0 \\ 
0 & -\frac{1}{2}\beta _{2}\beta _{4}^{-1} & \beta _{4}^{-1}%
\end{bmatrix}%
\in Aut\left( J_{3,3}\right) \mbox{ and }\phi W=W_{3}.
\end{equation*}

\item If $\beta _{1}\beta _{2}\beta _{4}\neq 0$, we have%
\begin{equation*}
\phi =%
\begin{bmatrix}
\beta _{1}^{-1}\beta _{4} & 0 & 0 \\ 
0 & 1 & 0 \\ 
0 & -\frac{1}{2}\beta _{2}\beta _{4}^{-1} & \beta _{1}^{-1}\beta _{4}%
\end{bmatrix}%
\in Aut\left( J_{3,3}\right) \mbox{ and }\phi W=W_{5}.
\end{equation*}
\end{enumerate}

\textsc{Case 2.} $W=\left\langle \left[ \delta _{b,b}\right] +\left[ \delta
_{a,c}\right] ,\beta _{1}\left[ \delta _{a,a}\right] +\beta _{2}\left[
\delta _{b,b}\right] +\beta _{4}\left[ \delta _{b,c}\right] \right\rangle $.

\begin{enumerate}
\item If $\beta _{1}\neq 0,\beta _{2}=\beta _{4}=0$, then $W=\left\langle %
\left[ \delta _{b,b}\right] +\left[ \delta _{a,c}\right] ,\left[ \delta
_{a,a}\right] \right\rangle $. So we have a representative $%
W_{6}=\left\langle \left[ \delta _{b,b}\right] +\left[ \delta _{a,c}\right] ,%
\left[ \delta _{a,a}\right] \right\rangle $.

\item If $\beta _{2}\neq 0,\beta _{1}=\beta _{4}=0$, then $W=W_{2}$.

\item If $\beta _{4}\neq 0,\beta _{1}=\beta _{2}=0$, then $W=\left\langle %
\left[ \delta _{b,b}\right] +\left[ \delta _{a,c}\right] ,\left[ \delta
_{b,c}\right] \right\rangle $. Let $\phi $ be as follows: 
\begin{equation*}
\phi =\left[ 
\begin{array}{ccc}
0 & 1 & 0 \\ 
1 & 0 & 0 \\ 
0 & 0 & 1%
\end{array}%
\right] .
\end{equation*}%
Then $\phi W=W_{5}$.

\item If $\beta _{1}\beta _{2}\neq 0,\beta _{4}=0$, then%
\begin{equation*}
\phi =%
\begin{bmatrix}
\sqrt{\beta _{1}\beta _{2}^{-1}} & 0 & 0 \\ 
0 & \beta _{1}\beta _{2}^{-1} & 0 \\ 
\frac{1}{2}\beta _{1}\beta _{2}^{-1}\sqrt{\beta _{1}\beta _{2}^{-1}} & 0 & 
\beta _{1}\beta _{2}^{-1}\sqrt{\beta _{1}\beta _{2}^{-1}}%
\end{bmatrix}%
\in Aut\left( J_{3,3}\right) \mbox{ and }\phi W=W_{4}.
\end{equation*}

\item If $\beta _{1}\beta _{4}\neq 0,\beta _{2}=0$, we have%
\begin{equation*}
\phi =\left[ 
\begin{array}{ccc}
\sqrt[3]{\beta _{1}\beta _{4}^{-1}} & 0 & 0 \\ 
0 & \sqrt[3]{\beta _{1}^{2}\beta _{4}^{-2}} & 0 \\ 
0 & 0 & \beta _{1}\beta _{4}^{-1}%
\end{array}%
\right] \in Aut\left( J_{3,3}\right) \mbox{ and }\phi W=\left\langle \left[
\delta _{b,b}\right] +\left[ \delta _{a,c}\right] ,\left[ \delta _{a,a}%
\right] +\left[ \delta _{b,c}\right] \right\rangle .
\end{equation*}%
So we get a representative $W_{7}=\left\langle \left[ \delta _{b,b}\right] +%
\left[ \delta _{a,c}\right] ,\left[ \delta _{a,a}\right] +\left[ \delta
_{b,c}\right] \right\rangle $.

\item If $\beta _{2}\beta _{4}\neq 0,\beta _{1}=0$, we have%
\begin{equation*}
\phi =\left[ 
\begin{array}{ccc}
0 & 1 & 0 \\ 
1 & 0 & 0 \\ 
-\frac{1}{2}\beta _{2}\beta _{4}^{-1} & 0 & 1%
\end{array}%
\right] \in Aut\left( J_{3,3}\right) \mbox{ and }\phi W=W_{5}.
\end{equation*}

\item If $\beta _{1}\beta _{2}\beta _{4}\neq 0$, we have%
\begin{equation*}
\phi =\left[ 
\begin{array}{ccc}
\sqrt[3]{\beta _{1}\beta _{4}^{-1}} & 0 & 0 \\ 
0 & \sqrt[3]{\beta _{1}^{2}\beta _{4}^{-2}} & 0 \\ 
0 & -\frac{1}{2}\beta _{2}\beta _{4}^{-1}\sqrt[3]{\beta _{1}^{2}\beta
_{4}^{-2}} & \beta _{1}\beta _{4}^{-1}%
\end{array}%
\right] \in Aut\left( J_{3,3}\right) \mbox{ and }\phi W=W_{7}.
\end{equation*}
\end{enumerate}

$\allowbreak $\textsc{Case 3.} $W=\left\langle \left[ \delta _{a,c}\right] +%
\left[ \delta _{b,c}\right] ,\beta _{1}\left[ \delta _{a,a}\right] +\beta
_{2}\left[ \delta _{b,b}\right] +\beta _{4}\left[ \delta _{b,c}\right]
\right\rangle $.

\begin{enumerate}
\item Suppose first that $\beta _{4}\neq 0$.

\begin{enumerate}
\item If $\beta _{1}\neq 0,\beta _{2}=0$, we have%
\begin{equation*}
\phi =\left[ 
\begin{array}{ccc}
\beta _{1}\beta _{4}^{-1} & 0 & 0 \\ 
0 & \beta _{1}\beta _{4}^{-1} & 0 \\ 
\frac{1}{2}\beta _{1}^{2}\beta _{4}^{-2} & 0 & \beta _{1}^{2}\beta _{4}^{-2}%
\end{array}%
\right] \in Aut\left( J_{3,3}\right) \mbox{ and }\phi W=W_{5}.
\end{equation*}

\item If $\beta _{1}=0,\beta _{2}\neq 0$, we have%
\begin{equation*}
\phi =\left[ 
\begin{array}{ccc}
0 & -\beta _{2}\beta _{4}^{-1} & 0 \\ 
-\beta _{2}\beta _{4}^{-1} & 0 & 0 \\ 
\frac{1}{2}\beta _{2}^{2}\beta _{4}^{-2} & 0 & \beta _{2}^{2}\beta _{4}^{-2}%
\end{array}%
\right] \in Aut\left( J_{3,3}\right) \mbox{ and }\phi W=W_{5}.
\end{equation*}

\item If $\beta _{1}\beta _{2}\neq 0$, we have%
\begin{equation*}
\phi =\left[ 
\begin{array}{ccc}
\beta _{4}^{-1}\sqrt[3]{\beta _{1}\beta _{2}^{2}} & 0 & 0 \\ 
0 & -\beta _{4}^{-1}\sqrt[3]{\beta _{1}^{2}\beta _{2}} & 0 \\ 
\frac{1}{2}\beta _{4}^{-2}\sqrt[3]{\beta _{1}^{4}\beta _{2}^{2}} & \frac{1}{2%
}\beta _{4}^{-2}\sqrt[3]{\beta _{1}^{2}\beta _{2}^{4}} & -\beta _{1}\beta
_{2}\beta _{4}^{-2}%
\end{array}%
\right] \in Aut\left( J_{3,3}\right) \mbox{ and }\phi W=W_{7}.
\end{equation*}

\item If $\beta _{1}=\beta _{2}=0$, then $W=W_{3}$.
\end{enumerate}

\item Suppose now that $\beta _{4}=0$.

\begin{enumerate}
\item If $\beta _{1}\neq 0$, we have%
\begin{equation*}
\phi =\left[ 
\begin{array}{ccc}
\sqrt{\beta _{1}^{-1}} & 0 & 0 \\ 
0 & \sqrt{\beta _{1}^{-1}} & 0 \\ 
0 & 0 & \beta _{1}^{-1}%
\end{array}%
\right] \in Aut\left( J_{3,3}\right) \text{.}
\end{equation*}%
Moreover, $\phi W=\left\langle \left[ \delta _{a,c}\right] +\left[ \delta
_{b,c}\right] ,\left[ \delta _{a,a}\right] +\alpha \left[ \delta _{b,b}%
\right] \right\rangle \mbox{ where }\alpha \in \mathbb{F}$. So we get the
representative $W_{8}^{\alpha \in \mathbb{F}}=\left\langle \left[ \delta
_{a,c}\right] +\left[ \delta _{b,c}\right] ,\left[ \delta _{a,a}\right]
+\alpha \left[ \delta _{b,b}\right] \right\rangle $.

\item If $\beta _{1}=0$, we have%
\begin{equation*}
\phi =\left[ 
\begin{array}{ccc}
0 & \sqrt{\beta _{2}^{-1}} & 0 \\ 
\sqrt{\beta _{2}^{-1}} & 0 & 0 \\ 
0 & 0 & \beta _{2}^{-1}%
\end{array}%
\right] \in Aut\left( J_{3,3}\right) \mbox{ and }\phi W=W_{8}^{\alpha =0}%
\text{.}
\end{equation*}
\end{enumerate}
\end{enumerate}

To summarize, we have the following representatives:%
\begin{table}[H] \centering%
$%
\begin{tabular}{|c|c|c|}
\hline
$W$ & $\Psi \left( W\right) $ & Type of corresponding algebra \\ \hline
$W_{1}=\left\langle \left[ \delta _{a,c}\right] ,\left[ \delta _{a,a}\right]
\right\rangle $ & $\left( 2,1\right) $ & $\left( 2,2,1\right) $ \\ \hline
$W_{2}=\left\langle \left[ \delta _{a,c}\right] ,\left[ \delta _{b,b}\right]
\right\rangle $ & $\left( 2,1\right) $ & $\left( 2,2,1\right) $ \\ \hline
$W_{3}=\left\langle \left[ \delta _{a,c}\right] ,\left[ \delta _{b,c}\right]
\right\rangle $ & $\left( 1,1\right) $ & $\left( 2,1,2\right) $ \\ \hline
$W_{4}=\left\langle \left[ \delta _{a,c}\right] ,\left[ \delta _{a,a}\right]
+\left[ \delta _{b,b}\right] \right\rangle $ & $\left( 1,1\right) $ & $%
\left( 2,2,1\right) $ \\ \hline
$W_{5}=\left\langle \left[ \delta _{a,c}\right] ,\left[ \delta _{a,a}\right]
+\left[ \delta _{b,c}\right] \right\rangle $ & $\left( 1,1\right) $ & $%
\left( 2,1,2\right) $ \\ \hline
$W_{6}=\left\langle \left[ \delta _{b,b}\right] +\left[ \delta _{a,c}\right]
,\left[ \delta _{a,a}\right] \right\rangle $ & $\left( 2,0\right) $ & $%
\left( 2,2,1\right) $ \\ \hline
$W_{7}=\left\langle \left[ \delta _{b,b}\right] +\left[ \delta _{a,c}\right]
,\left[ \delta _{a,a}\right] +\left[ \delta _{b,c}\right] \right\rangle $ & $%
\left( 1,1\right) $ & $\left( 2,1,2\right) $ \\ \hline
$W_{8}^{\alpha =0}=\left\langle \left[ \delta _{a,c}\right] +\left[ \delta
_{b,c}\right] ,\left[ \delta _{a,a}\right] \right\rangle $ & $\left(
2,1\right) $ & $\left( 2,2,1\right) $ \\ \hline
$W_{8}^{\alpha \in \mathbb{F}^{\ast }}=\left\langle \left[ \delta _{a,c}%
\right] +\left[ \delta _{b,c}\right] ,\left[ \delta _{a,a}\right] +\alpha %
\left[ \delta _{b,b}\right] \right\rangle $ & $\left( 1,1\right) $ & $\left(
2,2,1\right) $ \\ \hline
\end{tabular}%
$%
\caption{The list of representatives.}%
\label{tab3}%
\end{table}%
Let $S_{1}=\left\{ W_{1},W_{2},W_{8}^{\alpha =0}\right\} $, $S_{2}=\left\{
W_{3},W_{5},W_{7}\right\} $ and $S_{3}=\left\{ W_{4},W_{8}^{\alpha \in 
\mathbb{F}^{\ast }}\right\} $. Then, from Table \ref{tab3}, we have $%
\mbox{Orb}\left( X\right) \cap $ $\mbox{Orb}\left( Y\right) =\emptyset $ for 
$X\in S_{i},Y\in S_{j}$, $i\neq j$. Further we claim that $\mbox{Orb}\left(
W_{4}\right) \cap $ $\mbox{Orb}\left( W_{8}^{\alpha \in \mathbb{F}^{\ast
}}\right) =\emptyset $ and $\mbox{Orb}\left( X\right) \cap $ $\mbox{Orb}%
\left( Y\right) =\emptyset $ for $X,Y\in S_{i}$, $i=1,2$. To see this,
consider any $\phi =\big(a_{ij}\big)\in $ $Aut\left( J_{3,3}\right) $. Then
we either have $a_{11}a_{22}\neq 0,a_{12}=a_{21}=0$ or $%
a_{11}=a_{22}=0,a_{12}a_{21}\neq 0$. Assume first that $a_{11}a_{22}\neq
0,a_{12}=a_{21}=0$. Then 
\begin{eqnarray*}
\phi W_{1} &=&\left\langle 2a_{11}a_{31}\left[ \delta _{a,a}\right]
+a_{11}^{2}a_{22}\left[ \delta _{a,c}\right] ,a_{11}^{2}\left[ \delta _{a,a}%
\right] \right\rangle , \\
\phi W_{2} &=&\left\langle 2a_{11}a_{31}\left[ \delta _{a,a}\right]
+a_{11}^{2}a_{22}\left[ \delta _{a,c}\right] ,a_{22}^{2}\left[ \delta _{b,b}%
\right] \right\rangle , \\
\phi W_{3} &=&\left\langle 2a_{11}a_{31}\left[ \delta _{a,a}\right]
+a_{11}^{2}a_{22}\left[ \delta _{a,c}\right] ,2a_{22}a_{32}\left[ \delta
_{b,b}\right] +a_{11}a_{22}^{2}\left[ \delta _{b,c}\right] \right\rangle , \\
\phi W_{4} &=&\left\langle 2a_{11}a_{31}\left[ \delta _{a,a}\right]
+a_{11}^{2}a_{22}\left[ \delta _{a,c}\right] ,a_{11}^{2}\left[ \delta _{a,a}%
\right] +a_{22}^{2}\left[ \delta _{b,b}\right] \right\rangle , \\
\phi W_{5} &=&\left\langle 2a_{11}a_{31}\left[ \delta _{a,a}\right]
+a_{11}^{2}a_{22}\left[ \delta _{a,c}\right] ,a_{11}^{2}\left[ \delta _{a,a}%
\right] +2a_{22}a_{32}\left[ \delta _{b,b}\right] +a_{11}a_{22}^{2}\left[
\delta _{b,c}\right] \right\rangle .
\end{eqnarray*}%
As $\left[ \delta _{a,a}\right] +\left[ \delta _{b,c}\right] \notin $ $\phi
W_{3}$ we have $\phi W_{3}\neq W_{i}$ for $i=5,7$. Moreover, $\left[ \delta
_{b,b}\right] +\left[ \delta _{a,c}\right] \notin \phi W_{5}$ and therefore $%
\phi W_{5}\neq W_{7}$. Also, $\phi W_{1}\neq W_{2}$, $\phi W_{1}\neq
W_{8}^{\alpha =0}$ and $\phi W_{2}\neq W_{8}^{\alpha =0}$ since $\delta
_{b,b},\delta _{a,c}+\delta _{b,c}\notin $ $\phi W_{1}$ and $\delta
_{a,c}+\delta _{b,c}\notin \phi W_{2}$. Furthermore, $\phi W_{4}\neq
W_{8}^{\alpha \in \mathbb{F}^{\ast }}$ since $\delta _{a,c}+\delta
_{b,c}\notin \phi W_{4}$.

Assume now that $a_{11}=a_{22}=0,a_{12}a_{21}\neq 0$. Then%
\begin{eqnarray*}
\phi W_{1} &=&\left\langle 2a_{12}a_{32}\left[ \delta _{b,b}\right]
+a_{12}^{2}a_{21}\left[ \delta _{b,c}\right] ,a_{12}^{2}\left[ \delta _{b,b}%
\right] \right\rangle , \\
\phi W_{2} &=&\left\langle 2a_{12}a_{32}\left[ \delta _{b,b}\right]
+a_{12}^{2}a_{21}\left[ \delta _{b,c}\right] ,a_{21}^{2}\left[ \delta _{a,a}%
\right] \right\rangle , \\
\phi W_{3} &=&\left\langle 2a_{12}a_{32}\left[ \delta _{b,b}\right]
+a_{12}^{2}a_{21}\left[ \delta _{b,c}\right] ,2a_{21}a_{31}\left[ \delta
_{a,a}\right] +a_{12}a_{21}^{2}\left[ \delta _{a,c}\right] \right\rangle , \\
\phi W_{4} &=&\left\langle 2a_{12}a_{32}\left[ \delta _{b,b}\right]
+a_{12}^{2}a_{21}\left[ \delta _{b,c}\right] ,a_{21}^{2}\left[ \delta _{a,a}%
\right] +a_{12}^{2}\left[ \delta _{b,b}\right] \right\rangle , \\
\phi W_{5} &=&\left\langle 2a_{12}a_{32}\left[ \delta _{b,b}\right]
+a_{12}^{2}a_{21}\left[ \delta _{b,c}\right] ,2a_{21}a_{31}\left[ \delta
_{a,a}\right] +a_{12}^{2}\left[ \delta _{b,b}\right] +a_{12}a_{21}^{2}\left[
\delta _{a,c}\right] \right\rangle .
\end{eqnarray*}%
As $\left[ \delta _{a,a}\right] +\left[ \delta _{b,c}\right] \notin $ $\phi
W_{3}$ we have $\phi W_{3}\neq W_{i}$ for $i=5,7$. Moreover, $\left[ \delta
_{b,b}\right] +\left[ \delta _{a,c}\right] \notin \phi W_{5}$ and therefore $%
\phi W_{5}\neq W_{7}$. Also, $\phi W_{1}\neq W_{2}$, $\phi W_{1}\neq
W_{8}^{\alpha =0}$ and $\phi W_{2}\neq W_{8}^{\alpha =0}$ since $\delta
_{a,c},\delta _{a,c}+\delta _{b,c}\notin $ $\phi W_{1}$ and $\delta
_{a,c}+\delta _{b,c}\notin \phi W_{2}$. Furthermore, $\phi W_{4}\neq
W_{8}^{\alpha \in \mathbb{F}^{\ast }}$ since $\delta _{a,c}+\delta
_{b,c}\notin \phi W_{4}$. Therefore, our claim follows.

It remains to determine the possible orbits among the representatives $%
W_{8}^{\alpha \in \mathbb{F}^{\ast }}$. We claim that $\mbox{Orb}\left(
W_{8}^{\alpha }\right) =\mbox{Orb}\left( W_{8}^{\beta }\right) $ if and only
if $\alpha =\beta $ or $\alpha =\beta ^{-1}$. To prove this, suppose first
that $\phi W_{8}^{\alpha }=W_{8}^{\beta }$ for some $\phi =\big(a_{ij}\big)%
\in $ $Aut\left( J_{3,3}\right) $. If $a_{12}=a_{21}=0$, then $\phi
W_{8}^{\alpha }=W_{8}^{\beta }$ yields to%
\begin{eqnarray*}
2a_{11}a_{31}\left[ \delta _{a,a}\right] +2a_{22}a_{32}\left[ \delta _{b,b}%
\right] +a_{11}^{2}a_{22}\left[ \delta _{a,c}\right] +a_{11}a_{22}^{2}\left[
\delta _{b,c}\right] &=&b_{11}\left( \left[ \delta _{a,c}\right] +\left[
\delta _{b,c}\right] \right) +b_{21}\left( \left[ \delta _{a,a}\right]
+\beta \left[ \delta _{b,b}\right] \right) , \\
a_{11}^{2}\left[ \delta _{a,a}\right] +\alpha a_{22}^{2}\left[ \delta _{b,b}%
\right] &=&b_{12}\left( \left[ \delta _{a,c}\right] +\left[ \delta _{b,c}%
\right] \right) +b_{22}\left( \left[ \delta _{a,a}\right] +\beta \left[
\delta _{b,b}\right] \right) .
\end{eqnarray*}%
with $\det \big(b_{ij}\big)\neq 0$. The last equation implies $\alpha
a_{22}^{2}=a_{11}^{2}\beta $, while the first implies $a_{11}=a_{22}$. This
ensures that $\alpha =\beta $. Further if $a_{11}=a_{22}=0$, then $\phi
W_{8}^{\alpha }=W_{8}^{\beta }$ yields to%
\begin{eqnarray*}
2a_{21}a_{31}\left[ \delta _{a,a}\right] +2a_{12}a_{32}\left[ \delta _{b,b}%
\right] +a_{12}a_{21}^{2}\left[ \delta _{a,c}\right] +a_{12}^{2}a_{21}\left[
\delta _{b,c}\right] &=&b_{11}\left( \left[ \delta _{a,c}\right] +\left[
\delta _{b,c}\right] \right) +b_{21}\left( \left[ \delta _{a,a}\right]
+\beta \left[ \delta _{b,b}\right] \right) , \\
\alpha a_{21}^{2}\left[ \delta _{a,a}\right] +a_{12}^{2}\left[ \delta _{b,b}%
\right] &=&b_{12}\left( \left[ \delta _{a,c}\right] +\left[ \delta _{b,c}%
\right] \right) +b_{22}\left( \left[ \delta _{a,a}\right] +\beta \left[
\delta _{b,b}\right] \right) ,
\end{eqnarray*}%
with $\det \big(b_{ij}\big)\neq 0$. The last equation implies $\alpha
a_{21}^{2}=\beta ^{-1}a_{12}^{2}$, while the first implies $a_{12}=a_{21}$.
This ensures that $\alpha =\beta ^{-1}$. Conversely, suppose that $\alpha
=\beta ^{-1}$. Let $\phi $ be the following automorphism:%
\begin{equation*}
\phi =%
\begin{bmatrix}
0 & 1 & 0 \\ 
1 & 0 & 0 \\ 
0 & 0 & 1%
\end{bmatrix}%
.
\end{equation*}%
Then $\phi W_{8}^{\beta ^{-1}}=W_{8}^{\beta }$. This finishes the proof of
the claim. So we get the following algebras:

\begin{itemize}
\item $J_{5,34}:a\circ b=c,a\circ c=d,a^{2}=e.$

\item $J_{5,35}:a\circ b=c,a\circ c=d,b^{2}=e.$

\item $J_{5,36}:a\circ b=c,a\circ c=d,b\circ c=e.$

\item $J_{5,37}:a\circ b=c,a\circ c=d,a^{2}=e,b^{2}=e.$

\item $J_{5,38}:a\circ b=c,a\circ c=d,a^{2}=e,b\circ c=e.$

\item $J_{5,39}:a\circ b=c,b^{2}=d,a\circ c=d,a^{2}=e.$

\item $J_{5,40}:a\circ b=c,b^{2}=d,a\circ c=d,a^{2}=e,b\circ c=e.$

\item $J_{5,41}^{\alpha \in \mathbb{F}}:a\circ b=c,a\circ c=d,b\circ
c=d,a^{2}=e,b^{2}=\alpha e$. Isomorphism: $J_{5,41}^{\alpha }\cong
J_{5,41}^{\beta }$ if and only if $\beta \left( \alpha ^{2}+1\right) =\alpha
\left( \beta ^{2}+1\right) $. To describe the isomorphism when $\beta
=\alpha ^{-1}$, let $\varphi (a)=b,\varphi (b)=a,\varphi (c)=c,\varphi
(d)=d,\varphi (e)=\alpha ^{-1}e$. Then $\varphi :J_{5,41}^{\alpha
}\longrightarrow J_{5,41}^{\alpha ^{-1}}$ is an isomorphism.
\end{itemize}

\begin{thm}
\label{Major}There is an infinite number of five-dimensional nilpotent
Jordan algebras over an algebraically closed fields of characteristic $\neq
2 $, and any nilpotent non-associative Jordan algebra of dimension five over
an algebraically closed field $\mathbb{F}$ of characteristic $\neq 2$ is
isomorphic to one of the following algebras:

\begin{itemize}
\item $J_{5,1}:a^{2}=b,b\circ c=d.$

\item $J_{5,2}:a\circ b=c,a\circ c=d.$

\item $J_{5,3}:a\circ b=c,a\circ c=d,b^{2}=d.$

\item $J_{5,4}:a\circ b=c,a\circ c=d,b\circ c=d.$

\item $J_{5,5}:a^{2}=b,d^{2}=e,b\circ c=e.$

\item $J_{5,6}:a^{2}=b,a\circ d=e,b\circ c=e.$

\item $J_{5,7}:a\circ b=c,c\circ d=e.$

\item $J_{5,8}:a\circ b=c,c\circ d=e,a^{2}=e.$

\item $J_{5,9}:a\circ b=c,c\circ d=e,a^{2}=e,b^{2}=e.$

\item $J_{5,10}:a\circ b=c,a\circ c=e,d^{2}=e,b\circ c=e.$

\item $J_{5,11}:a\circ b=c,a\circ c=e,d^{2}=e.$

\item $J_{5,12}:a\circ b=c,a\circ c=e,d^{2}=e,b^{2}=e.$

\item $J_{5,13}:a\circ b=c,a\circ c=e,a\circ d=e,b\circ c=e.$

\item $J_{5,14}:a\circ b=c,a\circ c=e,b\circ d=e.$

\item $J_{5,15}^{\alpha \in \mathbb{F}}:a^{2}=b,a\circ b=c,a\circ
c=e,b^{2}=e,b\circ d=e,d^{2}=\alpha e.$

\item $J_{5,16}:a\circ b=d,c^{2}=d,a\circ d=e.$

\item $J_{5,17}:a\circ b=d,c^{2}=d,a\circ d=e,b\circ c=e.$

\item $J_{5,18}:a\circ b=d,c^{2}=d,a\circ d=e,b^{2}=e.$

\item $J_{5,19}:a\circ b=d,c^{2}=d,c\circ d=e.$

\item $J_{5,20}:a\circ b=d,c^{2}=d,c\circ d=e,a^{2}=e.$

\item $J_{5,21}:a\circ b=d,c^{2}=d,c\circ d=e,a^{2}=e,b^{2}=e.$

\item $J_{5,22}:a^{2}=b,a\circ b=d,c^{2}=d,b^{2}=e,a\circ d=e.$

\item $J_{5,23}^{\alpha \in \mathbb{F}}:a^{2}=c,a\circ b=d,a\circ c=e,b\circ
d=e,b\circ c=\alpha e.$

\item $J_{5,24}^{\alpha \in \mathbb{F}^{\ast }-\left\{ 1\right\}
}:a^{2}=c,a\circ b=d,b\circ c=e,a\circ d=\alpha e.$

\item $J_{5,25}:a^{2}=c,a\circ b=d,a\circ c=e,b\circ c=-2e,a\circ
d=e.\allowbreak $

\item $J_{5,26}^{\alpha \in \mathbb{F}^{\ast }}:a^{2}=c,b^{2}=d,a\circ
c=e,a\circ d=e,b\circ d=\alpha e.$

\item $J_{5,27}^{\alpha ,\beta \in \mathbb{F}}\left( \alpha \beta \neq
1\right) :a^{2}=c,b^{2}=d,b\circ c=e,a\circ d=e,b\circ d=\alpha e,a\circ
c=\beta e.$

\item $J_{5,28}:a^{2}=b,b\circ c=d,c^{2}=e.$

\item $J_{5,29}:a^{2}=b,b\circ c=d,a\circ c=e.$

\item $J_{5,30}:a^{2}=b,b\circ c=d,a\circ c=e,c^{2}=e.$

\item $J_{5,31}:a^{2}=b,b\circ c=d,a\circ b=e.$

\item $J_{5,32}:a^{2}=b,b\circ c=d,a\circ b=e,a\circ c=e.$

\item $J_{5,33}:a^{2}=b,b\circ c=d,a\circ b=e,c^{2}=e.$

\item $J_{5,34}:a\circ b=c,a\circ c=d,a^{2}=e.$

\item $J_{5,35}:a\circ b=c,a\circ c=d,b^{2}=e.$

\item $J_{5,36}:a\circ b=c,a\circ c=d,b\circ c=e.$

\item $J_{5,37}:a\circ b=c,a\circ c=d,a^{2}=e,b^{2}=e.$

\item $J_{5,38}:a\circ b=c,a\circ c=d,a^{2}=e,b\circ c=e.$

\item $J_{5,39}:a\circ b=c,b^{2}=d,a\circ c=d,a^{2}=e.$

\item $J_{5,40}:a\circ b=c,b^{2}=d,a\circ c=d,a^{2}=e,b\circ c=e.$

\item $J_{5,41}^{\alpha \in \mathbb{F}}:a\circ b=c,a\circ c=d,b\circ
c=d,a^{2}=e,b^{2}=\alpha e.$

\item $M_{5,1}:a^{2}=b,a\circ b=d+e,c^{2}=d,b^{2}=e,a\circ d=e.$

\item $M_{5,2}:a^{2}=b,b\circ c=d,a\circ b=e,a\circ c=e,c^{2}=e.$
\end{itemize}

Among these algebras there are precisely the following isomorphisms:

\begin{itemize}
\item $J_{5,23}^{\alpha }\cong J_{5,23}^{\beta }$ if and only if $\alpha
^{2}-\beta ^{2}=0.$

\item $J_{5,26}\cong J_{5,26}^{\beta }$ if and only if $\alpha ^{2}-\beta
^{2}=0.$

\item $J_{5,27}^{\alpha ,\beta }\cong J_{5,27}^{\alpha ^{\prime },\beta
^{\prime }}$ if and only if $\left( \alpha ,\beta \right) =\left( \alpha
^{\prime },\beta ^{\prime }\right) $ or $\left( \alpha ,\beta \right)
=\left( \beta ^{\prime },\alpha ^{\prime }\right) .$

\item $J_{5,41}^{\alpha }\cong J_{5,41}^{\beta }$ if and only if $\beta
\left( \alpha ^{2}+1\right) =\alpha \left( \beta ^{2}+1\right) .$

\item $M_{5,1}\cong J_{5,22}$ if and only if characteristic $\mathbb{F}$ $%
\neq 3.$

\item $M_{5,2}\cong J_{5,33}$ if and only if characteristic $\mathbb{F}$ $%
\neq 3.$
\end{itemize}
\end{thm}

\end{document}